\documentclass[11pt,reqno]{amsart}

\usepackage{amssymb,amsmath,amsthm,amsfonts}
\usepackage{hyperref, multirow, bm, color, marginnote, graphicx, wrapfig, subcaption,esint}
\usepackage{diagbox}

\makeatletter
\renewcommand{\tocsection}[3]{%
  \indentlabel{\@ifnotempty{#2}{\bfseries\ignorespaces#1 #2\quad}}\bfseries#3}
\renewcommand{\tocsubsection}[3]{%
  \indentlabel{\@ifnotempty{#2}{\ignorespaces#1 #2\quad}}#3}

\newcommand\@dotsep{4.5}
\def\@tocline#1#2#3#4#5#6#7{\relax
  \ifnum #1>\c@tocdepth 
  \else
    \par \addpenalty\@secpenalty\addvspace{#2}%
    \begingroup \hyphenpenalty\@M
    \@ifempty{#4}{%
      \@tempdima\csname r@tocindent\number#1\endcsname\relax
    }{%
      \@tempdima#4\relax
    }%
    \parindent\z@ \leftskip#3\relax \advance\leftskip\@tempdima\relax
    \rightskip\@pnumwidth plus1em \parfillskip-\@pnumwidth
    #5\leavevmode\hskip-\@tempdima{#6}\nobreak
    \leaders\hbox{$\m@th\mkern \@dotsep mu\hbox{.}\mkern \@dotsep mu$}\hfill
    \nobreak
    \hbox to\@pnumwidth{\@tocpagenum{\ifnum#1=1\bfseries\fi#7}}\par
    \nobreak
    \endgroup
  \fi}
\AtBeginDocument{%
\expandafter\renewcommand\csname r@tocindent0\endcsname{0pt}
}
\def\l@subsection{\@tocline{2}{0pt}{2.5pc}{5pc}{}}
\makeatother

\usepackage[margin=1in]{geometry}

\newcommand{\R}{\mathbb{R}}
\newcommand{\Z}{\mathbb{Z}}
\newcommand{\N}{\mathbb{N}}

\newcommand{\A}{\mathcal{A}}

\newcommand{\bR}{\bm{R}}

\newcommand{\bars}{\overline s}

\newcommand{\bu}{\bm{u}}

\newcommand{\bx}{\bm{x}}
\newcommand{\by}{\bm{y}}

\newcommand{\X}{\bm{X}}
\newcommand{\Y}{\bm{Y}}
\newcommand{\be}{\bm{e}}

\newcommand{\p}{\partial}
\renewcommand{\div}{{\rm{div}\,}}

\newcommand{\abs}[1]{\left\lvert #1 \right\rvert}
\newcommand{\norm}[1]{\left\lVert #1 \right\rVert}
\newcommand{\wh}[1]{\widehat{#1}}
\newcommand{\wt}[1]{\widetilde{#1}}
\newcommand{\mc}[1]{\mathcal{#1}}

\newtheorem{theorem}{Theorem}[section]
\newtheorem{lemma}[theorem]{Lemma}

\theoremstyle{definition}
\newtheorem{remark}{Remark}

\allowdisplaybreaks

\numberwithin{equation}{section}

\begin{document}
\title{A hierarchy of blood vessel models, Part I: 3D-1D to 1D}

\author{Laurel Ohm}
\address{Department of Mathematics, University of Wisconsin - Madison, Madison, WI 53706}
\email{lohm2@wisc.edu}

\author{Sarah Strikwerda}
\address{Department of Mathematics, University of Wisconsin - Madison, Madison, WI 53706}
\email{sstrikwerda@wisc.edu}

\begin{abstract}
We propose and analyze a family of models describing blood perfusion through a tissue surrounding a thin blood vessel. Our goal is to rigorously establish convergence results among 3D-3D Darcy--Stokes, 3D-1D Darcy--Poiseuille, and 1D Green's function methods commonly used to model this process. 
In Part I, we propose a 3D-1D Darcy--Poiseuille system where the coupling across the permeable vessel surface involves an angle-averaged Neumann boundary condition coupled with a geometrically constrained Robin boundary condition. We show that this model is well-posed and moreover limits to a 1D Green's function model as the maximum vessel radius $\epsilon\to 0$. 
In the 1D model, the exterior blood pressure is given by an explicit Green's function expression involving the interior blood pressure. The interior pressure satisfies a novel 1D integrodifferential equation in which the integral term incorporates the effects of the exterior pressure and the vessel geometry.
Much of this paper is devoted to analyzing this integrodifferential equation. Using the \emph{a priori} bounds obtained here, we show that the solution to the 1D model converges to the 3D-1D solution with a rate proportional to $\epsilon^{1/2}\abs{\log\epsilon}$.
In Part II \cite{partII}, we rely on the 1D estimates to show that both the 1D and 3D-1D models converge to a coupled 3D-3D Darcy-Stokes system as $\epsilon\to 0$, thereby establishing a convergence chain among all hierarchy levels.
\end{abstract}

\maketitle

\tableofcontents

\section{Introduction}
We propose and analyze a family of models describing blood perfusion through a tissue surrounding a thin blood vessel. Perfusion, the process by which blood delivers oxygen and nutrients to tissues and organs by way of a complex vascular system, is vital to maintaining functionality of bodily processes. The geometry of microvascular networks plays an important role in this process \cite{vitullo2023sensitivity, notaro2016mixed, pries2008blood,qi2021control,qi2022microvascular}, including the location and shape of arterioles and venules, thin vessels which deliver and collect, respectively, blood from dense capillary beds within the tissue. 
To understand the effects of possibly complicated arteriole or venule geometries on perfusion, reduced-dimensional models which capture the essential physics are important, as the reduced dimensionality facilitates consideration of more complex shapes. Obtaining effective reduced equations is generally an important aspect of blood vessel modeling \cite{canic2003effective, tambavca2005effective, zunino2016integrated,possenti2021mesoscale,laurino2019derivation,secomb2004green,grinberg2011modeling,d2008coupling,vcanic2005two,boulakia2024mathematical,berntsson2018one,blanco2009potentialities,fritz20221d, blanco20103d,bociu2025analysis,formaggia2001coupling}.

In this two-part series, we analyze a hierarchy of modeling frameworks treating this coupled process of blood flow within the arteriole or venule and perfusion through the surrounding capillary network at three levels of complexity. Our goal is to rigorously establish convergence results among 3D-3D Darcy--Stokes, 3D-1D Darcy--Poiseuille, and 1D Green's function methods, at least for a specific choice of coupling. For all of the models, we consider a fixed geometric setup wherein the capillary network surrounding the vessel is treated as an isotropic porous medium occupying the upper half space $\R^3_+$. The vessel itself is taken to have one end embedded at the tissue boundary at $z=0$ and one end free in the medium, and may be highly curved with a necessarily non-constant radius that decays to zero (i.e. becomes indistinguishable from the surrounding capillary bed) at the free end (see figure \ref{fig:vessel}). The maximum vessel radius is assumed to be a small parameter $\epsilon>0$. We note that the free end is both a convenience and a challenge. We work in the upper half space because we rely on an explicit Green's function to formulate and study the 1D model. To allow the vessel centerline geometry to be as general as possible, we choose to consider a free end within the medium, where the vessel radius becomes so small that it is indistinguishable from the surrounding capillaries. The decaying radius poses analytical challenges, however, and our treatment of the endpoint limits the convergence.

We further remark that it is due to the presence of the free end that we choose to present these models at the level of an arteriole or venule rather than at the level of a capillary, even though the modeling setup is very similar. It is perhaps more common to model perfusion at the capillary level, with the interstitium serving as the surrounding porous medium \cite{berrone2023optimization,notaro2016mixed,possenti2021mesoscale,hyde2013parameterisation}. 
However, the general error analysis framework here can be readily extended to an infinite slab with the vessel connected to walls at $z=0$ and $z=L$ and a prescribed constant pressure at both ends (see Remark \ref{rem:remark}). The Green's function in this setting, while explicit, is less convenient for practical implementation. We also note that some animals, including arthropods and some mollusks, do have open circulatory systems where vessels drain directly into a cavity known as the hemocoel in which blood mixes with interstitial fluid \cite{molnar202221}.

Here in Part I, we propose a coupled 3D-1D Darcy--Poiseuille system for the blood pressure outside and inside the vessel, respectively, in which the coupling across the vessel boundary involves a novel angle-averaged Neumann boundary condition coupled with a geometrically constrained Robin boundary condition. This is similar to the boundary conditions arising in a PDE justification of slender body theories for impermeable filaments in Stokes flow \cite{closed_loop, free_ends,rigid,inverse,laplace}. Many variations of 3D-1D Darcy--Poiseuille systems with different couplings have been proposed and analyzed \cite{berrone2023optimization, koch2020modeling, d2007multiscale, d2012finite, d2008coupling, blanco2007unified, formaggia2007stability, koppl20203d, nobile2009coupling, formaggia2006coupling, blanco2009potentialities, notaro2016mixed, possenti2021mesoscale, laurino2019derivation, kuchta2021analysis, boulakia2024mathematical}, and have seen great success in exploring the effects of complicated vessel geometries (again, these are often used at the level of a capillary). One benefit of our particular choice of coupling is that we show that the 3D-1D system itself naturally limits to a 1D `slender-body'-type approximation. In the 1D model, the exterior pressure is given by an explicit Green's function expression, similar to the methods in \cite{secomb2004green, hsu1989green, pries2008blood, fry2018predicting, qi2021control}, where the interior pressure determines the potential. Here the interior pressure satisfies a novel 1D integrodifferential equation along the length of the vessel where the integral term incorporates the effects of both the exterior pressure as well as the vessel geometry. Due to the decaying radius at the vessel free end, the 1D model is degenerate, leading to additional analytical challenges. However, the model is still flexible enough to easily handle complicated vessel centerline geometries, as demonstrated by numerics in section \ref{subsec:numerics}.

This paper is devoted to a detailed analysis of the 1D model, which allows us to prove a rigorous convergence result between the 3D-1D system and 1D model at a rate proportional to $\epsilon^{1/2}\abs{\log\epsilon}$. Here the loss of $\epsilon^{-1/2}$ over, for example, the analogous bound for Stokes slender body theory \cite{closed_loop, free_ends,inverse} is due to properties of the integral operator in the degenerate 1D integrodifferential equation which limit the \emph{a priori} estimates we are able to obtain.
In Part II \cite{partII} of our program, we show that both the 3D-1D model and the 1D approximation may be derived from a coupled 3D Darcy--3D Stokes system for the flow outside and within the vessel, respectively. There we obtain a convergence result between the 3D-3D system and the 1D model at a rate proportional to $\epsilon^{1/6}\abs{\log\epsilon}$, where the additional loss in $\epsilon$ is due entirely to the effect of the free end. 
More detailed asymptotics for the tip boundary layer would likely yield an improvement here, at the expense of more complicated 3D-1D and 1D models. Here we opt for simplicity of the reduced models, noting that related considerations for outlet boundary conditions are typical of pipe flows \cite{ghosh2021modified,canic2003effective,castineira2019rigorous}. 
In any case, the \emph{a priori} estimates obtained here for the 1D model are essential for the convergence results of Part II.

A benefit of the 1D model proposed here is that it captures the essential behavior of the 3D-3D system in a quantifiable way, especially the effects of vessel geometry and the interplay between the interior and exterior of the vessel. At the same time, the model is simple enough to probe more complex questions. Our 1D model is very similar to the model proposed in \cite{qi2021control,qi2024hemodynamic,qi2022microvascular} and can be similarly generalized to include multiple interacting vessels for exploring the role of collective vascular architectures in perfusion. The 1D model may also be useful in probing questions of shape optimization or inverse problems related to vessel arrangement. 
Furthermore, the hierarchy of models proposed here may be used to rigorously justify effective 1D descriptions of branching, an important aspect of blood vessel geometry \cite{kozlov2017one,berrone2023optimization,possenti2021mesoscale,notaro2016mixed,maruvsic2003rigorous,vitullo2023sensitivity,grinberg2011modeling,shipley2020hybrid,cassot2010branching} that should be incorporated into questions of shape optimization.

In many ways, this modeling hierarchy is a starting point. At all levels, the models are missing important features such as the elasticity of the vessel itself. The framework introduced here could be combined with reduced order models which incorporate the effects of vessel wall elasticity, such as
\cite{grinberg2011modeling, coccarelli2021framework, berntsson2018one, panasenko2020three, zunino2016integrated, vcanic2006blood, canic2003effective, mikelic2007fluid, vcanic2005two, tambavca2005effective}.
Furthermore, we consider only static models -- a single snapshot in time. Incorporating a time-dependent incoming pressure, vessel radius, and vessel length would allow important considerations including vessel growth \cite{berrone2023optimization,pries2008blood} to be addressed.

\subsection{Vessel geometry}
We consider a dense network of capillaries, modeled as a porous medium, occupying the upper half space $\R^3_+$, with a boundary at $z=0$. Let $\Gamma_0$ denote a unit-length, $C^2$ curve in $\R^3_+$ with one end attached to the tissue boundary at $z=0$ and the other end free in the porous medium. Let $\X:[0,1]\to \R^3_+$ denote the arclength parameterization of $\Gamma_0$ with parameter $s$. For convenience, we will require that $\frac{d\X}{ds}\big|_{s=0}=\X_s(0)$ is perpendicular to the $z=0$ plane so that the reflection of $\X(s)$ across $z=0$ has the same regularity as $\X$. We will also require that 
\begin{equation}
\min\bigg\{\inf_{s_1\neq s_2} \frac{\abs{\X(s_1)-\X(s_2)}}{\abs{s_1-s_2}}\,, \quad \inf_{0<s\le 1}\frac{{\rm dist}(\X(s),\{z=0\})}{s}\bigg\} = c_\Gamma>0\,,
\end{equation}
i.e. the vessel centerline does not self-intersect and, away from the $s=0$ cross section, the vessel does not intersect the wall at $z=0$.

The curve $\Gamma_0$ is the centerline of a blood vessel with circular cross sections of radius $\epsilon a(s)$, $0<\epsilon\ll1$, which may vary along the length of the vessel. We define an admissible radius function $a: [0,1]\to [0,1]$ by the following criteria (see \cite{free_ends}):
\begin{enumerate}
  \item (Regularity) The radius function satisfies $a(s)\in C^2[0,1)$ with 
  \begin{equation}\label{eq:astar}
  a_\star=\sup_{s\in[0,1)}\abs{a(s)a'(s)} <\infty\,, \quad a_{\star\star}=\sup_{s\in[0,1)}\abs{a^3(s)a''(s)} <\infty\,.
  \end{equation}
  \item (Upper/lower bounds) We have $\norm{a}_{L^\infty}=1$, and there exists $0<\delta\ll1$ independent of $\epsilon$ such that 
  \begin{equation}\label{eq:delta}
  a(s)\ge a_0>0 \qquad \text{for }0\le s\le 1-\delta\,.
  \end{equation}
  \item (Spheroidal ends) Given $\delta$ as above, for $s>1-\delta$, the radius function satisfies
  \begin{equation}\label{eq:spheroidal}
  \abs{a(s)-\sqrt{1-s^2}} \le C\epsilon^2\sqrt{1-s^2}\,.
  \end{equation}
  In particular, $a(1)=0$, and we will require the decay to be monotonic as $s\to 1$.
\end{enumerate}

As a prototypical example, we may consider $a$ to be the upper hemisphere 
\begin{equation}
a(s) = \sqrt{1-s^2}\,, \quad s\in[0,1]\,,
\end{equation}
so that $\epsilon a(s)$ is the equatorial radius of the upper half of a slender prolate spheroid.

\begin{figure}[!ht]
\centering
\includegraphics[scale=0.4]{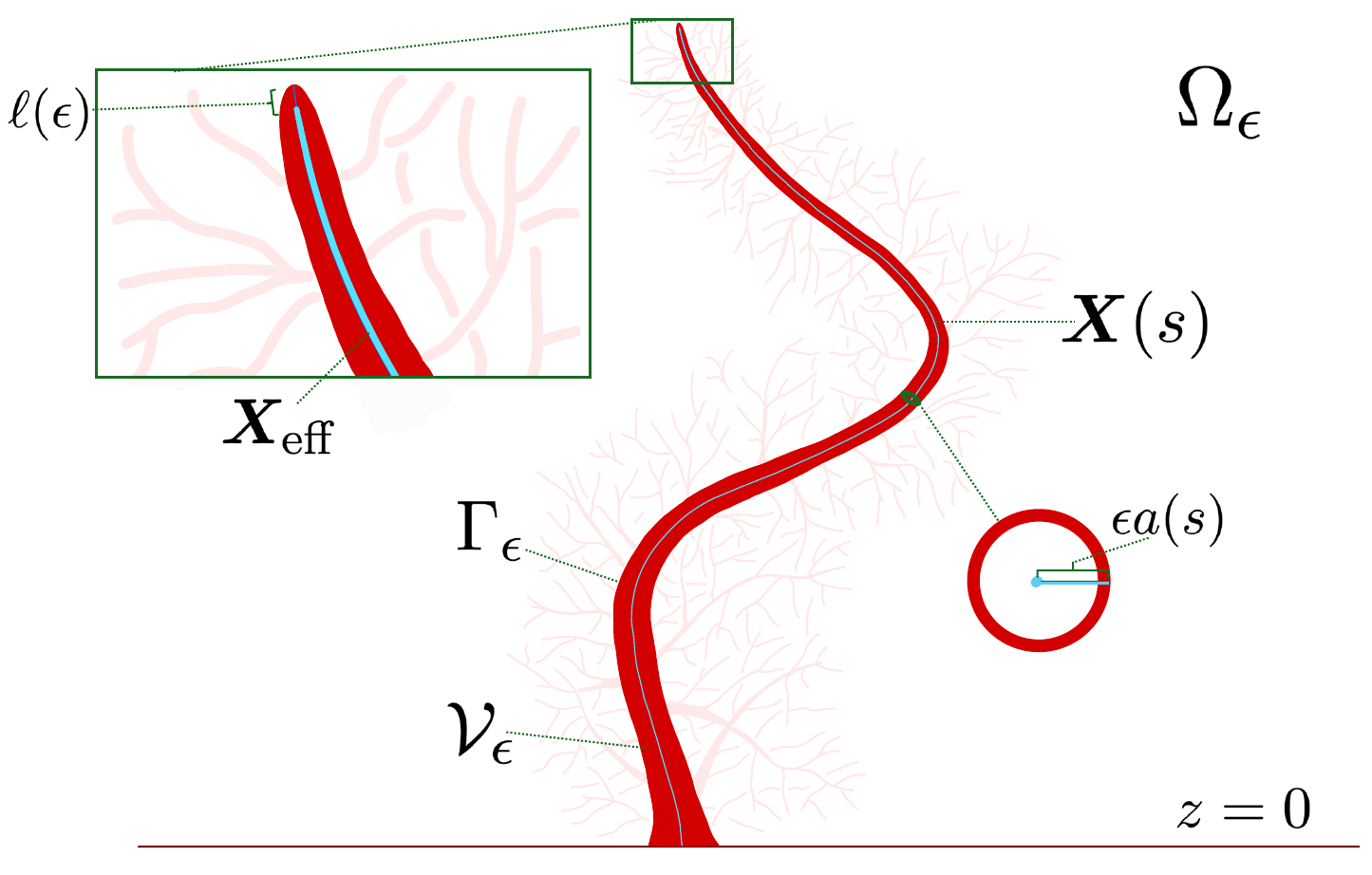}
\caption{An example of the blood vessel geometry $\mc{V}_\epsilon$ considered in this analysis. }
\label{fig:vessel}
\end{figure}

For each $s$, we may consider the plane $\X_s^\perp(s)$ perpendicular to $\X(s)$ and define the blood vessel $\mc{V}_\epsilon$ and its surface $\Gamma_\epsilon$ as
\begin{align}
\mc{V}_\epsilon &= \{\bx\in \X_s^\perp(s) \,:\, {\rm dist}(\bx,\X(s))<\epsilon a(s)\,, \; 0\le s\le 1  \} \,, \label{eq:Veps} \\
\Gamma_\epsilon &= \{\bx\in \X_s^\perp(s) \,:\, {\rm dist}(\bx,\X(s))=\epsilon a(s)\,, \; 0\le s\le 1  \}\,. 
\end{align}
Here we assume that $\epsilon$ is sufficiently small with respect to the vessel curvature that each $\bx\in\mc{V}_\epsilon$ belongs to a unique cross section, and, besides the $s=0$ cross section, $\mc{V}_\epsilon\cap\{z=0\}=\varnothing$.  
An example of the blood vessel geometry is pictured in figure \ref{fig:vessel}.

Throughout, we will think of the surface $\Gamma_\epsilon$ as a function of arclength $s$ and angle $\theta$ over each cross section. Along $\Gamma_\epsilon$, we then denote the surface element $dS_\epsilon$ as 
\begin{equation}\label{eq:dS}
dS_\epsilon = \mc{J}_\epsilon(s,\theta)\,d\theta ds\,,
\end{equation}
where $\mc{J}_\epsilon$ is a Jacobian factor that we will parameterize explicitly when needed.

\subsection{The 3D-1D Darcy--Poiseuille model and the 1D slender body model}
We begin by introducing the coupled 3D-1D Darcy--Poiseuille system. Let $\Omega_\epsilon=\R^3_+\backslash\overline{\mc{V}_\epsilon}$ denote the porous medium surrounding the vessel.
Throughout $\Omega_\epsilon$, the blood flow is given by Darcy's law coupled with an incompressibility constraint; in particular the blood pressure $q$ satisfies the Laplace equation \eqref{eq:darcy}. This must be coupled with an equation for the pressure $p$ within the thin blood vessel, which we will assume to be a function of arclength only, $p=p(s)$; in other words, the interior pressure is constant across each cross section of the vessel. This assumption arises from Poiseuille flow asymptotics considered in Part II \cite{partII}. One tricky aspect of formulating a boundary value problem for $q$ and $p$ is giving meaning to the coupling between the 3D exterior pressure field and the 1D interior pressure in a well-defined way. Other methods have been proposed in \cite{d2008coupling, nobile2009coupling,
formaggia2006coupling, koppl20203d, formaggia2007stability, blanco2007unified, d2012finite, laurino2019derivation, notaro2016mixed, kuchta2021analysis}, but we emphasize that the formulation considered here lends itself particularly well to a fully 1D reduction. We consider $(q,p):\Omega_\epsilon\times[0,1]\to\R\times\R$ satisfying
\begin{subequations}
\begin{align}
\Delta q &= 0 \hspace{2.8cm} \text{in }\Omega_\epsilon  \label{eq:darcy} \\
\frac{\p q}{\p \bm{n}} &= \frac{\omega}{\epsilon} (p(s)-q) \qquad \hspace{0.25cm} \text{on }\Gamma_\epsilon  \label{eq:robin}\\
\int_0^{2\pi}\frac{\p q}{\p \bm{n}}\,\mc{J}_\epsilon(s,\theta)\,d\theta &= \eta\frac{d}{ds}\bigg(a^4(s)\frac{dp}{ds}\bigg) \hspace{0.3cm} \text{ on }\Gamma_\epsilon  \label{eq:RtoN} \\
\frac{\p q}{\p\bm{n}} &= 0 \; \text{ on }\p\Omega_\epsilon\backslash\Gamma_\epsilon\,, \quad q\to 0 \text{ as }\abs{\bx}\to\infty \label{eq:extBCs}
\end{align}
\end{subequations}
with incoming pressure data $p(0)=p_0\in \R$ on the tissue boundary. We understand $(q,p)$ as being measured with respect to a constant background pressure level, so that $p_0<0$ corresponds to a venule (outgoing) and $p_0>0$ corresponds to an arteriole (incoming). Here $\bm{n}$ points out of the domain $\Omega_\epsilon$, into the vessel $\mc{V}_\epsilon$.
The boundary condition \eqref{eq:robin} states that the pressure flux through the boundary of the blood vessel $\frac{\p q}{\p \bm{n}}\big|_{\Gamma_\epsilon}$ is proportional to the pressure difference $p-q$ inside and outside the vessel, respectively, where the interior pressure $p=p(s)$ is unknown but, as mentioned, constant across each cross section of the vessel. 
The right hand side of the boundary condition \eqref{eq:RtoN} is the variable-radius Poiseuille's law for the pressure $p(s)$ within $\mc{V}_\epsilon$, which is balanced by the total flux of $q$ across each cross section. Due to the degeneracy of \eqref{eq:RtoN} as $s\to 1$, $p$ does not have a well-defined boundary value in the tissue at $s=1$. The parameter $\eta$ corresponds to the relative conductivity of the interior of the vessel compared to the exterior, with larger $\eta$ describing a more permeable interior. The parameter $\omega$ in the pressure balance boundary condition \eqref{eq:robin} corresponds to the conductivity of the vessel walls, with large $\omega$ serving to equalize the interior and exterior pressures across $\Gamma_\epsilon$.

Multiplying \eqref{eq:robin} by the Jacobian factor $\mc{J}_\epsilon$ and integrating in $\theta$, the boundary condition \eqref{eq:RtoN} may be seen to be equivalent to 
\begin{equation}\label{eq:RtoD}
-\eta\frac{d}{ds}\bigg(a^4(s)\frac{dp}{ds}\bigg)+ j_\epsilon(s) \,\frac{\omega}{\epsilon} \,p(s) = \frac{\omega}{\epsilon}\int_0^{2\pi} q(s,\theta)\,\mc{J}_\epsilon(s,\theta)\,d\theta   \qquad \text{on }\Gamma_\epsilon\,,  
\end{equation}
where $j_\epsilon(s) =\int_0^{2\pi}\mc{J}_\epsilon(s,\theta)d\theta\approx 2\pi \epsilon a(s)$. Both \eqref{eq:RtoN} and \eqref{eq:RtoD} will be useful for formulating the 1D model.

To construct the 1D model, we will first need to introduce some notation.  
Let $\mc{G}_N(\bx,\by)$ denote the Green's function for the Neumann problem for the Laplace equation in the upper half space $\R^3_+$, i.e. satisfying 
\begin{equation}\label{eq:GreenNeumann}
\begin{aligned}
-\Delta_x \mc{G}_N(\bx,\by) &= \delta(\bx-\by)\,, \quad \bx,\by\in \R^3_+\,,\\
\frac{\p \mc{G}_N}{\p\bm{n}_x} &= 0\,,  \hspace{1.6cm} \bx = (x_1,x_2,0)\,, 
\end{aligned}
\end{equation}
where here $\delta(\bx)$ denotes the Dirac delta. The expression for $\mc{G}_N$ may be derived by correcting the full space Green's function in $\R^3$ to satisfy the Neumann boundary condition at $z=0$ (see, e.g., \cite{craig2018course}). It is given by
\begin{equation}\label{eq:GN}
\mc{G}_N(\bx,\by) = \frac{1}{4\pi}\bigg(\frac{1}{\abs{\bx-\by}}+\frac{1}{\abs{\bx-\by^*}} \bigg)\,,
\end{equation}
where, for $\by=(y_1,y_2,y_3)$, the point $\by^*$ is given by reflection across $y_3=0$: $\by^* = (y_1,y_2,-y_3)$.

To build the 1D approximation to \eqref{eq:darcy}-\eqref{eq:extBCs}, we will need some additional definitions. As in \cite{free_ends}, for some $\ell(\epsilon)$ with $\frac{1}{C}\epsilon^2\le \ell(\epsilon)\le C\epsilon^2$, we define the \emph{effective centerline} of the vessel to be 
\begin{equation}\label{eq:Xeff}
\X_{\rm eff} = \{ \X(s) \,:\, 0\le s\le 1-\ell(\epsilon) \}\,,
\end{equation}
i.e. extending up to $O(\epsilon^2)$ away from the free endpoint of the vessel at $\X(1)$ (see figure \ref{fig:vessel}). For the prolate hemispheroid $\epsilon a(s)=\epsilon\sqrt{1-s^2}$, for example, $\ell(\epsilon)$ may be chosen such that $\X_{\rm eff}$ extends up to the focus of the spheroid at $s=\sqrt{1-\epsilon^2}$. 
Along with the effective centerline, we define a stretch operator
\begin{equation}\label{eq:varphi}
\varphi_\epsilon(s) = \frac{s}{1-\ell(\epsilon)}
\end{equation}
mapping $s\in [0,1-\ell(\epsilon)]$ to $[0,1]$. 
Given $f:[0,1]\to \R$, for $\bx\in \R^3_+$, we define the operator
\begin{equation}\label{eq:SN_op}
\mc{S}_N[f](\bx) = \int_0^1\mc{G}_N\big(\bx,\X(\varphi_\epsilon^{-1}(t))\big)\,f(t)\,dt\,.
\end{equation}
Here the inverse stretch operator $\varphi_\epsilon^{-1}$ appears in the second argument of $\mc{G}_N$ to ensure that $\mc{S}_N$ is well defined everywhere along the blood vessel surface $\Gamma_\epsilon$, up to and including the endpoint at $\X(1)$. We may then define an approximation $q^{\rm SB}(\bx)$ to the exterior pressure $q(\bx)$ by
\begin{equation}\label{eq:qSB}
q^{\rm SB}(\bx) = \eta\,\mc{S}_N\bigg[\frac{d}{ds}\bigg(a^4\frac{dp^{\rm SB}}{ds}\bigg)\bigg](\bx)\,, \qquad \bx\in \Omega_\epsilon\,,
\end{equation}
where $p^{\rm SB}(s)$ satisfies the integrodifferential equation
\begin{equation}\label{eq:pSB}
\begin{aligned}
-\eta\frac{d}{ds}\bigg(a^4\frac{dp^{\rm SB}}{ds}\bigg)+ 2\pi a(s)\omega \,p^{\rm SB} &= \omega\int_0^{2\pi} q^{\rm SB}\big|_{\Gamma_\epsilon}\, a(s)\,d\theta  \\
&= a(s)\omega \eta\int_0^{2\pi} \mc{S}_N\bigg[\frac{d}{ds}\bigg(a^4\frac{dp^{\rm SB}}{ds}\bigg)\bigg]\bigg|_{\Gamma_\epsilon}\,d\theta
\end{aligned}
\end{equation}
with $p^{\rm SB}(0)=p_0$. As in the 3D-1D system, $p^{\rm SB}(1)$ is not well-defined due to the degenerate $a^4$ coefficient.

We use the superscript `SB' because the expression \eqref{eq:qSB} for $q^{\rm SB}$ is analogous to classical slender body theories for fluid flows around elongated impermeable objects, from potential flow around a slender watercraft \cite{moran1963line, tuck1964some, handelsman1967axially,newman1970applications} to elastic filaments in Stokes flow \cite{keller1976slender,johnson1980improved,gotz2000interactions, tornberg2004simulating,lighthill1976flagellar,lauga2009hydrodynamics}.
One of the main novelties here is that rather than prescribing a known 1D density as data in $\mc{S}_N$, the density $(a^4p^{\rm SB}_s)_s$ comes from solving a closed-form 1D integrodifferential equation along the length of the vessel. 

The equation \eqref{eq:pSB} incorporates the effects of the vessel geometry, including its centerline shape and radius, via the integral operator $\mc{S}_N$ over the surface of the vessel. These geometric considerations, particularly the behavior of the operator toward the endpoint of the vessel, make this equation especially interesting from an analysis perspective. Note that the kernel of the operator $\mc{S}_N$ is smooth along $\Gamma_\epsilon$ by construction, and the integration in $\theta$ comes directly from the analogous expression \eqref{eq:RtoD} in the coupled 3D-1D model. In particular, the pressure $p^{\rm SB}$ along the vessel is well defined and avoids the logarithmic singularity that would arise from attempting to evaluate $\mc{S}_N$ directly on $\X(s)$.\footnote{For numerical methods, rather than evaluating $\mc{S}_N|_{\Gamma_\epsilon}$ and then integrating in $\theta$ to obtain a $\theta$-independent expression, one could perform an expansion of $\mc{S}_N|_{\Gamma_\epsilon}$ about $\epsilon=0$ to obtain a limiting centerline approximation in the style of Keller--Rubinow \cite{keller1976slender}; see also \cite{gotz2000interactions, tornberg2004simulating}. This may be more convenient computationally, but is much less convenient for analysis, as we rely heavily on favorable features of the $\theta$-integrated integral operator in \eqref{eq:pSB} that are not satisfied by such a limiting centerline approximation. } 

Since \eqref{eq:pSB} is a self-contained 1D equation, it is simpler to use in practice than the coupled 3D-1D system \eqref{eq:darcy}-\eqref{eq:extBCs}, which is in turn simpler than the full 3D-3D system \eqref{eq:3D3D}. It can handle complicated vessel geometries (see section \ref{subsec:numerics}) and can be generalized naturally to incorporate multiple blood vessels. Furthermore, the 1D model \eqref{eq:qSB}-\eqref{eq:pSB} can be used to formulate optimization and inverse problems related to vessel shape, and to consider time-dependent variants.

\subsection{Main results}
We first require well-posedness results and corresponding estimates for both the 3D-1D system (Theorem \ref{thm:3d1d}) and the 1D model (Theorem \ref{thm:pSB}). The estimates of Theorem \ref{thm:pSB} will be crucial not just here but also in Part II, where we use these results as a black box to compare the 1D model directly to a full 3D-3D Darcy--Stokes system. The main result of this paper is Theorem \ref{thm:3D1Dto1D}, an error estimate between the solution $(q,p)$ of the coupled 3D-1D Darcy--Poiseuille system and the solution $(q^{\rm SB},p^{\rm SB})$ of the 1D slender body model as $\epsilon\to 0$. 

We begin with a statement of well-posedness for the 3D-1D model \eqref{eq:darcy}-\eqref{eq:extBCs}, which is also explored in Part II \cite{partII}. Due to the degeneracy of the right hand side of \eqref{eq:RtoN} as $s\to 1$, we first define a type of weighted $H^1$ space and corresponding norm. We define 
\begin{equation}\label{eq:Ha_def}
\begin{aligned}
  \mc{H}^a(0,1) = \big\{ u\in L^2(0,1)\,:\, a^2u_s \in L^2(0,1) \big\}\,, 
  \quad
  \norm{u}_{\mc{H}^a(0,1)} := \norm{u}_{L^2(0,1)}+\norm{a^2u_s}_{L^2(0,1)}\,.
\end{aligned}
\end{equation}
We may then state the following.
\begin{theorem}[Well-posedness of 3D-1D model]\label{thm:3d1d}
Let $D^{1,2}(\Omega_\epsilon) = L^6(\Omega_\epsilon)\cap \dot H^1(\Omega_\epsilon)$. Given a curve $\X(s)$, a radius function $a(s)$ as in \eqref{eq:astar}-\eqref{eq:spheroidal}, and $p_0\in\R$, there exists a unique solution $(q,p)\in D^{1,2}(\Omega_\epsilon)\times \mc{H}^a(0,1)$ to the coupled 3D-1D system \eqref{eq:darcy}-\eqref{eq:extBCs} satisfying 
\begin{equation}\label{eq:energy}
\norm{q}_{D^{1,2}(\Omega_\epsilon)}+ \frac{1}{\abs{\Gamma_\epsilon}^{1/2}}\norm{p-q}_{\Gamma_\epsilon}+\norm{p}_{\mc{H}^a(0,1)} \le C\abs{p_0}
\end{equation}
for $C$ independent of $\epsilon$ as $\epsilon\to 0$.
\end{theorem}
The proof of Theorem \ref{thm:3d1d} appears in Appendix \ref{sec:app}.

The bulk of the analysis in this paper is in obtaining estimates for the novel 1D integrodifferential equation \eqref{eq:pSB} for $p^{\rm SB}(s)$. Much of the difficulty lies in the fact that different weights are needed as $s\to 1$ to account for behavior at the vessel endpoint. In addition, we need explicit knowledge of the $\epsilon$-dependence introduced by the integral operator $\mc{S}_N$ along $\Gamma_\epsilon$. 
\begin{theorem}[Estimates for the $p^{\rm SB}$ equation]\label{thm:pSB}
Given a curve $\X(s)$, a radius function $a(s)$ as in \eqref{eq:astar}-\eqref{eq:spheroidal}, pressure data $p_0\in\R$, and $\epsilon$ sufficiently small, there exists a unique solution $p^{\rm SB}(s)$ to the integrodifferential equation \eqref{eq:pSB} satisfying 
\begin{equation}\label{eq:pSB_basic_ests}
  \norm{p^{\rm SB}}_{\mc{H}^a(0,1)}\le C\abs{p_0}\,, \qquad \norm{p^{\rm SB}}_{L^\infty(0,1)} \le C\epsilon^{-1/2}\abs{p_0}\,, \qquad \norm{a\,p^{\rm SB}_s}_{L^\infty(0,1)} \le C\epsilon^{-1/2}\abs{p_0}\,.
\end{equation}
In addition, higher derivatives of $p^{\rm SB}$ satisfy the following weighted bounds:
\begin{equation}
\begin{aligned}
\|a^{-1/2}\big(a^4p^{\rm SB}_s\big)_s\|_{L^2(0,1)} &\le C\abs{p_0}\,, \qquad 
\norm{a^{-1}\big(a^4p^{\rm SB}_s\big)_s}_{L^\infty(0,1)} \le C\epsilon^{-1/2}\abs{p_0}\\
\|a^{3/2}(a^4p^{\rm SB}_s)_{ss}\|_{L^2(0,1)} &\le C\abs{\log\epsilon}\abs{p_0}\,, \qquad
\norm{a\,(a^4p^{\rm SB}_s)_{ss}}_{L^\infty(0,1)} \le C\epsilon^{-1/2}\abs{\log\epsilon}\abs{p_0}\,.
\end{aligned}
\end{equation}
Here the constants $C$ are bounded independent of $\epsilon$ as $\epsilon\to 0$.
\end{theorem}
The proof of Theorem \ref{thm:pSB} appears in section \ref{sec:pSB}. The main technical tools are a series of kernel estimates for the integral operator appearing in \eqref{eq:pSB}, which are shown in section \ref{subsec:kernel}. These in turn make use of a series of integral estimates from \cite{free_ends}, which can be largely reused modulo minor modification and are introduced in section \ref{sec:prelim}. The estimates for $p^{\rm SB}$ are especially tricky because we need to ensure that the $L^\infty$-based bounds do not grow too quickly as $\epsilon\to 0$ while balancing the decay/growth requirements at the vessel endpoint.

The $L^\infty$ bound for $p^{\rm SB}$ in Theorem \ref{thm:pSB} comes from a maximum principle satisfied by the left hand side terms of the $p^{\rm SB}$ equation \eqref{eq:pSB}, where the integral term is treated as forcing and is estimated using the $L^2$ bound for $(a^4p^{\rm SB}_s)_s$. Estimating the integral term in this way gives rise to the $\epsilon^{-1/2}$ loss. Without an endpoint, it may be possible to incorporate the integral itself into a maximum principle argument to obtain a uniform-in-$\epsilon$ $L^\infty$ bound for $p^{\rm SB}$. With the endpoint, we find numerical evidence that the full $p^{\rm SB}$ equation does not satisfy a maximum principle (see section \ref{subsec:numerics}). However, it may still be possible to improve upon the $\epsilon^{-1/2}$ bound for $\norm{p^{\rm SB}}_{L^\infty}$ by other means.

Given the bounds of Theorem \ref{thm:pSB} for $p^{\rm SB}(s)$, we obtain the following error estimate between the solution $(q,p)$ to the 3D-1D model \eqref{eq:darcy}-\eqref{eq:extBCs} and the 1D approximation $(q^{\rm SB},p^{\rm SB})$ as $\epsilon\to0$. 
\begin{theorem}[3D-1D to 1D error estimate]\label{thm:3D1Dto1D}
Given a curve $\X(s)$ and a radius function $a(s)$ as in \eqref{eq:astar}-\eqref{eq:spheroidal}, let $\mc{V}_\epsilon$ be a blood vessel as in \eqref{eq:Veps} and let $\Omega_\epsilon = \R^3_+\backslash\overline{\mc{V}_\epsilon}$. Given incoming pressure data $p_0\in\R$, let $(q,p)\in D^{1,2}(\Omega_\epsilon)\times \mc{H}^a(0,1)$ denote the solution to the coupled 3D-1D system \eqref{eq:darcy}-\eqref{eq:extBCs} in $\Omega_\epsilon\times [0,1]$, and let $(q^{\rm SB},p^{\rm SB})$ denote the 1D approximation \eqref{eq:qSB}, \eqref{eq:pSB}. The difference $(q^{\rm SB}-q,p^{\rm SB}-p)$ satisfies
\begin{equation}
\begin{aligned}
\norm{q^{\rm SB}-q}_{D^{1,2}(\Omega_\epsilon)} + &\,\epsilon^{-1/2}\norm{(p^{\rm SB}-p)-(q^{\rm SB}-q)}_{L^2(\Gamma_\epsilon)} + \norm{p^{\rm SB}-p}_{\mc{H}^a(0,1)} \\
&\le C \epsilon^{1/2}\abs{\log\epsilon}\,\abs{p_0}\,,
\end{aligned}
\end{equation}
where $C$ is independent of $\epsilon$ as $\epsilon\to 0$.
\end{theorem}

The proof of Theorem \ref{thm:3D1Dto1D} appears in section \ref{sec:error}, based on residual calculations for $q^{\rm SB}(\bx)$ in section \ref{sec:residual}.  

\begin{remark}\label{rem:remark}
As mentioned in the introduction, the convergence analysis here can be extended to the case of an infinite slab where the vessel connects to walls at $z=0$ and $z=L$, with constant pressure data supplied at both ends. This is perhaps more representative of the setting for human capillaries. In this case, the domain still admits an explicit Green's function, but instead of one reflection, an infinite sum of reflections is needed, effectively leading to an infinitely long vessel in the 1D model. 
  A technical consideration in this setting is the behavior of the exterior pressure as $\abs{\bx}\to\infty$. This may be addressed by requiring ``zero absorption'' at infinity, i.e. $\lim_{R\to\infty}\frac{\p q}{\p\bm{n}}\big|_{\p B_R}\to 0$, where $B_R$ is a ball of radius $R$ centered about the vessel. In particular, fluid may only enter and leave the tissue through the vessel. This amounts to the condition $\int_{\Gamma_\epsilon}(q-p)\,dS_\epsilon=0$, the net pressure difference across the vessel surface is zero. Analogous conditions are also needed in the 3D-3D and 1D models. 
\end{remark}

\subsection{Relation to 3D-3D Darcy--Stokes system}
In Part II of our program \cite{partII}, we consider a coupled 3D Darcy--3D Stokes system for vessel and surrounding capillary network. The blood flow inside the vessel is modeled by the 3D Stokes equations, which are coupled across $\Gamma_\epsilon$ to 3D Darcy flow outside. Letting $(\bu^\epsilon,p^\epsilon):\mc{V}_\epsilon\times \mc{V}_\epsilon \to \R^3\times \R$ denote the fluid velocity and pressure within the blood vessel, and letting $\Sigma^\epsilon = \mu(\nabla \bu^\epsilon+(\nabla\bu^\epsilon)^{\rm T})-p^\epsilon{\bf I}$ denote the corresponding stress tensor, the full 3D-3D Darcy--Stokes system is given by 
\begin{equation}\label{eq:3D3D}
\begin{aligned}
\Delta q^\epsilon &= 0 \hspace{5.75cm} \text{in }\Omega_\epsilon\\
-\mu\Delta\bu^\epsilon +\nabla p^\epsilon &= 0\,, \quad \div \bu^\epsilon=0 \hspace{3.45cm} \text{in } \mc{V}_\epsilon\\
\bu^\epsilon\cdot\bm{n} &= -\zeta\epsilon^4\frac{\p q^\epsilon}{\p\bm{n}}= k\epsilon^3\big((\Sigma^\epsilon\bm{n})\cdot\bm{n}+q^\epsilon\big) \qquad \text{on }\Gamma_\epsilon\\
\bu^\epsilon-(\bu^\epsilon\cdot\bm{n})\bm{n}&=0 \hspace{5.8cm} \text{on }\Gamma_\epsilon\\
\frac{\p q^\epsilon}{\p\bm{n}} &= 0 \; \text{ on }\p\Omega_\epsilon\backslash\Gamma_\epsilon\,, \quad q^\epsilon\to 0 \text{ as }\abs{\bx}\to\infty
\end{aligned}
\end{equation}
with $\Sigma^\epsilon\bm{n} = -p_0\bm{n}$ at the vessel base $\overline{\mc{V}_\epsilon}\cap\{z=0\}$. Note that here the pressure $p^\epsilon$ is not constant on cross sections. The parameters $\zeta$, $k$, and $\mu$ correspond to the 3D-1D model parameters via $\omega=\frac{k}{\zeta}$ and $\eta=\frac{\pi}{8\zeta\mu}$.

 In Part II, we rely on the explicit properties of the 1D model shown in Theorem \ref{thm:pSB} to show an error bound between the full 3D-3D system \eqref{eq:3D3D} and the 1D model. 
 In particular, using $p^{\rm SB}(s)$ and the curvilinear cylindrical coordinates $\bx=\bx(r,\theta,s)$ defined in \eqref{eq:param}, we may build the following velocity ansatz:
\begin{equation}\label{eq:Udef0}
\bm{U}(\bx) = -\frac{1}{16\mu}\p_s\big( \phi_\epsilon(s)(r^3 -2(\epsilon a)^2r) p^{\rm SB}_s \big)\be_r(s,\theta)+ \frac{1}{4\mu}\phi_\epsilon(s)\big(r^2-(\epsilon a)^2\big)p^{\rm SB}_s\,\be_{\rm t}(s)\,.
\end{equation} 
Here $\phi_\epsilon(s)$ is a smooth cutoff function in $s$ satisfying 
\begin{equation}\label{eq:phieps_def0}
  \phi_\epsilon(s) = \begin{cases}
    1\,, & 0\le s\le 1-2\epsilon^{4/3}\\
    0 \,, & 1-\epsilon^{4/3}\le s\le 1\,,
  \end{cases}
  \qquad \abs{\p_s\phi_\epsilon}\le C\epsilon^{-4/3}\,.
\end{equation}
This choice of ansatz and reasoning for the cutoff function at the tip is explored in detail in Part II. Using \eqref{eq:Udef0}, we show the following:
\begin{theorem}[3D-3D to 1D convergence \cite{partII}]\label{thm:PartII}
Given $\mc{V}_\epsilon$ as in \eqref{eq:Veps} with $\Gamma_0\in C^3$, pressure data $p_0\in \R$, and $\epsilon$ sufficiently small, let $(\bu^\epsilon,p^\epsilon,q^\epsilon)$ satisfy the 3D-3D system \eqref{eq:3D3D} and let $(\bm{U},p^{\rm SB},q^{\rm SB})$ satisfy the 1D model \eqref{eq:qSB}-\eqref{eq:pSB} with $\bm{U}$ given by \eqref{eq:Udef0}. The differences $(\bu^\epsilon-\bm{U},p^\epsilon-p^{\rm SB},q^\epsilon-q^{\rm SB})$ satisfy 
\begin{equation}
\begin{aligned}
&\frac{1}{|\mc{V}_\epsilon|^{1/2}}\bigg(\norm{\epsilon^{-2}(\bu^\epsilon-\bm{U})}_{L^2(\mc{V}_\epsilon)}+\norm{\epsilon^{-1}\nabla(\bu^\epsilon-\bm{U})}_{L^2(\mc{V}_\epsilon)} + \norm{p^\epsilon-p^{\rm SB}}_{L^2(\mc{V}_\epsilon)}\bigg)  \\
&\qquad  + \frac{1}{\epsilon^3|\Gamma_\epsilon|^{1/2}}\norm{(\bu^\epsilon-\bm{U})\cdot\bm{n}}_{L^2(\Gamma_\epsilon)} + \norm{q^\epsilon-q^{\rm SB}}_{D^{1,2}(\Omega_\epsilon)} \le C\epsilon^{1/6}\abs{\log\epsilon}\abs{p_0}
\end{aligned}
\end{equation}
for $C$ independent of $\epsilon$ as $\epsilon\to 0$.
\end{theorem}
The bound of Theorem \ref{thm:3D1Dto1D}, coupled with Theorem \ref{thm:PartII} from Part II, completes a chain of convergence among these three models describing perfusion outside of a thin blood vessel at three levels of detail. This makes precise the sense in which the simple model \eqref{eq:qSB}-\eqref{eq:pSB} approximates both the 3D-3D and 3D-1D models, and provides insight into the tradeoffs among implementing and analyzing each of the three descriptions. As emphasized in Part II, the largest error source in Theorem \ref{thm:PartII} is the free end, which results in the $\epsilon^{1/6}$ bound instead of the $\epsilon^{1/2}$ bound of Theorem \ref{thm:3D1Dto1D}. However, we emphasize that this bound holds up to the tip, despite the relative simplicity of the way in which we model the endpoint. More detailed tip asymptotics at the level of the 3D-3D model would likely yield an improvement in convergence here, similar to the outlet boundary layer solution constructed in \cite{canic2003effective}, at the expense of more complicated 3D-1D and 1D models at the tip.

\subsection{Numerics}\label{subsec:numerics}
As a brief demonstration of the utility of the 1D model \eqref{eq:qSB}-\eqref{eq:pSB} and an exploration of some of its features, we include some numerical examples. 

\begin{figure}[!ht]
\centering
\includegraphics[scale=0.5]{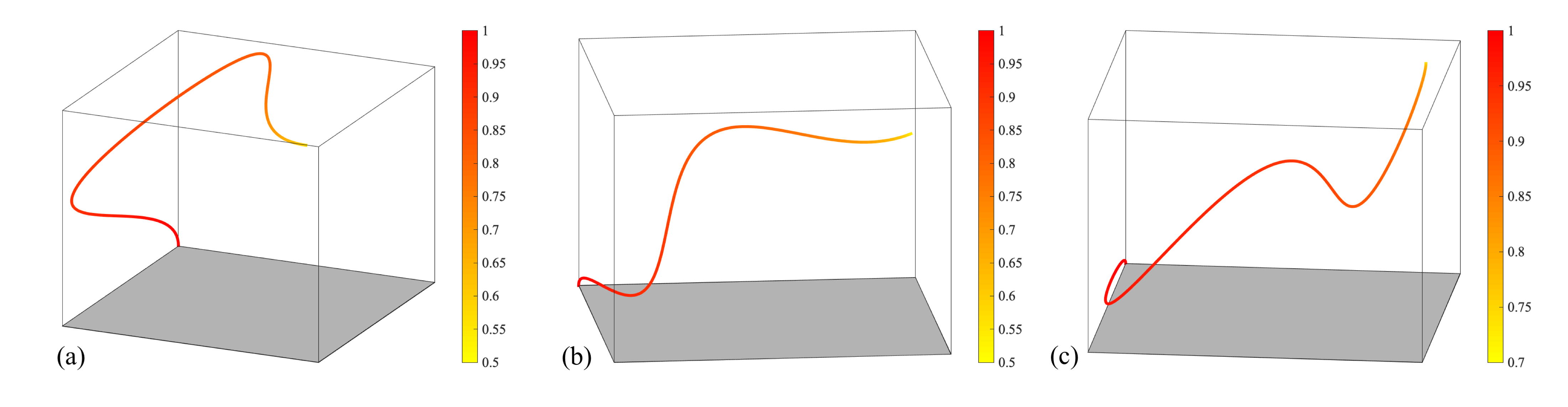}
\caption{Interior pressure $p^{\rm SB}(s)$ depicted along three different 3D vessel centerline geometries. In all cases, $p_0=1$, $\eta=0.05$, $\omega=10$, $\epsilon=0.01$, and the radius function $a(s)$ is spheroidal \eqref{eq:Las}. Here the colorbar corresponds to the value of $p^{\rm SB}(s)$. }
\label{fig:funstuff1}
\end{figure}

\begin{figure}[!ht]
\centering
\includegraphics[scale=0.5]{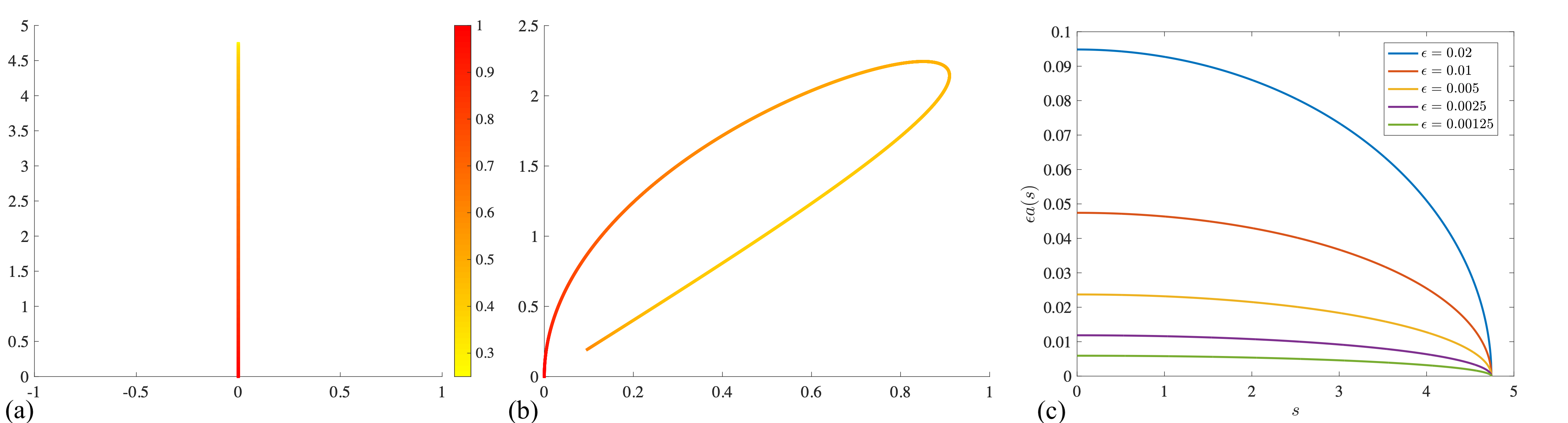}
\caption{(a),(b) The vessel centerline geometries considered in the $\epsilon$-scaling tests of figure \ref{fig:the_tests}. The color corresponds to the value of $p^{\rm SB}(s)$ using $\epsilon=0.02$, $\eta=0.05$, $\omega=10$ (i.e. the blue curve in figures \ref{fig:the_tests}a,b). (c) The profile of the radius function $\epsilon a(s)$ for $L=4.7426$ and each value of $\epsilon$ considered. }
\label{fig:test_geoms}
\end{figure}

We begin in figure \ref{fig:funstuff1} with a display of the interior pressure $p^{\rm SB}(s)$ along three different 3D vessels. 
Note that the equation \eqref{eq:pSB} is written for a unit length curve in arclength parameterization, and so must be adjusted for curves of arbitrary length $L$ which are not parameterized by arclength in order to be generally useful. Here vessel (a) has length $L=6.4708$, vessel (b) has $L=7.1502$, and vessel (c) has length $L=11.6179$. 
In all cases, the radius function $a(s)$ is given with respect to arclength $s$ by the prolate spheroid 
\begin{equation}\label{eq:Las}
  a(s) = \sqrt{L^2-s^2}
\end{equation}
with length-to-maximum-width ratio $\epsilon=0.01$ (see figure \ref{fig:test_geoms}c for a depiction of the radius profile with $L=4.7426$). For all three vessels, we take incoming pressure data $p_0=1$ and parameters $\eta=0.05$ and $\omega=10$ to allow for greater variation in pressure along the vessel. Recall that $\eta$ corresponds to the relative conductivity of the vessel interior, with small $\eta$ less permeable, and large $\omega$ serves to equalize the interior and exterior pressures across $\Gamma_\epsilon$. 
Note that a longer vessel corresponds to slower pressure decay along the length of the vessel, since here the cross sectional area scales like $\epsilon^2L^2$, so by Darcy's law the flow rate within the vessel grows linearly with $L$. The shape of the radius function profile is also more consistent for a longer vessel.

\begin{figure}[!ht]
\centering
\includegraphics[scale=0.5]{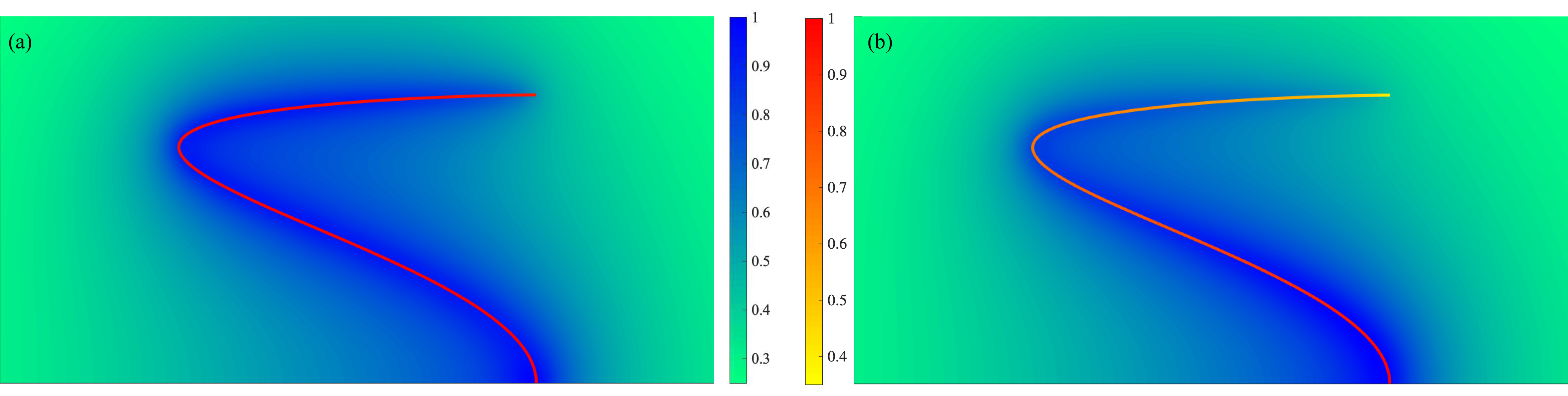}
\caption{Visualizations of a slice of the exterior pressure field $q^{\rm SB}(\bx)$ for the same vessel geometry in the case (a) $\eta=\omega=1$ and (b) $\eta=0.05$, $\omega=10$. The cool colors correspond to the value of $q^{\rm SB}$ and the warm colors correspond to the value of $p^{\rm SB}$ along the vessel. }
\label{fig:funstuff2}
\end{figure}

Next, figure \ref{fig:funstuff2} includes the exterior pressure field $q^{\rm SB}(\bx)$ \eqref{eq:qSB} surrounding an example vessel. Here for ease of visualization, the vessel centerline is planar and we plot a 2D slice of $q^{\rm SB}(\bx)$ in 3D, sampled slightly away from the vessel in the out-of-plane direction to avoid intersecting the vessel. For both examples, $p_0=1$, $\epsilon=0.01$, $L=4.5213$, and the radius function is given by \eqref{eq:Las}. The left vessel depicts the pressure fields for $\eta=\omega=1$ while the right displays the pressures for $\eta=0.05$, $\omega=10$. 
In the first case, the interior pressure $p^{\rm SB}$ remains relatively close to 1 up to the tip, and, near the vessel, the exterior pressure $q^{\rm SB}$ likewise remains close to 1 up to the tip. In the second case, the interior of the vessel is less permeable relative to the exterior, leading to a more rapid decay in pressure $p^{\rm SB}$ along the length of the vessel. Additionally, the high value of $\omega$ equalizes the interior and exterior pressures along the vessel, leading to a similar decay in $q^{\rm SB}$.

Finally, we attempt to numerically gain a better understanding of the sharpness of the bounds for $p^{\rm SB}(s)$ in Theorem \ref{thm:pSB}, particularly the $\epsilon^{-1/2}$ loss in $L^\infty$. 
We explore the $\epsilon$-dependence in $p^{\rm SB}$ for two planar vessel centerline geometries: a straight vessel and a vessel which nearly loops back to its base point (see figure \ref{fig:test_geoms}a,b). The second vessel is constructed to ensure that it does not intersect itself or the wall for any of the values of $\epsilon$ considered. Both vessels have length $L=4.7426$ and use the same radius function $a(s)$ \eqref{eq:Las}, pictured in figure \ref{fig:test_geoms}c for each of the values of $\epsilon$ considered. 

\begin{figure}[!ht]
\centering
\includegraphics[scale=0.5]{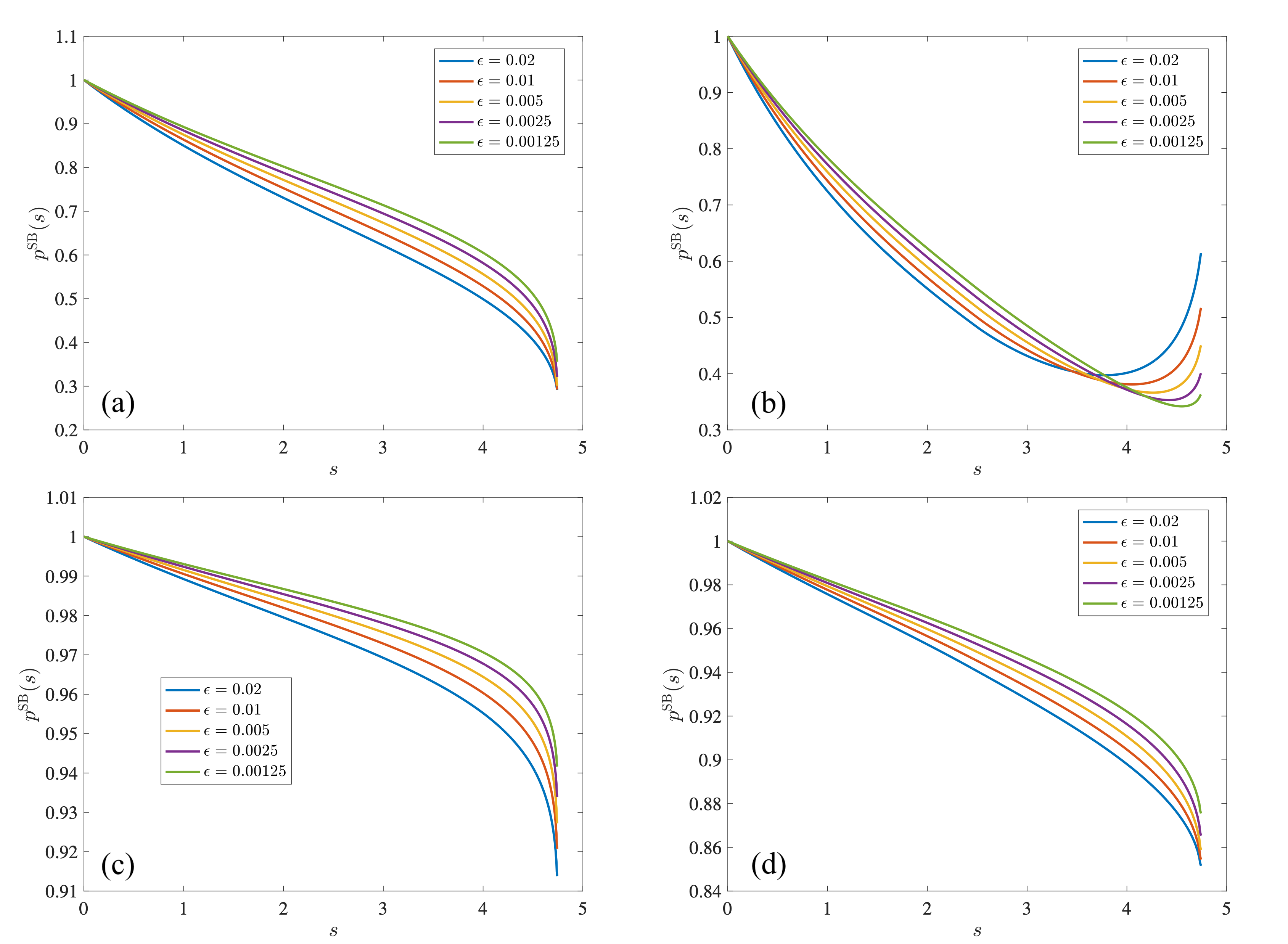}
\caption{Plots of $p^{\rm SB}(s)$ versus arclength $s$ for different values of length-to-maximum-width ratio $\epsilon$. Figures (a) and (c) are the straight vessel and (b) and (d) are the near-self-intersecting vessel pictured in figure \ref{fig:test_geoms}. In (a) and (b), the vessel walls are very permeable relative to the vessel itself: $\eta=0.05$ and $\omega=10$. In (c) and (d), the vessel walls are relatively impermeable: $\eta=1$ and $\omega=1$. In all cases, $p_0=1$.}
\label{fig:the_tests}
\end{figure}

In figure \ref{fig:the_tests}, we plot $p^{\rm SB}(s)$ versus arclength $s$ for five different values of $\epsilon$. The straight vessel is on the left, in figures (a) and (c), while the near-intersecting vessel is on the right in (b) and (d). The top figures (a) and (b) correspond to a highly permeable vessel with $\eta=0.05$, $\omega=10$. The bottom figures (c) and (d) are a less permeable vessel with $\eta=1$, $\omega=1$. 

In all cases, although the behavior at the very tip becomes more difficult to resolve for small $\epsilon$, we note a steepening of the derivative of $p^{\rm SB}$ toward the endpoint as $\epsilon$ decreases, possibly suggestive of the $\epsilon^{-1/2}$ scaling seen in the $L^\infty$ bounds for derivatives of $p^{\rm SB}$. However, it appears that $p^{\rm SB}$ itself may satisfy a uniform-in-$\epsilon$ $L^\infty$ bound, suggesting that perhaps the $\epsilon$-dependence in the lowest bound \eqref{eq:pSB_basic_ests} could be improved. 
Most notably, we are able to find non-monotonic behavior in $s$ for the near-self-intersecting vessel in the highly permeable case (figure \ref{fig:the_tests}b). In this case, the high pressure at the base of the vessel impacts the pressure in the surrounding tissue and consequently, the nearby tip of the vessel. This suggests that, while $p^{\rm SB}$ may satisfy uniform-in-$\epsilon$ $L^\infty$ bounds, there is not a maximum principle for the full free-end $p^{\rm SB}$ equation \eqref{eq:pSB} including the integral term.

\section{Preliminary expressions and estimates}\label{sec:prelim}
The proofs of Theorems \ref{thm:pSB} and \ref{thm:3D1Dto1D} will rely on the explicit representation \eqref{eq:qSB} of $q^{\rm SB}(\bx)$, which we will explore here in greater detail. In particular, we will rely on a series of nearly-singular integral estimates presented in section \ref{subsec:integral}, many of which can be directly adapted from the very similar geometric setup of \cite{free_ends}. Some important refinements will be necessary, however, to deal with the endpoint behavior in the integrodifferential equation \eqref{eq:pSB} for $p^{\rm SB}(s)$.

We begin with an explicit parameterization of the vessel and the expression \eqref{eq:qSB}.
Letting $\be_{\rm t}(s)=\X_s(s)$ denote the unit tangent vector to the vessel centerline and choosing an initial normal vector $\be_1(0)\perp\be_{\rm t}(0)$, we may define a $C^1$ orthonormal Bishop frame \cite{bishop1975there} about $\X(s)$ as $(\be_{\rm t}(s),\be_1(s),\be_2(s))$ satisfying
\begin{equation}
\frac{d}{ds}\begin{pmatrix}
\be_{\rm t}\\
\be_1\\
\be_2\\
\end{pmatrix}
= \begin{pmatrix}
0 & \kappa_1(s) &\kappa_2(s) \\
-\kappa_1(s) & 0 & 0\\
-\kappa_2(s) & 0 & 0
\end{pmatrix}\,.
\end{equation}
Here the coefficients satisfy $\kappa_1^2+\kappa_2^2=\kappa^2$, where $\kappa^2=\abs{\X_{ss}}^2$ is the centerline curvature. We define 
\begin{equation}
\kappa_\star = \norm{\kappa}_{L^\infty}\,.
\end{equation}
We also define the radial and angular unit vectors
\begin{equation}
\begin{aligned}
\be_r(s,\theta) &= \cos\theta\be_1(s) + \sin\theta\be_2(s)\,,\\
\be_\theta(s,\theta) &= -\sin\theta\be_1(s) + \cos\theta\be_2(s)\,.
\end{aligned}
\end{equation}
We may uniquely describe points $\bx$ in a neighborhood of $\Gamma_\epsilon$ using the curved cylindrical coordinate system $(r,\theta,s)$:
\begin{equation}\label{eq:param}
\bx = \X(s) + r\be_r(s,\theta)\,.
\end{equation}

Using the Bishop frame, we may also parameterize the Jacobian factor \eqref{eq:dS} as 
\begin{equation}\label{eq:free_jacfac}
\begin{aligned}
\mc{J}_\epsilon(s,\theta) &= \epsilon a(s) \sqrt{(1-\epsilon a(s)\wh\kappa(s,\theta))^2+\epsilon^2(a'(s))^2}\,,\\
\wh\kappa(s,\theta) &= \kappa_1(s)\cos\theta+\kappa_2(s)\sin\theta\,.
\end{aligned}
\end{equation}
We note that $\mc{J}_\epsilon(s,\theta)\approx \epsilon a(s)$ in the following sense: 
\begin{equation}\label{eq:jac_est}
\begin{aligned}
\abs{\mc{J}_\epsilon(s,\theta) - \epsilon a(s)} &= \epsilon a(s) \abs{\frac{(1-\epsilon a \wh\kappa)^2+\epsilon^2(a')^2-1}{\sqrt{(1-\epsilon a\wh\kappa)^2+\epsilon^2(a')^2}+1}}\\
&\le \epsilon a \bigg(\frac{2\epsilon a \abs{\wh\kappa} + \epsilon^2a^2\wh\kappa^2}{2-\epsilon a\wh\kappa} + \frac{\epsilon^2(a')^2}{1+\epsilon \abs{a'}} \bigg)
\le C(\kappa_\star,a_\star)\,\epsilon^2\,.
\end{aligned}
\end{equation}
Here $a_\star$ is as in \eqref{eq:astar}.

\begin{figure}[!ht]
\centering
\includegraphics[scale=0.4]{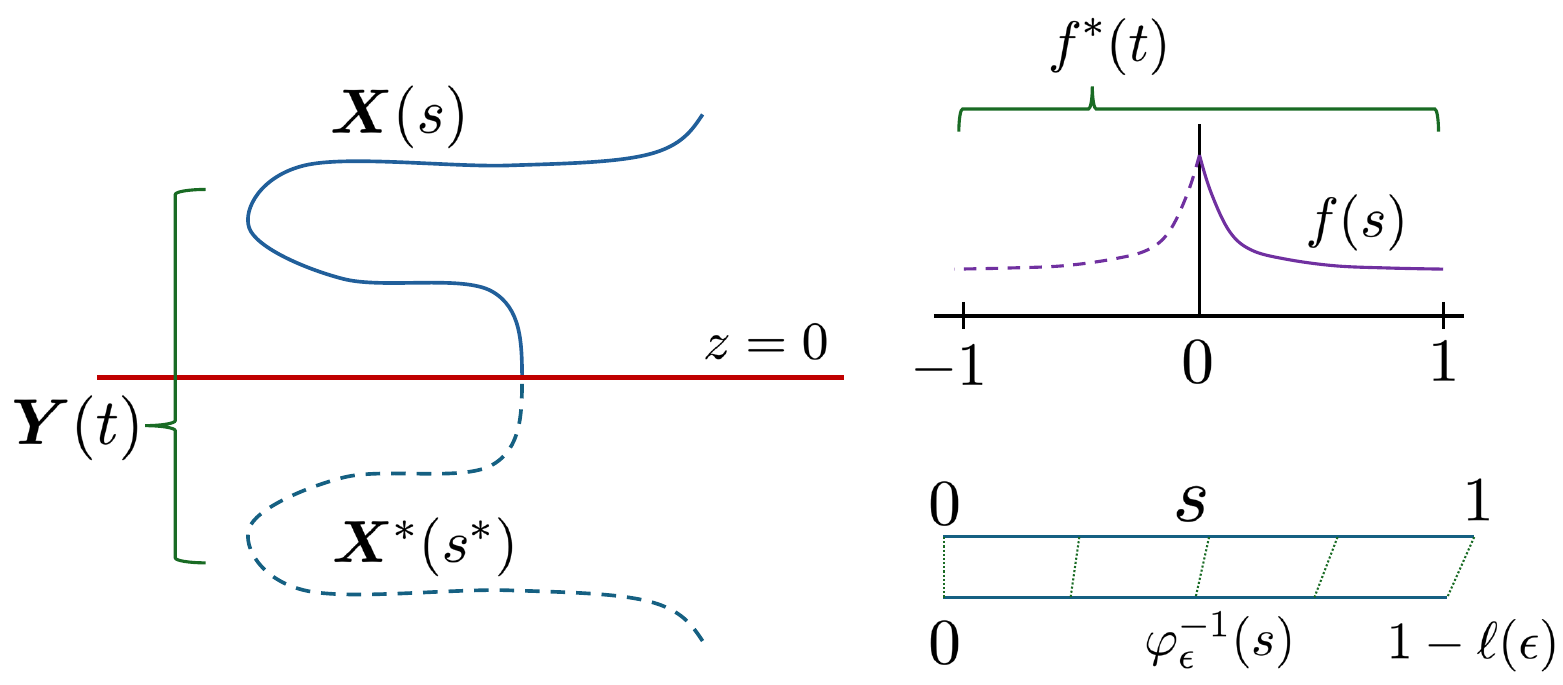}
\caption{Sketches of the reflection $\X^*$ of $\X$ across $z=0$; the extension $f^*$ of $f(s)$, $s\in[0,1]$, to $t\in[-1,1]$; and the stretch operator $\varphi_\epsilon^{-1}(s)$. }
\label{fig:reflection}
\end{figure}

Let $\X^*$ denote the reflection of the curve $\X(s)$ across the $z=0$ plane, which we parameterize with $s^*=-s\in[-1,0]$, i.e.
\begin{equation}
 \X^*(s^*) = \X(s)-2(\X(s)\cdot\be_z)\be_z\,.
\end{equation} 
 Then for $t\in[-1,1]$, we may define (see figure \ref{fig:reflection}) 
\begin{equation}\label{eq:Yt}
\Y(t) = \X^*\cup \X\,.
\end{equation}
Note that $\Y(t)$ belongs to $C^2[-1,1]$ due to the requirement that $\X_s(0)\perp\{z=0\}$. We may also extend the definition of the stretch operator $\varphi_\epsilon$ as in \eqref{eq:varphi} to the full interval $[-1,1]$ by defining
\begin{equation}\label{eq:varphi_ext}
\varphi_\epsilon(t) = \frac{t}{1-\ell(\epsilon)} \,, \qquad t\in [-1,1]\,.
\end{equation}
For $f:[0,1]\to \R$, we also define the extended function $f^*:[-1,1]\to \R$ by
\begin{equation}\label{eq:fstar}
f^*(t) = \begin{cases}
f(t)\,, &t\in [0,1]\\
f(-t)\,, &t\in[-1,0]\,.
\end{cases} 
\end{equation}
For $\bx\in \R^3_+$, we may then rewrite the operator $\mc{S}_N$ given by \eqref{eq:SN_op} as
\begin{equation}\label{eq:SN_ref}
\mc{S}_N[f](\bx) = \frac{1}{4\pi}\int_{-1}^{1}\frac{1}{\abs{\bx-\Y(\varphi_\epsilon^{-1}(t))}}\,f^*(t)\,dt\,.
\end{equation}
Note that for $f\in C^1[0,1]$, the extended function $f^*$ is continuous but not necessarily $C^1$ at $t=0$. However, we will use that the following expansions hold for any $s\in [0,1]$. For $t<0$, we may expand
\begin{equation}\label{eq:exp1}
f^*(t) = f(-t) = f(s) +(-t-s)\mc{R}_{-}(s,t)\,, \quad \abs{\mc{R}_{-}}\le C\norm{f}_{C^1[0,1]}\,.
\end{equation}
Similarly, for $t\ge 0$, we may expand
\begin{equation}\label{eq:exp2}
f^*(t) = f(t) = f(s) +(t-s)\mc{R}_{+}(s,t)\,, \quad \abs{\mc{R}_{+}}\le C\norm{f}_{C^1[0,1]}\,.
\end{equation}

Now, returning to the parameterization \eqref{eq:param}, in a neighborhood of $\Gamma_\epsilon$, we may write
\begin{equation}
\bm{R}(\bx,t) := \bx-\Y(\varphi_\epsilon^{-1}(t)) = \X(s) - \Y(\varphi_\epsilon^{-1}(t)) + r\be_r(s,\theta)\,.
\end{equation}
Within $\Omega_\epsilon$, near $\Gamma_\epsilon$, we may then calculate
\begin{equation}\label{eq:Rderivs}
\frac{\p\bR}{\p r} = \be_r(s,\theta)\,, \quad 
\frac{\p\bR}{\p\theta} = r\be_\theta(s,\theta)\,, \quad
\frac{\p\bR}{\p s} = (1-r\wh\kappa)\be_{\rm t}(s)\,.
\end{equation}
For $\bx\in \Gamma_\epsilon$, we may write
\begin{equation}
\bx(s,\theta) = \X(s) + \epsilon a(s) \be_r(s,\theta)\,,
\end{equation}
and define 
\begin{equation}\label{eq:Reps_def}
\bm{R}_\epsilon(s,t,\theta) := \bR\big|_{\Gamma_\epsilon} = \X(s) - \Y(\varphi_\epsilon^{-1}(t)) + \epsilon a(s)\be_r(s,\theta)\,.
\end{equation}
Since $\X\in C^2[0,1]$ and, by construction, $\Y\in C^2[-1,1]$, we may write
\begin{equation}\label{eq:Qdef}
\X(s) - \Y(\varphi_\epsilon^{-1}(t)) = (s-\varphi_\epsilon^{-1}(t))\be_{\rm t}(s) + (s-\varphi_\epsilon^{-1}(t))^2\bm{Q}(s,t)\,, \quad \abs{\bm{Q}(s,t)} \le \frac{\kappa_\star}{2}\,,
\end{equation}
and thus we have 
\begin{equation}\label{eq:bars}
\bm{R}_\epsilon = \bars \be_{\rm t}(s) + \epsilon a(s) \be_r(s,\theta) + \bars^2\bm{Q}(s,\bars)\,, \qquad \bars := s-\varphi_\epsilon^{-1}(t)\,.
\end{equation}
When working with integral expressions of the form \eqref{eq:SN_ref} along $\Gamma_\epsilon$, it will often be convenient to work in terms of the variable $\bars = s-\varphi_\epsilon^{-1}(t)$ instead of $t$. In particular, for $\bx\in\Gamma_\epsilon$, we may write $\mc{S}_N[f]\big|_{\Gamma_\epsilon}$ as 
\begin{equation}\label{eq:SN_Gamma}
\mc{S}_N[f]\big|_{\Gamma_\epsilon}= \frac{1}{4\pi(1-\ell(\epsilon))}\int_{s-\varphi_\epsilon^{-1}(1)}^{s+\varphi_\epsilon^{-1}(1)}\frac{1}{\abs{\bR_\epsilon}}\,f^*(\varphi_\epsilon(s-\bars))\,d\bars\,.
\end{equation}
We further note that, by Lemma 3.1 in \cite{free_ends}, the following upper and lower bounds for $\abs{\bR_\epsilon}$ hold:
\begin{equation}\label{eq:Reps_UL}
 \abs{\abs{\bR_\epsilon}-\sqrt{\bars^2+\epsilon^2a^2(s)}} \le \frac{\kappa_\star}{2}\bars^2\,, \qquad 
 \abs{\bR_\epsilon}\ge C\sqrt{\bars^2+\epsilon^2a^2(s)}\,,
\end{equation} 
where $C$ depends only on $\kappa_\star$ and $c_\Gamma$.

\subsection{Collection of integral lemmas}\label{subsec:integral}
Here we recall some general integral estimates from \cite{free_ends}, which are in turn based on results from \cite{closed_loop}. These are all generalizations of the simple calculus inequality
\begin{equation}\label{eq:calculus}
\int_{s-\varphi_\epsilon^{-1}(1)}^{s+\varphi_\epsilon^{-1}(1)} \frac{\abs{\bars}^m}{(\bars^2+\epsilon^2 a^2)^{n/2}} \le
\begin{cases}
4\abs{\log(\epsilon a)} \,,& n=m+1\\
\pi(\epsilon a)^{m+1-n}\,, & n\ge m+2
\end{cases}
\end{equation}
for integers $n>m\ge 0$ and $\epsilon,a>0$.
With the exception of Lemma \ref{lem:at_bound}, which is a significant modification of results in \cite{free_ends}, we state the following lemmas without proof; for details, see \cite[section 3.1]{free_ends}.

We begin with a basic integral estimate. The statement is a combination of Lemmas 3.3 and 3.4 in \cite{free_ends}; note that the $n=m+1$ and $n\ge m+2$ cases are erroneously swapped in Lemma 3.4 of \cite{free_ends}.

\begin{lemma}[Basic integral estimate]\label{lem:int1}
Given a radius function $a(s)$ as in \eqref{eq:astar}-\eqref{eq:spheroidal} and $\bR_\epsilon$ as in \eqref{eq:Reps_def}, for each $s\in [0,1]$, integers $m\ge 0$ and $n\ge m+1$, and $\epsilon$ sufficiently small, we have
\begin{equation}\label{eq:int1_form}
\int_{s-\varphi_\epsilon^{-1}(1)}^{s+\varphi_\epsilon^{-1}(1)}  \frac{\abs{\bars}^m}{\abs{\bR_\epsilon}^n}\,d\bars 
\le \begin{cases}
C \abs{\log\epsilon}\,, & n = m+1\\
C \min\{(\epsilon a)^{m+1-n}\,, \epsilon^{m-n}a^{m+2-n} \}\,, & n\ge m+2\,,
\end{cases}
\end{equation}
where $C$ depends only on $n$, $c_\Gamma$, and $\kappa_\star$.  
\end{lemma}

We next state a straightforward modification of Lemma 3.5 from \cite{free_ends}, a pointwise bound in which decay in a function $g(s):[0,1]\to \R$ as $s\to 1$ is traded for a factor of $\epsilon$ near the endpoint of the vessel (i.e. for $t=\varphi_\epsilon^{-1}(1)=1-\ell(\epsilon)$). 
\begin{lemma}[Trading decay for $\epsilon$]\label{lem:decay}
Given a radius function $a(s)$ as in \eqref{eq:astar}-\eqref{eq:spheroidal}, and a function $g:[0,1]\to \R$ with $a^{-2}(s)g(s)\in C[0,1]$, for each $s\in[0,1]$ and $j\ge1$, we have
\begin{equation}\label{eq:ptwisebd}
\abs{\frac{g(s)}{(s-\varphi_\epsilon^{-1}(1)+\epsilon^2a^2(s))^{j/2}}} \le C\epsilon^{-j+1}a^{-j+1}\norm{a^{-2}g}_{L^\infty(0,1)}\,,
\end{equation}
where $C$ depends on $a_\star$.
\end{lemma}

Note that the statement of Lemma 3.5 in \cite{free_ends} uses the bound $C\epsilon^{-j+1}a^{-j}\norm{a^{-1}g}_{L^\infty(0,1)}$; however, a straightforward modification of the proof yields Lemma \ref{lem:decay}. In particular, the bound \eqref{eq:ptwisebd} arises from the pointwise bound $C\epsilon^{-j+1}a^{-j-1}\abs{g}$ on the right hand side, so we may move a factor of $a^{-1}(s)$ from outside the $\norm{\cdot}_{L^\infty}$ to inside, as long as $g(s)$ decays sufficiently fast as $s\to 1$. Note that this modification will filter through the statements of the next two lemmas. In addition, we will make use of the fact that, by Lemma \ref{lem:int1}, we can trade a power of $\epsilon$ for better decay at the vessel endpoint when necessary. This tradeoff was not used in the statements of the following two lemmas in \cite{free_ends}, but we will need it here, and the proofs are not significantly altered. To summarize, the statements of the following two lemmas (\ref{lem:int2} and \ref{lem:int3}) differ slightly from their analogues in \cite{free_ends}, due to the following minor modifications: (1) Use of Lemma \ref{lem:decay} rather than Lemma 3.5 in \cite{free_ends}, and (2) Use of the tradeoff between an additional factor of $\epsilon$ versus better decay at the vessel endpoint, as stated in Lemma \ref{lem:int1}.

Our next integral estimate, taken from Lemma 3.6 in \cite{free_ends}, takes advantage of symmetry cancellations for kernels of the form \eqref{eq:int1_form} for odd powers of $\bars$, at the expense of a higher regularity requirement on the function against which the kernel is integrated. Given $g\in C^1[0,1]$, we will be applying this estimate to the extended function $g^*$ given by \eqref{eq:fstar}, which is not necessarily $C^1$ across $t=0$. Here we rely on the expansions \eqref{eq:exp1} for $t<0$ (i.e. $\bars>s$) and \eqref{eq:exp2} for $t\ge 0$ (i.e. $\bars\le s$), which both hold for $g\in C^1[0,1]$. 
\begin{lemma}[Odd cancellation]\label{lem:int2}
Given a radius function $a(s)$ as in \eqref{eq:astar}-\eqref{eq:spheroidal}, $\bR_\epsilon$ as in \eqref{eq:Reps_def}, and given a function $g\in C^1[0,1]$ with $a^{-2}(s)g(s)\in C[0,1]$, for each $s\in [0,1]$, odd integer $m\ge 0$, integer $n\ge m+2$, and $\epsilon$ sufficiently small, we have
\begin{equation}
\begin{aligned}
&\bigg|\int_{s-\varphi_\epsilon^{-1}(1)}^{s+\varphi_\epsilon^{-1}(1)}  \frac{\bars^m}{\abs{\bR_\epsilon}^n}g^*(\varphi_\epsilon(s-\bars))\,d\bars \bigg|\\
&\qquad\le \begin{cases}
C \big(\norm{g}_{C^1}\abs{\log(\epsilon a)}+\norm{a^{-2}g}_{L^\infty}\big)\,, & n = m+2\\
C \min\{(\epsilon a)^{m+2-n},\epsilon^{m+1-n}a^{m+3-n}\}\big(\norm{g}_{C^1}+\norm{a^{-2}g}_{L^\infty}\big)\,, & n\ge m+3\,,
\end{cases}
\end{aligned}
\end{equation}
where $C$ depends on $n$, $c_\Gamma$, $\kappa_\star$, and $a_\star$.
\end{lemma}
Note that the decay requirement on $g$ is stronger than in \cite{free_ends} due to Lemma \ref{lem:decay}, but in exchange we obtain better control of the growth as $s\to 1$.

The next estimate that we will need is taken from Lemma 3.7 in \cite{free_ends}, which will allow us to roughly evaluate certain types of integrals. The result again relies on the expansions \eqref{eq:exp1} and \eqref{eq:exp2} for the extension of a $C^1[0,1]$ function, and we again use the modified Lemma \ref{lem:decay} which requires stronger decay on $g$ in exchange for better control of the error term as $s\to 1$. 
\begin{lemma}[Even evaluation]\label{lem:int3}
Given a radius function $a(s)$ as in \eqref{eq:astar}-\eqref{eq:spheroidal}, $\bR_\epsilon$ as in \eqref{eq:Reps_def}, and given a function $g\in C^1[0,1]$ with $a^{-2}(s)g(s)\in C[0,1]$, for each $s\in [0,1]$, even integer $m\ge 0$, odd integer $n\ge m+3$, and $\epsilon$ sufficiently small, we have
\begin{equation}
\begin{aligned}
\bigg|\int_{s-\varphi_\epsilon^{-1}(1)}^{s+\varphi_\epsilon^{-1}(1)}  &\frac{\bars^m}{\abs{\bR_\epsilon}^n}g^*(\varphi_\epsilon(s-\bars))\,d\bars - (\epsilon a)^{m+1-n}d_{mn}\,g(s) \bigg| \\
&\le 
C \min\{(\epsilon a)^{m+2-n},\epsilon^{m+1-n}a^{m+3-n}\}\big(\norm{g}_{C^1}+\norm{a^{-2}g}_{L^\infty}\big)\,, \\
d_{mn} &= \int_{-\infty}^\infty\frac{\tau^m}{(\tau^2+1)^{n/2}}\,d\tau\,; \quad d_{03} = 2\,,
\end{aligned}
\end{equation}
where $C$ depends on $n$, $c_\Gamma$, $\kappa_\star$, and $a_\star$.
\end{lemma}

In the final lemma from \cite{free_ends}, Lemma 3.14, we make use of additional cancellations from integrating in $\theta$. In particular, the lemma applies to kernels of the form $\bars^m/\abs{\bR_\epsilon}^n$, further integrated in $\theta$ against $\cos\theta$ or $\sin\theta$. 
\begin{lemma}[Cancellation from $\theta$-integration]\label{lem:int4}
Consider $\bR_\epsilon$ as in \eqref{eq:Reps_def}, a function $g\in C[0,1]$, and integer $m>0$ and $n=m+1$ or $n=m+2$. For each $s\in [0,1]$, $0\neq k\in\Z$, $\theta_0\in\R$, and $\epsilon$ sufficiently small, we have
\begin{equation}
\begin{aligned}
&\abs{\int_{s-\varphi_\epsilon^{-1}(1)}^{s+\varphi_\epsilon^{-1}(1)} \int_0^{2\pi} \frac{\bars^mg^*(s-\bars)}{\abs{\bR_\epsilon}^n}\cos(k(\theta+\theta_0))\,d\theta\,d\bars} \\
&\qquad \le 
\begin{cases}
C\,\epsilon\abs{\log\epsilon}\norm{g}_{L^\infty}\,, & n=m+1\\
C\,\norm{g}_{L^\infty}\,, & n=m+2\,,
\end{cases}
\end{aligned}
\end{equation}
where $C$ depends on $n$, $c_\Gamma$, and $\kappa_\star$.
\end{lemma}

In addition to the previous lemmas, which are straightforward modifications of results from \cite{free_ends}, we will require one significant refinement to deal with the endpoint behavior of the $p^{\rm SB}$ equation \eqref{eq:pSB}. In particular, the following lemma trades a factor of the radius function $a^*(t)$ inside the integral for a factor of $a(s)$ outside. Here $a^*(t)$ is the extension-by-reflection \eqref{eq:fstar} to $t\in[-1,1]$. 
\begin{lemma}[Trading integrated decay for growth control]\label{lem:at_bound}
Given a radius function $a(s)$ as in \eqref{eq:astar}-\eqref{eq:spheroidal}, let $a^*(t)$ be as in \eqref{eq:fstar}. For $\bR_\epsilon$ as in \eqref{eq:Reps_def} and $\epsilon$ sufficiently small, we have
\begin{equation}\label{eq:at_bd1}
\int_{-1}^1\frac{a^*(t)}{\abs{\bR_\epsilon(s,t,\theta)}^2}\,dt \le C\epsilon^{-1}\,,
\end{equation}
where $C$ depends on $c_\Gamma$, $\kappa_\star$, and $a_0$. 
\end{lemma}

\begin{proof}
Let $\delta$ be as in \eqref{eq:delta}. For $s\le1-\delta$, the bound \eqref{eq:at_bd1} follows from Lemma \ref{lem:int1} and condition \eqref{eq:delta} since $a(s)\ge a_0$ for $a_0$ independent of $\epsilon$. 

For $s>1-\delta$, we will rely on the spheroidal nature \eqref{eq:spheroidal} of the radius function. First, by \eqref{eq:Reps_UL}, we have
\begin{equation}
\abs{\bR_\epsilon}^2\ge C\big((s-\varphi_\epsilon^{-1}(t))^2+\epsilon^2a^2(s)\big) =: D(t,s)\,.
\end{equation}
For all $s>1-\delta$, there exists a $k$ such that $s=1-\ell(\epsilon)^k$. We can see $0<\frac{\log\delta}{\log(\ell(\epsilon))}<k\in\R$, so we have
\begin{equation}
\begin{aligned}
D(t,s) &= (1-\ell(\epsilon))^2\bigg(\frac{s}{1-\ell(\epsilon)}-t\bigg)^2+\epsilon^2a^2(s) \\
&= (1-\ell(\epsilon))^2\bigg(1-t + \frac{\ell(\epsilon)-\ell(\epsilon)^k}{1-\ell(\epsilon)}\bigg)^2
+\epsilon^2\ell(\epsilon)^k\big(2-\ell(\epsilon)^k\big) + A_\epsilon \\
&= (1-\ell(\epsilon))^2(1-t)^2 + 2(1-\ell(\epsilon))(\ell(\epsilon)-\ell(\epsilon)^k)(1-t) + (\ell(\epsilon)-\ell(\epsilon)^k)^2\\
&\qquad +\epsilon^2\ell(\epsilon)^k\big(2-\ell(\epsilon)^k\big) + A_\epsilon\,,
\end{aligned}
\end{equation}
where $\abs{A_\epsilon}\le C\epsilon^{4+2k}$.
Now, for $k\ge1$, we have $\ell(\epsilon)-\ell(\epsilon)^k\ge 0$ and thus
\begin{equation}
D(t,s) \ge (1-\ell(\epsilon))^2(1-t)^2 + \ell(\epsilon)^2\big(1-\ell(\epsilon)^{k-1}\big)^2 +\epsilon^{2+2k} 
\ge C\big((1-t)^2 + \epsilon^4\big)\,.
\end{equation}
Using the spheroidal form \eqref{eq:spheroidal} of $a^*(t)$ as $t\to 1$, we then have 
\begin{equation}
\begin{aligned}
\int_{-1}^1\frac{a^*(t)}{\abs{\bR_\epsilon}^2}\,dt &\le 
C\int_{-1}^1\frac{\sqrt{1-t^2}}{D(t,s)}\,dt 
\le C\int_{-1}^1\frac{1}{\big[(1-t)^2 + C\epsilon^4\big]^{3/4}}\,dt 
\le C\epsilon^{-1}\,.
\end{aligned}
\end{equation}
Here we have used that the final integral can be evaluated explicitly as the hypergeometric function $C\epsilon^{-3}\,_2F_1\big(\frac{1}{2},\frac{3}{4};\frac{3}{2};-\frac{C}{\epsilon^4}\big)$, where $_2F_1\big(c_1,c_2;c_3;-\frac{C}{\epsilon^4}\big)\le C\epsilon^2$ as $\epsilon\to 0$.

For $\frac{\log\delta}{\log(\ell(\epsilon))}<k<1$, we have
\begin{equation}
D(t,s)\ge C\big[\big(1-t - m(\epsilon)\big)^2 + \epsilon^{2+2k}\big]\,, \qquad m(\epsilon) = \frac{\ell(\epsilon)^k-\ell(\epsilon)}{1-\ell(\epsilon)}>0\,.
\end{equation}
Note that $m(\epsilon)\le C\epsilon^{2k}$.
Furthermore, using the spheroidal endpoint condition \eqref{eq:spheroidal}, we have
\begin{equation}
a^*(t)\le C\sqrt{1-t^2}\le C\big( \sqrt{1-t-m(\epsilon)}+\sqrt{m(\epsilon)}\big)\,.
\end{equation}
Thus
\begin{equation}
\begin{aligned}
\int_{-1}^1\frac{a^*(t)}{\abs{\bR_\epsilon}^2}\,dt &\le 
C\int_{-1}^1\frac{\sqrt{1-t-m(\epsilon)}}{D(t,s)}\,dt + C\sqrt{m(\epsilon)}\int_{-1}^1\frac{1}{D(t,s)}\,dt \\
&\le C\int_{-1}^1\frac{1}{\big[(1-m(\epsilon) -t)^2+\epsilon^{2+2k} \big]^{3/4}}\,dt + C\int_{-1}^1\frac{\sqrt{m(\epsilon)}}{(1-m(\epsilon) -t)^2+\epsilon^{2+2k}}\,dt \\
&\le C\epsilon^{-(1+k)/2} + C\sqrt{m(\epsilon)}\epsilon^{-(1+k)} 
\le C\epsilon^{-1}\,,
\end{aligned}
\end{equation}
where we have used that $k<1$ and $\sqrt{m(\epsilon)}\le C\epsilon^k$. Here we have again used that both integrals can be evaluated directly, the first as the hypergeometric function $C\epsilon^{-3(1+k)/2}\,_2F_1\big(\frac{1}{2},\frac{3}{4};\frac{3}{2};-\frac{C}{\epsilon^{2(1+k)}}\big)$ as before, and the second as $C\epsilon^{-(1+k)}\big[\arctan(\frac{2-\sqrt{m(\epsilon)}}{\epsilon})+\arctan(\frac{\sqrt{m(\epsilon)}}{\epsilon})\big]$.
\end{proof}

\section{Estimates for the 1D integrodifferential equation}\label{sec:pSB}
We now consider the integrodifferential equation \eqref{eq:pSB} for $p^{\rm SB}:[0,1]\to\R$. In order to show the well-posedness results and $\epsilon$-dependent bounds of Theorem \ref{thm:pSB}, we will first need to use the integral estimates of section \ref{subsec:integral} to obtain a series of lemmas regarding the behavior of the integral operator. We provide the statements of the following three lemmas first, and then proceed to prove each in the following subsections.

\subsection{Kernel lemmas}\label{subsec:kernel}
For $\bR_\epsilon$ as in \eqref{eq:Reps_def}, we define the kernel $K_\epsilon(s,t)$ by 
\begin{equation}\label{eq:Keps}
K_\epsilon(s,t) = \frac{\eta}{8\pi^2}\int_0^{2\pi}\frac{1}{\abs{\bR_\epsilon(s,t,\theta)}}\, d\theta\,,
\end{equation}
and note that this is the kernel of the integral operator appearing in the $p^{\rm SB}$ equation \eqref{eq:pSB}. The behavior of this kernel dictates the behavior of the equation for $p^{\rm SB}$, particularly its $\epsilon$-dependence, so we need a detailed understanding of its properties.

First, just to obtain existence and uniqueness of a solution to the $p^{\rm SB}$ equation, we will require the following lemma. 
\begin{lemma}[Near-positivity of integral operator]\label{lem:pos_op}
Let $a(s)$ be a radius function as in \eqref{eq:astar}-\eqref{eq:spheroidal} and let $K_\epsilon(s,t)$ be given by \eqref{eq:Keps}. Given $f:[0,1]\to \R$ with $a^{-1/2}f\in L^2(0,1)$, let $f^*(t)$ be as in \eqref{eq:fstar}. We then have
\begin{equation}\label{eq:pos_op}
\begin{aligned}
\int_0^1\int_{-1}^1 K_\epsilon(s,t)\, f^*(t)f(s)\,dt\,ds \ge -C\epsilon^{1/2}\abs{\log\epsilon}^{1/2}\|a^{-1/2}f\|_{L^2(0,1)}^2\,,
\end{aligned}
\end{equation}
where $C$ is bounded independent of $\epsilon$ as $\epsilon\to 0$. 
\end{lemma}
In particular, we will rely on the fact that the operator \eqref{eq:pos_op} is nearly positive (and in fact would be positive if the vessel had no endpoint). 
The proof of Lemma \ref{lem:pos_op} appears below in section \ref{subsubsec:near_pos}.

To obtain the bounds for $(a^4p^{\rm SB}_s)_{ss}$ stated in Theorem \ref{thm:pSB}, we will rely on the fact that the kernel $K_\epsilon(s,t)$ is nearly a convolution kernel in $s$ and $t$, so that its partial derivatives with respect to $s$ and $t$ are nearly antisymmetric. 
\begin{lemma}[Near-antisymmetry of kernel derivatives]\label{lem:sym_op}
Let $a(s)$ be a radius function as in \eqref{eq:astar}-\eqref{eq:spheroidal} and let $K_\epsilon(s,t)$ be given by \eqref{eq:Keps}. Given $f:[0,1]\to \R$ with $a^{-1}f\in L^\infty(0,1)$ and $af_s\in L^\infty(0,1)$, let $f^*(t)$ be as in \eqref{eq:fstar}. We then have
\begin{equation}\label{eq:sym_decomp}
\begin{aligned}
\int_{-1}^1\p_sK_\epsilon(s,t)\,f^*(t)\,dt &= -\int_{-1}^1\p_tK_\epsilon(s,t)\,f^*(t)\,dt + H_\epsilon[f](s)\,,
\end{aligned}
\end{equation}
where the remainder term satisfies
\begin{equation}
\begin{aligned}
\|a^2H_\epsilon[f]\|_{L^2} 
&\le C\big(\|a^{-1/2}f\|_{L^2} + \epsilon\norm{a^{-1}f}_{L^\infty}+ \epsilon\norm{a\,f_s}_{L^\infty}\big)\,, \\
\norm{a^2H_\epsilon[f]}_{L^\infty} 
&\le C\big(\norm{a^{-1}f}_{L^\infty} + \epsilon\norm{a\,f_s}_{L^\infty} \big)\,.
\end{aligned}
\end{equation}
\end{lemma}
The proof of Lemma \ref{lem:sym_op} is below in section \ref{subsubsec:sym_op}.
The near-antisymmetry of Lemma \ref{lem:sym_op} is necessary but not quite sufficient to obtain estimates on higher derivatives of $p^{\rm SB}$, due to the endpoint weights. We will therefore need an additional lemma so that integration by parts within the integral operator makes sense up to the vessel endpoint. 
\begin{lemma}[Weighted integration by parts]\label{lem:kernelIBP}
Let $a(s)$ be a radius function as in \eqref{eq:astar}-\eqref{eq:spheroidal} and let $K_\epsilon(s,t)$ be given by \eqref{eq:Keps}. Given $f:[0,1]\to \R$ with $a^{-1}f\in L^\infty(0,1)$ and $a^2f_s\in L^1(0,1)$, let $a^*(t)$, $f^*(t)$ be as in \eqref{eq:fstar}. We then have
\begin{equation}
 a^2(s)\int_{-1}^1 \p_tK_\epsilon(s,t)\, f^*(t)\,dt = 
 -\int_{-1}^1 a^*(t)^2K_\epsilon(s,t)\, f_t^*(t)\,dt
 + G_\epsilon[f](s) 
\end{equation} 
where
\begin{align}
\norm{a^{1/2}G_\epsilon[f]}_{L^2}&\le C\abs{\log\epsilon}\|a^{-1/2}f\|_{L^2(0,1)}\,, \label{eq:GepsL2}\\
\norm{G_\epsilon[f]}_{L^\infty} &\le C\epsilon^{-1/2}\|a^{-1/2}f\|_{L^2(0,1)} + C\epsilon\norm{a^{-1}f}_{L^\infty(0,1)} \label{eq:GepsLinfty}
\end{align}
for constants $C$ independent of $\epsilon$.
\end{lemma}
The proof of Lemma \ref{lem:kernelIBP} appears in section \ref{subsubsec:IBP}.

%
\subsubsection{Proof of Lemma \ref{lem:pos_op}: near-positivity of integral operator}\label{subsubsec:near_pos}
We begin by recalling the definition \eqref{eq:Yt} of the extended curve $\bm{Y}(t)$, where the extension is by reflection of $\X(s)$ across the $z=0$ plane (recall figure \ref{fig:reflection}). Along $\Y(t)$, for $-1\le t<0$, we take $\be_r(t,\theta)=\be_r(-t,\theta)-2\big( \be_r(-t,\theta)\cdot\be_z\big)\be_z$ to be the reflection of the radial vector $\be_r(t,\theta)$ for $t\ge 0$. 
For $\tau,t\in [-1,1]$, we define the extended kernel 
\begin{equation}\label{eq:KY}
K_{\rm Y}(\tau,t) = \frac{\eta}{8\pi^2}\int_0^{2\pi}\frac{1}{\abs{\Y(\tau) - \Y(\varphi_\epsilon^{-1}(t))+\epsilon a^*(\tau)\be_r(\tau,\theta)}}\,d\theta\,,
\end{equation}
where $a^*$ is as in \eqref{eq:fstar}.
In addition, we define the integral
\begin{equation}
I_\epsilon:= \frac{1}{2}\int_{-1}^1\int_{-1}^1\bigg(K_{\rm Y}(\tau,t) + K_{\rm Y}(t,\tau) \bigg) f^*(\tau)f^*(t)\,dt\,d\tau 
\end{equation}
of the symmetric part of $K_{\rm Y}(\tau,t)$. Note that we have 
\begin{equation}
\begin{aligned}
I_\epsilon = \int_{-1}^1\int_{-1}^1 K_{\rm Y}(\tau,t) \,f^*(t)f^*(\tau)\,dt\,d\tau = 2\int_0^1\int_{-1}^1 K_\epsilon(s,t)\, f^*(t)f(s)\,dt\,ds\,.
\end{aligned}
\end{equation}
The first equality follows since $\frac{1}{2}\big(K_{\rm Y}(\tau,t) + K_{\rm Y}(t,\tau) \big)=K_{\rm Y}(\tau,t) - \frac{1}{2}\big(K_{\rm Y}(\tau,t) - K_{\rm Y}(t,\tau) \big)$ and
\begin{equation}\label{eq:antisym}
\int_{-1}^1\int_{-1}^1 \big(K_{\rm Y}(\tau,t) - K_{\rm Y}(t,\tau) \big) \,f^*(t)f^*(\tau)\,dt\,d\tau = 0\,.
\end{equation}
The second equality is due to even reflections across $z=0$ and $t=0$. In particular, using that the reflected curve $\X^*$ satisfies $\abs{\X^*(t)}=\abs{\X(-t)}$ for $-1\le t<0$, we have 
\begin{equation}
\begin{aligned}
&\frac{8\pi^2}{\eta}\int_{-1}^0\int_{-1}^1 K_{\rm Y}(\tau,t) \,f^*(t)f^*(\tau)\,dt\,d\tau\\
%
&\quad =  \int_{-1}^0\int_{-1}^0\int_0^{2\pi}\frac{f(-t)f(-\tau)}{\abs{\X(-\tau) - \X(-\varphi_\epsilon^{-1}(t))+\epsilon a(-\tau)\be_r(-\tau,\theta)}}\,d\theta\,dt\,d\tau \\
&\qquad + \int_{-1}^0\int_{0}^1\int_0^{2\pi}\frac{f(t)f(-\tau)}{\abs{\X(-\tau) - \X^*(-\varphi_\epsilon^{-1}(t))+\epsilon a(-\tau)\be_r(-\tau,\theta)}}\,d\theta\,dt\,d\tau \\
&\quad =  \int_{0}^1\int_{0}^1\int_0^{2\pi}\frac{f(t)f(\tau)}{\abs{\X(\tau) - \X(\varphi_\epsilon^{-1}(t))+\epsilon a(\tau)\be_r(\tau,\theta)}}\,d\theta\,dt\,d\tau \\
&\qquad + \int_{0}^1\int_{-1}^0\int_0^{2\pi}\frac{f(-t)f(\tau)}{\abs{\X(\tau) - \X^*(\varphi_\epsilon^{-1}(t))+\epsilon a(\tau)\be_r(\tau,\theta)}}\,d\theta\,dt\,d\tau \\
&\quad = \frac{8\pi^2}{\eta}\int_0^1\int_{-1}^1 K_\epsilon(s,t)\, f^*(t)f(s)\,dt\,ds\,.
\end{aligned}
\end{equation}
Thus, to prove Lemma \ref{lem:pos_op}, it suffices to show that 
\begin{equation}\label{eq:I_toshow}
I_\epsilon \ge -C\epsilon^{1/2}\abs{\log\epsilon}^{1/2}\|a^{-1/2}f\|_{L^2(0,1)}^2\,.
\end{equation}

Recalling that $\varphi_\epsilon^{-1}(t) = (1-\ell(\epsilon))t$ for $\ell(\epsilon)$ as in \eqref{eq:varphi}, we denote $L_\epsilon := 1-\ell(\epsilon)$ and consider the interval $[-L_\epsilon,L_\epsilon]$, i.e. up to $O(\epsilon^2)$ from the extended vessel endpoints. We decompose $I_\epsilon$ as 
\begin{equation}
\begin{aligned}
I_\epsilon &= I_{\epsilon,0} + I_{\epsilon,1}+ I_{\epsilon,2} + I_{\epsilon,3} \,,\\
I_{\epsilon,0}&= \frac{1}{2}\int_{-L_\epsilon}^{L_\epsilon}\int_{-L_\epsilon}^{L_\epsilon} \bigg(K_{\rm Y}(\tau,\varphi_\epsilon(t)) + K_{\rm Y}(t,\varphi_\epsilon(\tau)) \bigg) \, f^*(t)f^*(\tau)\, dt\,d\tau\,, \\
I_{\epsilon,1}&= \int_{-L_\epsilon}^{L_\epsilon}\int_{-L_\epsilon}^{L_\epsilon} \bigg(K_{\rm Y}(\tau,t)- K_{\rm Y}(\tau,\varphi_\epsilon(t))\bigg) \, f^*(t)f^*(\tau)\, dt\,d\tau\,, \\
I_{\epsilon,2}&= \bigg(\int_{L_\epsilon< \abs{\tau}\le 1}\int_{-L_\epsilon}^{L_\epsilon} + \int_{-L_\epsilon}^{L_\epsilon}\int_{L_\epsilon<\abs{t}\le 1} \bigg)K_{\rm Y}(\tau,t) \, f^*(t)f^*(\tau)\, dt\,d\tau \,,\\
I_{\epsilon,3} &= \int_{L_\epsilon<\abs{\tau}\le 1}\int_{L_\epsilon< \abs{t}\le 1}K_{\rm Y}(\tau,t) \, f^*(t)f^*(\tau)\, dt\,d\tau\,.
\end{aligned}
\end{equation}
Here we note that the kernel in $I_{\epsilon,0}$ undoes the inverse stretch operation in the second argument, and therefore blows up at the vessel endpoint. We must therefore consider $I_{\epsilon,0}$ over $[-L_\epsilon,L_\epsilon]\times[-L_\epsilon,L_\epsilon]$. The three additional integrals are remainder terms that arise due to our treatment of the vessel endpoint. Note that by \eqref{eq:antisym}, we may just consider $K_{\rm Y}(\tau,t)$ in these terms instead of its symmetrization. 
We will show that the integral $I_{\epsilon,0}$ is non-negative, while the three remaining integrals are small in $\epsilon$. 

We begin by noting that the symmetric kernel $\frac{1}{2}(K_{\rm Y}(\tau,\varphi_\epsilon(t)) + K_{\rm Y}(t,\varphi_\epsilon(\tau)))$ is positive definite on $[-L_\epsilon,L_\epsilon]\times [-L_\epsilon,L_\epsilon]$ in the following sense (as stated in \cite{bochner1959lectures,fasshauer2011positive}). For any $n\in \N$, sequence $t_j$, $j=1,\dots,n$, in $[-L_\epsilon,L_\epsilon]$, and values $c_j\in \R$, $j=1,\dots,n$, we have
\begin{equation}\label{eq:posdef}
\frac{1}{2}\sum_{i=1}^n\sum_{j=1}^n\big( K_{\rm Y}(t_i,\varphi_\epsilon(t_j))+K_{\rm Y}(t_j,\varphi_\epsilon(t_i))\big)\,c_ic_j \ge 0\,.
\end{equation}
To see that $K_{\rm Y}$ satisfies \eqref{eq:posdef}, we note that since $\Y\in C^2[-1,1]$, we may expand $\Y(\tau) - \Y(t)=(\tau-t)\be_{\rm t}(\tau) + (\tau-t)^2\bm{Q}(\tau,t)$ for $\bm{Q}$ as in \eqref{eq:Qdef}.
Thus for $\tau,t \in[-L_\epsilon,L_\epsilon]$, we have
\begin{equation}
\begin{aligned}
\abs{\Y(\tau) - \Y(t)+\epsilon a^*(\tau)\be_r(\tau,\theta)}^2
&= \abs{\Y(\tau) - \Y(t)}^2 + \epsilon a^*(\tau)\,(\tau-t)^2\be_r\cdot\bm{Q} + \epsilon^2 (a^*)^2 \\
&\ge \epsilon^2 \big(a^*(\tau)\big)^2
\end{aligned}
\end{equation}
for $\epsilon$ sufficiently small. Then, for $\tau,t\in [-L_\epsilon,L_\epsilon]$, we have
\begin{equation}\label{eq:diag_dom}
\begin{aligned}
\frac{1}{2}\big(K_{\rm Y}(\tau,\varphi_\epsilon(t)) + K_{\rm Y}(t,\varphi_\epsilon(\tau))\big) &= \frac{\eta}{16\pi^2}\int_0^{2\pi}\bigg(\frac{1}{\abs{\Y(\tau) - \Y(t)+\epsilon a^*(\tau)\be_r(\tau,\theta)}}\\
&\hspace{2cm} + \frac{1}{\abs{\Y(t) - \Y(\tau)+\epsilon a^*(t)\be_r(t,\theta)}}\bigg)\,d\theta \\
&\le \frac{\eta}{8\pi}\frac{1}{\epsilon^2 a^*(\tau)a^*(t)} 
= K_{\rm Y}(\tau,\varphi_\epsilon(\tau))^{1/2} K_{\rm Y}(t,\varphi_\epsilon(t))^{1/2}\,.
\end{aligned}
\end{equation}
The diagonal dominance \eqref{eq:diag_dom} then implies \eqref{eq:posdef}. Due to the positive definiteness \eqref{eq:posdef} of $\frac{1}{2}(K_{\rm Y}(\tau,\varphi_\epsilon(t)) + K_{\rm Y}(t,\varphi_\epsilon(\tau)))$ on $[-L_\epsilon,L_\epsilon]\times[-L_\epsilon,L_\epsilon]$, we have (see \cite[Definitions 2 \& 3]{fasshauer2011positive} or \cite[Chapter 4, section 20]{bochner1959lectures})
\begin{equation}\label{eq:Ieps0}
I_{\epsilon,0}\ge 0\,.
\end{equation}

It remains to show that $I_{\epsilon,1}$, $I_{\epsilon,2}$, and $I_{\epsilon,3}$ are small in $\epsilon$. We begin with $I_{\epsilon,1}$, where smallness arises because the stretch operator $\varphi_\epsilon(t)$ as in \eqref{eq:varphi} is at worst $O(\epsilon^2)$ away from $t$.
As in \eqref{eq:Reps_def}, we will use $\bR_\epsilon(\tau,t,\theta)$ to denote 
\begin{equation}\label{eq:YReps}
\bR_\epsilon(\tau,t,\theta) =  \Y(\tau)- \Y(\varphi_\epsilon^{-1}(t)) +\epsilon a^*(\tau)\be_r(\tau,\theta)\,.
\end{equation}
We use the same notation as in \eqref{eq:Reps_def} since the only difference is that $\tau$ is now considered over the entire reflected curve $\Y$. We emphasize that the integral estimates of section \ref{subsec:integral} continue to hold.  
We additionally define 
\begin{equation}
\bR_{0,\epsilon}(\tau,t,\theta) =  \Y(\tau)- \Y(t) +\epsilon a^*(\tau)\be_r(\tau,\theta)\,.
\end{equation}

Expanding 
\begin{equation}
\Y(\varphi_\epsilon^{-1}(t)) = \Y(t) -\ell(\epsilon)t \,\bm{Q}_{\rm Y}(t)\,, \qquad  \abs{\bm{Q}_{\rm Y}(t)}=1\,,
\end{equation}
we note that
\begin{equation}
\bR_\epsilon = \bR_{0,\epsilon} + \ell(\epsilon)t \,\bm{Q}_{\rm Y}(t)\,.
\end{equation}
For $t,\tau\in [-L_\epsilon,L_\epsilon]$, we then have
\begin{equation}
\begin{aligned}
 \abs{K_{\rm Y}(\tau,t)- K_{\rm Y}(\tau,\varphi_\epsilon(t))} &= \abs{\frac{\eta}{8\pi^2}\int_0^{2\pi}\bigg(\frac{1}{\abs{\bR_\epsilon}} - \frac{1}{\abs{\bR_{0,\epsilon}}} \bigg)\,d\theta}\\
&= \abs{\frac{\eta}{8\pi^2}\int_0^{2\pi}\frac{-2\ell(\epsilon)t\bm{Q}_{\rm Y}\cdot\bR_{0,\epsilon}- \ell(\epsilon)^2t^2\abs{\bm{Q}}^2}{\abs{\bR_{0,\epsilon}}\abs{\bR_\epsilon}(\abs{\bR_{0,\epsilon}}+\abs{\bR_\epsilon})}  } \\
&\le C\int_0^{2\pi}\frac{\ell(\epsilon)+\ell(\epsilon)^{3/2}a^*(\tau)^{-1}}{\abs{\bR_\epsilon}^2}\,d\theta
\le \int_0^{2\pi}\frac{C\epsilon^2}{\abs{\bR_\epsilon}^2}\,d\theta\,,
 \end{aligned}
\end{equation}
where we have used the spheroidal endpoint condition \eqref{eq:spheroidal} to obtain $a^*(\tau)\ge C\epsilon$ for $\tau\in [-L_\epsilon,L_\epsilon]$.  
We may then estimate $I_{\epsilon,1}$ as
\begin{equation}
\begin{aligned}
\abs{I_{\epsilon,1}}&\le C\epsilon^2\int_{-L_\epsilon}^{L_\epsilon}\int_{-L_\epsilon}^{L_\epsilon} \int_0^{2\pi}\frac{a^*(\tau)^{1/2}a^*(t)^{1/2}}{\abs{\bR_\epsilon}^2} \, a^*(t)^{-1/2}f^*(t)\, a^*(\tau)^{-1/2}f^*(\tau)\,d\theta\, dt\,d\tau\\
&\le C\epsilon^2\|a^{-1/2}f\|_{L^2(0,1)}\bigg\| \underbrace{\int_0^{2\pi}\int_{-L_\epsilon}^{L_\epsilon} \frac{a^*(\cdot)^{1/2}a^*(t)^{1/2}}{\abs{\bR_\epsilon(\cdot,t,\theta)}^2} \, a^*(t)^{-1/2}f^*(t)\,dt\,d\theta }_{T_\epsilon[a^{-1/2}f]}  \bigg\|_{L^2(-L_{\epsilon},L_\epsilon)} \,.
\end{aligned}
\end{equation}
To estimate the integral operator $T_\epsilon[a^{-1/2}f]$, we will rely on the following version of Schur's test (see \cite{tao2004lecture}, \cite[Chapter 5]{halmos2012bounded}). Given two intervals $I_1,I_2\subset \R$ and a kernel $K:I_1\times I_2\to \R$ satisfying
\begin{equation}\label{eq:schur1}
\begin{aligned}
\sup_{\tau\in I_1}\int_{I_2}\abs{K(\tau,t)}\,dt \le A\,, \qquad 
\sup_{t\in I_2}\int_{I_1}\abs{K(\tau,t)}\,d\tau \le B
\end{aligned}
\end{equation}
for some $0<A,B<\infty$, then for all $f\in L^2(I_2)$, we have
\begin{equation}\label{eq:schur2}
\norm{\int_{I_2}K(\tau,t)\,f(t)\,dt}_{L^2(I_1)}
\le A^{1/2}B^{1/2}\norm{f}_{L^2(I_2)}\,.
\end{equation}

 We begin by applying the Cauchy-Schwarz inequality followed by Lemmas \ref{lem:int1} on the first term in the product and \ref{lem:at_bound} with $m=0$ and $n=2$ on the second term in the product to see
\begin{equation}
A_\epsilon:=\int_{-L_\epsilon}^{L_\epsilon}\frac{a^*(\tau)^{1/2}a^*(t)^{1/2}}{\abs{\bR_\epsilon(\tau,t,\theta)}^2}\,dt
 \le \bigg(\int_{-L_\epsilon}^{L_\epsilon}\frac{a^*(\tau)}{\abs{\bR_\epsilon}^2}\,dt\bigg)^{1/2}\bigg(\int_{-L_\epsilon}^{L_\epsilon}\frac{a^*(t)}{\abs{\bR_\epsilon}^2}\,dt\bigg)^{1/2} \le C\epsilon^{-1}\,.
\end{equation}
On the other hand, integrating with respect to $\tau$, we may use \eqref{eq:Reps_UL} to bound $\abs{\bR_\epsilon}^2\ge C((\tau-\varphi_\epsilon^{-1}(t))^2+\epsilon^4)$. Taking $\bars=\tau-\varphi_\epsilon^{-1}(t)$ as in \eqref{eq:bars}, and recalling the simple calculus estimate \eqref{eq:calculus}, we have
\begin{equation}
\begin{aligned}
B_\epsilon&:=\int_{-L_\epsilon}^{L_\epsilon}\frac{a^*(\tau)^{1/2}a^*(t)^{1/2}}{\abs{\bR_\epsilon(\tau,t,\theta)}^2}\,d\tau
 \le \bigg(\int_{-L_\epsilon}^{L_\epsilon}\frac{a^*(\tau)}{\abs{\bR_\epsilon}^2}\,d\tau\bigg)^{1/2}\bigg(\int_{-L_\epsilon}^{L_\epsilon}\frac{a^*(t)}{\abs{\bR_\epsilon}^2}\,d\tau\bigg)^{1/2} \\
 &\le C\bigg(\epsilon^{-1}\int_{-L_\epsilon-\varphi_\epsilon^{-1}(t)}^{L_\epsilon-\varphi_\epsilon^{-1}(t)}\frac{1}{\sqrt{\bars^2+\epsilon^4}}\,d\bars\bigg)^{1/2}\bigg(\int_{-L_\epsilon-\varphi_\epsilon^{-1}(t)}^{L_\epsilon-\varphi_\epsilon^{-1}(t)}\frac{1}{\bars^2+\epsilon^4}\,d\bars\bigg)^{1/2}
 \le C\epsilon^{-3/2}\abs{\log\epsilon}^{1/2}\,.
 \end{aligned}
\end{equation}
Then, using Schur's test \eqref{eq:schur2}, we may bound 
\begin{equation}
\norm{T_{\epsilon}[a^{-1/2}f]}_{L^2(-L_\epsilon,L_\epsilon)} \le C A_\epsilon^{1/2}B_\epsilon^{1/2}\|a^{-1/2}f\|_{L^2(0,1)}
\le C\epsilon^{-5/4}\abs{\log\epsilon}^{1/4}\|a^{-1/2}f\|_{L^2(0,1)}\,.
\end{equation}
We thus obtain
\begin{equation}
\abs{I_{\epsilon,1}} \le C\epsilon^{3/4}\abs{\log\epsilon}^{1/4}\|a^{-1/2}f\|_{L^2(0,1)}^2\,.
\end{equation}

For $I_{\epsilon,2}$ and $I_{\epsilon,3}$, we will rely on both the smallness of the interval $[L_\epsilon,1]$, recalling that $L_\epsilon-1=\ell(\epsilon)\le C\epsilon^2$, and the smallness of the vessel radius $a(s)$ over this interval. In particular, by the spheroidal endpoint condition \eqref{eq:spheroidal}, $a(s)\le C\epsilon$ for $s\in [L_\epsilon,1]$.

We being with the simpler term $I_{\epsilon,3}$. By the definition \eqref{eq:KY} of the kernel $K_{\rm Y}(\tau,t)$, we have
\begin{equation}\label{eq:Klower}
\abs{K_{\rm Y}(\tau,t)}\ge \frac{C}{\epsilon^2}
\end{equation}
for all $\tau,t\in [-1,1]$. We thus have
\begin{equation}
\begin{aligned}
\abs{I_{\epsilon,3}}&\le \frac{C}{\epsilon^2}\bigg(\int_{L_\epsilon<\abs{t}\le 1}a^*(t)^{1/2}\abs{a^*(t)^{-1/2}f^*(t)}\,dt\bigg)^2 \\
&\le \frac{C}{\epsilon^2}\|a^{-1/2}f\|_{L^2(0,1)}^2\int_{L_\epsilon<\abs{t}\le 1}a^*(t)\,dt 
\le C\epsilon\|a^{-1/2}f\|_{L^2(0,1)}^2\,,
\end{aligned}
\end{equation}
where we have used that the endpoint integral is $O(\epsilon^3)$.

For $I_{\epsilon,2}$, we may begin by bounding
\begin{equation}
\begin{aligned}
\abs{I_{\epsilon,2}}&\le C\|a^{-1/2}f\|_{L^2(0,1)}\bigg(\bigg\|\underbrace{\int_{L_\epsilon< \abs{\tau}\le 1}a^*(\cdot)^{1/2}K_{\rm Y}(\tau,\cdot)f^*(\tau)\,d\tau}_{T_{1,\epsilon}[a^{-1/2}f]}\bigg\|_{L^2(-L_\epsilon,L_\epsilon)}\\
&\qquad +\bigg\|\underbrace{\int_{L_\epsilon< \abs{t}\le 1}a^*(\cdot)^{1/2}K_{\rm Y}(\cdot,t)f^*(t)\,dt }_{T_{2,\epsilon}[a^{-1/2}f]}\bigg\|_{L^2(-L_\epsilon,L_\epsilon)} \bigg)\,.
\end{aligned}
\end{equation}
We will again rely on Schur's test \eqref{eq:schur1}-\eqref{eq:schur2} to estimate $T_{1,\epsilon}[a^{-1/2}f]$ and $T_{2,\epsilon}[a^{-1/2}f]$.

Using \eqref{eq:Klower} and that $1-L_\epsilon=\ell(\epsilon)\le C\epsilon^2$, we have
\begin{equation}
\begin{aligned}
A_{1,\epsilon} &:= \sup_{t\in[-L_\epsilon,L_\epsilon]}\int_{L_\epsilon< \abs{\tau}\le 1}\abs{K_{\rm Y}(\tau,t)}\,a^*(\tau)^{1/2}a^*(t)^{1/2}\,d\tau\\
&\le C\epsilon^{1/2}\sup_{t\in[-L_\epsilon,L_\epsilon]}\int_{L_\epsilon< \abs{\tau}\le 1}\abs{K_{\rm Y}(\tau,t)}\,d\tau
\le C\epsilon^{1/2}\int_{L_\epsilon< \abs{\tau}\le 1}\frac{1}{\epsilon^2}\,d\tau
\le C\epsilon^{1/2}\,.
\end{aligned}
\end{equation}
By an analogous computation, we also have
\begin{equation}
\begin{aligned}
A_{2,\epsilon} &:= \sup_{\tau\in[-L_\epsilon,L_\epsilon]}\int_{L_\epsilon< \abs{t}\le 1}\abs{K_{\rm Y}(\tau,t)}\,a^*(\tau)^{1/2}a^*(t)^{1/2}\,dt 
\le C\epsilon^{1/2}\,.
\end{aligned}
\end{equation}
Furthermore, using the form \eqref{eq:KY} of $K_{\rm Y}(\tau,t)$ and the definition \eqref{eq:YReps} of $\bR_\epsilon$, we may bound 
\begin{equation}
\begin{aligned}
B_{1,\epsilon} &:= \sup_{L_\epsilon< \abs{\tau}\le 1}\int_{-L_\epsilon}^{L_\epsilon}\abs{K_{\rm Y}(\tau,t)}\,a^*(\tau)^{1/2}a^*(t)^{1/2}\,dt\\
&\le C\epsilon^{1/2}\sup_{L_\epsilon< \abs{\tau}\le 1}\int_0^{2\pi}\int_{-L_\epsilon}^{L_\epsilon}\frac{1}{\abs{\bR_\epsilon(\tau,t,\theta)}}\,dt\,d\theta
\le C\epsilon^{1/2}\abs{\log\epsilon}\,,
\end{aligned}
\end{equation}
by Lemma \ref{lem:int1}. By a similar computation, but using \eqref{eq:Reps_UL} and the simple calculus estimate \eqref{eq:calculus} in place of Lemma \ref{lem:int1}, we have 
\begin{equation}
\begin{aligned}
B_{2,\epsilon} &= \sup_{L_\epsilon< \abs{t}\le 1}\int_{-L_\epsilon}^{L_\epsilon}\abs{K_{\rm Y}(\tau,t)}\,a^*(\tau)^{1/2}a^*(t)^{1/2}\,d\tau\\
&\le \epsilon^{1/2}\sup_{L_\epsilon< \abs{t}\le 1}\int_{-L_\epsilon-\varphi_\epsilon^{-1}(t)}^{L_\epsilon-\varphi_\epsilon^{-1}(t)}\frac{1}{\sqrt{\bars^2+\epsilon^4}}\,d\tau
\le C\epsilon^{1/2}\abs{\log\epsilon}\,.
\end{aligned}
\end{equation}
Thus, by Schur's test, we may bound
\begin{equation}
\begin{aligned}
\norm{T_{j,\epsilon}[a^{-1/2}f]}_{L^2(-L_\epsilon,L_\epsilon)} &\le CA_{j,\epsilon}^{1/2}B_{j,\epsilon}^{1/2}\|a^{-1/2}f\|_{L^2(0,1)}
\le C\epsilon^{1/2}\abs{\log\epsilon}^{1/2}\|a^{-1/2}f\|_{L^2(0,1)}
\end{aligned}
\end{equation}
for both $j=1,2$. Therefore, we obtain
\begin{equation}
\abs{I_{\epsilon,2}} \le C\epsilon^{1/2}\abs{\log\epsilon}^{1/2}\|a^{-1/2}f\|_{L^2(0,1)}^2\,.
\end{equation}

Combining the positivity of $I_{\epsilon,0}$ with the smallness of $I_{\epsilon,1}$, $I_{\epsilon,2}$, and $I_{\epsilon,3}$, we obtain \eqref{eq:I_toshow}, which yields Lemma \ref{lem:pos_op}.
\hfill\qedsymbol


\subsubsection{Proof of Lemma \ref{lem:sym_op}: near-antisymmetry of kernel derivatives}\label{subsubsec:sym_op}
Given the kernel $K_\epsilon(s,t)$ as in \eqref{eq:Keps}, we begin by writing
\begin{equation}
\p_sK_\epsilon(s,t) = -\p_tK_\epsilon(s,t) - K^{\rm r}_\epsilon(s,t)\,, \qquad K^{\rm r}_\epsilon(s,t):=\frac{\eta}{8\pi^2}\int_0^{2\pi}\frac{\bR_\epsilon\cdot(\p_s\bR_\epsilon+\p_t\bR_\epsilon)}{\abs{\bR_\epsilon}^3} \,d\theta\,.
\end{equation}
Using the expression \eqref{eq:Reps_def} for $\bR_\epsilon$, we note that
\begin{equation}\label{eq:dsR_dtR}
\begin{aligned}
\p_s\bR_\epsilon = (1-\epsilon a(s)\wh\kappa)\be_{\rm t}(s) + \epsilon a'(s)\be_r(s,\theta)  \,, \qquad
\p_t\bR_\epsilon = -(1-\ell(\epsilon))\be_{\rm t}(\varphi_\epsilon^{-1}(t))
\end{aligned}
\end{equation}
for $\wh\kappa(s,\theta)$ as in \eqref{eq:free_jacfac} and $\ell(\epsilon)$ as in \eqref{eq:Xeff}. Then, using the notation $\bars=s-\varphi_\epsilon^{-1}(t)$ as in \eqref{eq:bars}, we may write
\begin{equation}
 \be_{\rm t}(\varphi_\epsilon^{-1}(t)) = \be_{\rm t}(s) -\bars\bm{Q}(s,\bars)\,,
 \end{equation} 
where we note that $\be_{\rm t}(s)\cdot\bm{Q}(s,\bars) = \be_{\rm t}(s)\cdot\big(\kappa_1(\xi)\be_1(s) +\kappa_2(\xi)\be_2(s) + (\xi-s)\wt{\bm{Q}}(s,\xi)\big)$ for some $\xi\in[s,t]$; in particular, $\be_{\rm t}(s)\cdot\bm{Q}(s,\bars)=\bars Q_{\rm t}(s,\bars)$ for some $\abs{Q_{\rm t}}\le C\kappa_\star$.
Therefore, expanding
\begin{equation}
\begin{aligned}
\bR_\epsilon\cdot(\p_s\bR_\epsilon + \p_t\bR_\epsilon) &= (\bars\be_{\rm t}+\epsilon a \be_r + \bars^2\bm{Q})\cdot\big(\bars(1-\ell(\epsilon))\bm{Q} +(\ell(\epsilon) - \epsilon a\wh\kappa)\be_{\rm t} +\epsilon a'\be_r\big) \,,
\end{aligned}
\end{equation}
we may write, in the notation of \eqref{eq:SN_Gamma},
\begin{equation}\label{eq:Hepsdef}
\begin{aligned}
H_\epsilon[f]:=\int_{-1}^1K^{\rm r}_\epsilon(s,t)\, f^*(t)\,dt 
&= \frac{1}{1-\ell(\epsilon)}\int_{s-\varphi_\epsilon^{-1}(1)}^{s+\varphi_\epsilon^{-1}(1)}K^{\rm r}_\epsilon(s,\bars)\, f^*(\varphi_\epsilon(s-\bars))\,d\bars \\
&= \frac{\eta}{8\pi^2(1-\ell(\epsilon))}\big(J_1(s) + J_2(s) + J_3(s) + J_4(s) \big)\,,
\end{aligned}
\end{equation}
where 
\begin{equation}\label{eq:Js}
\begin{aligned}
J_1(s) &= \int_{s-\varphi_\epsilon^{-1}(1)}^{s+\varphi_\epsilon^{-1}(1)}\int_0^{2\pi}\frac{\epsilon^2a(s)a'(s)}{\abs{\bR_\epsilon}^3} \,d\theta\, f^*(\varphi_\epsilon(s-\bars))\,d\bars\,,\\
J_2(s) &= \int_{s-\varphi_\epsilon^{-1}(1)}^{s+\varphi_\epsilon^{-1}(1)}\int_0^{2\pi}\frac{\bars\ell(\epsilon)}{\abs{\bR_\epsilon}^3} \,d\theta\, f^*(\varphi_\epsilon(s-\bars))\,d\bars \,,\\
J_3(s) &= \int_{s-\varphi_\epsilon^{-1}(1)}^{s+\varphi_\epsilon^{-1}(1)}\int_0^{2\pi}\frac{\bars^3\big((1-\epsilon a\wh\kappa)Q_{\rm t}+ (1-\ell(\epsilon))\abs{\bm{Q}}^2\big)}{\abs{\bR_\epsilon}^3} \,d\theta\, f^*(\varphi_\epsilon(s-\bars))\,d\bars \,,\\
J_4(s) &= \int_{s-\varphi_\epsilon^{-1}(1)}^{s+\varphi_\epsilon^{-1}(1)}\int_0^{2\pi}\frac{\big
(\epsilon a\bars(1-\ell(\epsilon)) + \bars^2\epsilon a'\big)\bm{Q}\cdot\be_r(s,\theta) - \bars\epsilon a \wh\kappa(s,\theta)}{\abs{\bR_\epsilon}^3} \,d\theta\, f^*(\varphi_\epsilon(s-\bars))\,d\bars\,.
\end{aligned}
\end{equation}

We may estimate each $J_i$ in turn, beginning with $J_4$. Using that $\wh\kappa(s,\theta)=\cos\theta\kappa_1(s) + \sin\theta\kappa_2(s)$ and $\be_r(s,\theta)=\cos\theta\be_1(s)+\sin\theta\be_2(s)$, we may use Lemma \ref{lem:int4} to obtain the bound 
\begin{equation}
\abs{J_4} \le C\big(\epsilon a(s) + \epsilon^2\abs{\log\epsilon}|a'(s)| \big) \norm{f}_{L^\infty}\,.
\end{equation}
For $J_3$, letting
\begin{equation}
K_3(s,\bars) = \int_0^{2\pi}\frac{\bars^3\big((1-\epsilon a\wh\kappa)Q_{\rm t}+ (1-\ell(\epsilon))\abs{\bm{Q}}^2\big)}{\abs{\bR_\epsilon}^3} \,d\theta\,,
\end{equation}
we note that, by \eqref{eq:Reps_UL}, $\abs{K_3}\le C$ for $C$ independent of $\epsilon$. We thus have
\begin{equation}
\abs{J_3(s)} \le \norm{K_3(s,\cdot)}_{L^2}\norm{f}_{L^2} \le C\norm{f}_{L^2}\,.
\end{equation}
For $J_2$, we write the kernel as 
\begin{equation}
K_2(s,\bars) = \int_0^{2\pi}\frac{\bars\ell(\epsilon)}{\abs{\bR_\epsilon}^3} \,d\theta\,, \qquad \abs{K_2}\le \int_0^{2\pi}\frac{C\epsilon^2}{\abs{\bR_\epsilon}^2}\,d\theta \le C\epsilon a^{-1}(s)\int_0^{2\pi}\frac{1}{\abs{\bR_\epsilon}}\,d\theta \,.
\end{equation}
Using Lemma \ref{lem:at_bound}, we then note that
\begin{equation}
\norm{a^*(\cdot)^{1/2}K_2(s,\cdot)}_{L^2}^2 \le C\epsilon a^{-2}(s)\,,
\end{equation}
and therefore
\begin{equation}
\abs{J_2(s)} \le C\norm{a^{1/2}(\cdot)K_2(s,\cdot)}_{L^2}\|a^{-1/2}f\|_{L^2} \le C\epsilon^{1/2}a^{-1}(s)\|a^{-1/2}f\|_{L^2}\,.
\end{equation}

Finally, to estimate $J_1$, we begin by noting that for $s\in [0,1]$ and $t\in[-1,1]$, the radius function $a(s)$ admits the expansion 
\begin{equation}\label{eq:a_expand}
a(s) = a^*(t) + (s-t)\mc{R}^{(a)}(s,t)\,, 
\end{equation}
where, using the spheroidal endpoint condition \eqref{eq:spheroidal}, for $s,t>1-\delta$ we have
\begin{equation}\label{eq:Ra_bd}
 |\mc{R}^{(a)}(s,t)| \le C \abs{\frac{t+s}{\sqrt{1-s^2}+\sqrt{1-t^2}}} \le Ca^*(t)^{-1}\,.
 \end{equation} 
 In particular, 
\begin{equation}
\abs{\mc{R}^{(a)}(s,t)f^*(t)} \le C\norm{a^{-1}f}_{L^\infty}\,.
\end{equation}
Using $t$ and $s$ rather than $\bars$ and $s$, we may thus write
\begin{equation}
\frac{1}{1-\ell(\epsilon)}J_1(s) = \underbrace{\int_{-1}^{1}\int_0^{2\pi}\frac{\epsilon^2a'(s)}{\abs{\bR_\epsilon}^3} \,d\theta\, a^*(t)f^*(t)\, dt }_{J_{1a}(s)}
+ \underbrace{\int_{-1}^{1}\int_0^{2\pi}\frac{\epsilon^2a'(s)(s-t)}{\abs{\bR_\epsilon}^3} \,d\theta\, \mc{R}^{(a)}f^*(t)\, dt}_{J_{1b}(s)} \,.
\end{equation}
Using that $\abs{s-t}=\abs{\bars+\varphi_\epsilon^{-1}(t)-t}\le \abs{\bars}+\epsilon^2$, we may bound $a^2(s)J_{1b}(s)$ as
\begin{equation}
\begin{aligned}
\abs{a^2(s)J_{1b}(s)}&\le Ca_\star a(s)\norm{\mc{R}^{(a)}f}_{L^\infty}\int_0^{2\pi}\int_{s-\varphi_\epsilon^{-1}(1)}^{s+\varphi_\epsilon^{-1}(1)}\frac{\epsilon^2\abs{\bars}+\epsilon^4}{\abs{\bR_\epsilon}^3}\, d\bars\,d\theta 
\le C \epsilon \norm{a^{-1}f}_{L^\infty}\,,
\end{aligned}
\end{equation}
where we have used Lemma \ref{lem:int1}. 

To estimate $J_{1a}$, we first note that by Lemma \ref{lem:int3}, we may write
\begin{equation}
a^2(s)J_{1a}(s) = 4\pi a(s)a'(s) f(s) + \mc{R}_J[f](s)\,,
\end{equation}
where
\begin{equation}
\abs{\mc{R}_J[f](s)} \le Ca_\star\epsilon\big(\norm{a\,f_s}_{L^\infty}+\norm{a^{-1}f}_{L^\infty} \big) \,.
\end{equation}
We may then estimate: 
\begin{equation}
\begin{aligned}
\norm{a^2J_{1a}}_{L^2} &\le C\norm{f}_{L^2} + C\norm{\mc{R}_J[f]}_{L^\infty}
\le C\big(\norm{f}_{L^2} + \epsilon\norm{a\,f_s}_{L^\infty} + \epsilon\norm{a^{-1}f}_{L^\infty}\big)\,,\\
\norm{a^2J_{1a}}_{L^\infty} &\le C\norm{f}_{L^\infty} + C\norm{\mc{R}_J[f]}_{L^\infty}
\le C\big(\norm{a^{-1}f}_{L^\infty} + \epsilon\norm{a\,f_s}_{L^\infty} \big)\,.
\end{aligned}
\end{equation}

Collecting the above estimates and returning to the expression \eqref{eq:Hepsdef}, we may thus bound the remainder term $H_\epsilon[f]$ as
\begin{equation}\label{eq:Heps_ests}
\begin{aligned}
\|a^2H_\epsilon[f]\|_{L^2} &\le C\big(\norm{a^2J_1}_{L^2} + \norm{a^2J_2}_{L^\infty}+\norm{a^2J_3}_{L^\infty}+\norm{a^2J_4}_{L^\infty}\big)\\
&\le C\big(\|a^{-1/2}f\|_{L^2} + \epsilon\norm{a\,f_s}_{L^\infty} + \epsilon\norm{a^{-1}f}_{L^\infty}\big)\,, \\
\norm{a^2\,H_\epsilon[f]}_{L^\infty} &\le C\big(\norm{a^2J_1}_{L^\infty} + \norm{a^2J_2}_{L^\infty}+\norm{a^2J_3}_{L^\infty}+\norm{a^2J_4}_{L^\infty}\big)\\
&\le C\big(\norm{a^{-1}f}_{L^\infty} + \epsilon\norm{a\,f_s}_{L^\infty} \big)\,.
\end{aligned}
\end{equation}
\hfill \qedsymbol

\subsubsection{Proof of Lemma \ref{lem:kernelIBP}: weighted integration by parts}\label{subsubsec:IBP}
We begin by using the expansion \eqref{eq:a_expand} for $a(s)$ to write 
\begin{equation}
\begin{aligned}
a^2(s)\int_{-1}^1 \p_tK_\epsilon(s,t)\, f^*(t)\,dt &= 
\underbrace{\int_{-1}^1a^*(t)^2\p_tK_\epsilon(s,t)\, f^*(t)\,dt}_{G_1(s)} + 
\underbrace{\int_{-1}^1K^{(a)}(s,t)\, f^*(t)\,dt}_{G_2(s)} \,, \\
K^{(a)}(s,t) &= (s-t)\wt{\mc{R}}^{(a)}(s,t)\p_tK_\epsilon(s,t)\,,
\end{aligned}
\end{equation}
where $|\wt{\mc{R}}^{(a)}(s,t)|\le Ca_\star$. Using the expression \eqref{eq:dsR_dtR} for $\p_t\bR_\epsilon$ and that $\abs{s-t}\le \abs{s-\varphi_\epsilon^{-1}(t)}+C\epsilon^2$, we have 
\begin{equation}\label{eq:Kabd}
\abs{K^{(a)}(s,t)} \le C\int_0^{2\pi}\frac{\abs{s-\varphi_\epsilon^{-1}(t)}+\epsilon^2}{\abs{\bR_\epsilon}^2}\,d\theta \le \int_0^{2\pi}\frac{C}{\abs{\bR_\epsilon}}\,d\theta + \int_0^{2\pi}\frac{C\epsilon^2}{\abs{\bR_\epsilon}^2}\,d\theta\,.
\end{equation}

We will again rely on Schur's test \eqref{eq:schur1}-\eqref{eq:schur2} to obtain an $L^2$ estimate for $a^{1/2}(s)G_2(s)$.
First, we may estimate 
\begin{equation}\label{eq:Acomps}
\begin{aligned}
A_{g,1}&:=\sup_{s\in[0,1]}\int_{-1}^1\frac{a(s)^{1/2}a^*(t)^{1/2}}{\abs{\bR_\epsilon}}\,dt  \le C\abs{\log\epsilon}\,,\\
A_{g,2}&:= \epsilon^2\sup_{s\in[0,1]} \int_{-1}^1\frac{a(s)^{1/2}a^*(t)^{1/2}}{\abs{\bR_\epsilon}^2}\,dt\\
&\le \epsilon^2\sup_{s\in[0,1]}\bigg( \int_{-1}^1\frac{a(s)}{\abs{\bR_\epsilon}^2}\,dt\bigg)^{1/2}\bigg( \int_{-1}^1\frac{a^*(t)}{\abs{\bR_\epsilon}^2}\,dt\bigg)^{1/2}
\le C\epsilon\,,
\end{aligned}
\end{equation}
by Lemmas \ref{lem:int1} and \ref{lem:at_bound}. Furthermore, using the lower bound \eqref{eq:Reps_UL} for $\abs{\bR_\epsilon}$ in terms of $\bars=s-\varphi_\epsilon^{-1}(t)$ and that $\abs{\bR_\epsilon}\ge C\epsilon^2$, we may estimate
\begin{equation}\label{eq:Bcomps}
\begin{aligned}
B_{g,1}&:=\sup_{t\in[-1,1]}\int_0^1\frac{a(s)^{1/2}a^*(t)^{1/2}}{\abs{\bR_\epsilon}}\,ds \le \sup_{t\in[-1,1]}\bigg(\int_0^{1-\ell(\epsilon)}\frac{1}{\abs{\bR_\epsilon}}\,ds+\int_{1-\ell(\epsilon)}^1\frac{a(s)^{1/2}}{\abs{\bR_\epsilon}}\,ds \bigg)\\ 
&\le  \sup_{t\in[-1,1]}\bigg(\int_{-\varphi_\epsilon^{-1}(t)}^{\varphi_\epsilon^{-1}(1-t)}\frac{C}{\sqrt{\bars^2+\epsilon^4}}\,d\bars +\int_{1-\ell(\epsilon)}^1\frac{C\epsilon^{1/2}}{\abs{\bR_\epsilon}}\,ds\bigg)\\
&\le C\big(\abs{\log\epsilon} + \ell(\epsilon)\epsilon^{1/2}\epsilon^{-2}\big)  \le C\abs{\log\epsilon}\,,\\
B_{g,2}&:=\sup_{t\in[-1,1]}\int_0^1\frac{\epsilon^2a(s)^{1/2}a^*(t)^{1/2}}{\abs{\bR_\epsilon}^2}\,ds 
\le \sup_{t\in[-1,1]}\bigg(\int_0^{1-\ell(\epsilon)}\frac{\epsilon^2}{\abs{\bR_\epsilon}^2}\,ds + \int_{1-\ell(\epsilon)}^1\frac{\epsilon^2a(s)^{1/2}}{\abs{\bR_\epsilon}^2}\,ds \bigg)\\
&\le 
\sup_{t\in[-1,1]}\bigg(\int_{-\varphi_\epsilon^{-1}(t)}^{\varphi_\epsilon^{-1}(1-t)}\frac{\epsilon^2}{\bars^2+\epsilon^4}\,d\bars + \int_{1-\ell(\epsilon)}^1\frac{\epsilon^{5/2}}{\abs{\bR_\epsilon}^2}\,ds \bigg)\\
&\le 
 C\big(\epsilon^2(\epsilon^2)^{-1}+\ell(\epsilon)\epsilon^{5/2}\epsilon^{-4}\big)
 \le C\,.
\end{aligned}
\end{equation}
Here we have used the simple calculus bound \eqref{eq:calculus} as well as the fact that $a(s)\ge C\epsilon$ for $0\le s\le 1-\ell(\epsilon)$ and $a(s)\le C\epsilon$ for $s\ge 1-\ell(\epsilon)$. We also recall that $\ell(\epsilon)=O(\epsilon^2)$.

We may then use Schur's test \eqref{eq:schur2} to estimate $a^{1/2}(s)G_2(s)$ in $L^2$ as
\begin{equation}\label{eq:G2L2}
\begin{aligned}
\norm{a^{1/2}G_2}_{L^2(0,1)} &= \norm{a^{1/2}(s)\int_{-1}^1K^{(a)}(\cdot,t)a^*(t)^{1/2}\,a^*(t)^{-1/2}f^*(t)\,dt}_{L^2(0,1)} \\
&\le C(A_{g,1}^{1/2}B_{g,1}^{1/2}+A_{g,2}^{1/2}B_{g,2}^{1/2})\|a^{1/2}f\|_{L^2(0,1)} 
\le C\abs{\log\epsilon}\|a^{1/2}f\|_{L^2(0,1)} \,.
\end{aligned}
\end{equation}

To estimate $G_2$ in $L^\infty$, we may use the bound \eqref{eq:Kabd} for $K^{(a)}(s,t)$ to obtain
\begin{equation}\label{eq:G2}
\begin{aligned}
\abs{G_2(s)}&\le \int_{-1}^1\abs{K^{(a)}(s,t)}\abs{f^*(t)}\,dt \\
&\le C\int_0^{2\pi}\bigg(\int_{-1}^1\frac{a^*(t)}{\abs{\bR_\epsilon}^2}\,dt \bigg)^{1/2}\bigg(\int_{-1}^1a^*(t)^{-1}f^*(t)^2\,dt \bigg)^{1/2}\,d\theta
+ C\epsilon\norm{a^{-1}f}_{L^\infty} \\
&\le C\epsilon^{-1/2}\|a^{-1/2}f\|_{L^2(0,1)} + C\epsilon\norm{a^{-1}f}_{L^\infty(0,1)}\,, 
\end{aligned}
\end{equation}
where we have used Lemma \ref{lem:at_bound} to bound both terms.

For $G_1(s)$, we may use that $[a^{-1}f]^*\in L^\infty$, so $f^*\big|_{t=\pm 1}=0$, to integrate by parts\footnote{Note that $[af]^*$ may not be $C^1$ at $t=0$ due to reflection across 0, but $[af]^*\in C[-1,1]$ and $a(s)f(s)\in C^1[0,1]$, so $[af]^*\in H^1(-1,1)$.} in $t$. We obtain
\begin{equation}
G_1(s) = -\int_{-1}^1K_\epsilon(s,t)\,a^*(t)^2f_t^*(t)\,dt -\underbrace{2\int_{-1}^1K_\epsilon(s,t)\, a^*a^*_t(t)f^*(t)\,dt}_{G_{1a}(s)}\,.
\end{equation}
Using that, as in \eqref{eq:Acomps}, we may bound 
\begin{equation}
\begin{aligned}
\sup_{s\in [0,1]}\int_{-1}^1\abs{K_\epsilon(s,t)}\,dt \le \sup_{s\in [0,1]}\int_{-1}^1\frac{C}{\abs{\bR_\epsilon}}\,dt \le C\abs{\log\epsilon}\,,
\end{aligned}
\end{equation}
and as in \eqref{eq:Bcomps}, we may bound
\begin{equation}
\begin{aligned}
\sup_{t\in [-1,1]}\int_{0}^1\abs{K_\epsilon(s,t)}\,ds&\le C\sup_{t\in [-1,1]}\bigg(\int_{-\varphi_\epsilon^{-1}(t)}^{\varphi_\epsilon^{-1}(1-t)}\frac{1}{\sqrt{\bars^2+\epsilon^4}}\,d\bars
+ \int_{1-\ell(\epsilon)}^1\frac{1}{\abs{\bR_\epsilon}}\,ds \bigg)\\
&\le C\abs{\log\epsilon} + C\ell(\epsilon)\epsilon^{-2} \le C\abs{\log\epsilon}\,,
\end{aligned}
\end{equation}
by Schur's test \eqref{eq:schur2}, we have
\begin{equation}\label{eq:G1aL2}
\norm{G_{1a}}_{L^2}\le Ca_\star\abs{\log\epsilon}\norm{f}_{L^2(0,1)}\,.
\end{equation}
In addition, using Lemma \ref{lem:at_bound}, we have that $G_{1a}(s)$ satisfies the pointwise bound
\begin{equation}\label{eq:G1a}
\begin{aligned}
\abs{G_{1a}(s)} &\le C\bigg(\int_{-1}^1 a^*(t)K_\epsilon^2(s,t)\,dt\bigg)^{1/2} \bigg(\int_{-1}^1a^*(t)^{-1}f^*(t)^2\,dt\bigg)^{1/2} 
\le C\epsilon^{-1/2}\|a^{-1/2}f\|_{L^2} \,.
\end{aligned}
\end{equation}

Combining the $L^2$ bounds \eqref{eq:G2L2} and \eqref{eq:G1aL2}, we obtain \eqref{eq:GepsL2}, while the $L^\infty$ bounds \eqref{eq:G2} and \eqref{eq:G1a} may be combined to yield \eqref{eq:GepsLinfty}.
\hfill \qedsymbol

\subsection{Proof of Theorem \ref{thm:pSB}}\label{subsec:thm_pSB}
Equipped with the kernel estimates of Lemmas \ref{lem:pos_op}, \ref{lem:sym_op} and \ref{lem:kernelIBP}, we may proceed to the proof of Theorem \ref{thm:pSB} regarding the integrodifferential equation \eqref{eq:pSB} for $p^{\rm SB}(s)$.

Recalling the definition \eqref{eq:Keps} of $K_\epsilon(s,t)$ and letting $v=p^{\rm SB}-p_0$ and $\alpha(s)=2\pi a(s)\omega/\eta$, we have that $v$ satisfies the integrodifferential equation 
\begin{equation}\label{eq:v_eqn}
\begin{aligned}
(a^4(s)v_s)_s - \alpha(s)v + \alpha(s)\int_{-1}^1 K_\epsilon(s,t)\, [(a^4(t)v_t)_t]^*\,dt &= \alpha(s)p_0\,, \\ 
v(0) &= 0\,.
\end{aligned}
\end{equation}
Here we recall that the notation $[(a^4(t)v_t)_t]^*$, $t\in[-1,1]$, denotes the extension-by-reflection \eqref{eq:fstar} of $[(a^4(t)v_t)_t]$, $t\in[0,1]$, across $t=0$.

Recalling the definition \eqref{eq:Ha_def} of the weighted space $\mc{H}^a$, we further define $\mc{H}^a_0(0,1)= \{ u\in \mc{H}^a(0,1)\,:\, u\big|_{s=0}=0 \}$. Since $a(s)\sim\sqrt{1-s^2}$ as $s\to 1$, the following weighted Poincar\'e inequality holds on $\mc{H}^a_0$:
\begin{equation}\label{eq:Ha_poincare}
  \norm{u}_{L^2(0,1)} \le C\norm{a^2u_s}_{L^2(0,1)} \qquad \text{for all }u\in \mc{H}^a_0(0,1)\,.
\end{equation}
This is a version of Hardy's inequality ``at $\infty$'' (see, for example, \cite[Lemma 1]{mironescu2018role} or \cite[Lemma 3.14]{stein1971introduction}).
In particular, $\norm{a^2\p_s\cdot}_{L^2(0,1)}$ is a norm on $\mc{H}^a_0$.
Let $\mc{A}$ denote the weighted $H^2$-type function space 
\begin{equation}
\mc{A} = \big\{v\in \mc{H}^a_0(0,1)\;:\; a(s)^{-1/2}(a^4(s)v_s)_s\in L^2(0,1)\big\}
\end{equation}
with norm 
\begin{equation}
  \norm{v}_{\mc{A}} := \norm{a^{-1/2}(a^4v_s)_s}_{L^2(0,1)} + \norm{a^2v_s}_{L^2(0,1)}\,.
\end{equation}
We define the bilinear form $\mc{B}:\A\times\A\to \R$ by
\begin{equation}
\begin{aligned}
\mc{B}(v,w)&= \int_0^1\bigg( \frac{1}{\alpha(s)}(a^4v_s)_s(a^4w_s)_s+ a^4v_sw_s\bigg)\,ds\\
&\qquad + \int_0^1\int_{-1}^1 K_\epsilon(s,t)\, [(a^4v_t)_t]^*(t)(a^4w_s)_s(s)\,dt\,ds \,.
\end{aligned}
\end{equation}
Noting that $\abs{K_\epsilon}\le C\epsilon^{-2}$ and $\alpha(s)\propto a(s)$, we have
\begin{equation}
\begin{aligned}
\abs{\mc{B}(v,w)}&\le C\big(\|a^{-1/2}(a^4v_s)_s\|_{L^2}\|a^{-1/2}(a^4w_s)_s\|_{L^2} + \norm{a^2v_s}_{L^2}\norm{a^2w_s}_{L^2} \\
&\qquad + \epsilon^{-2}\norm{(a^4v_s)_s}_{L^2}\norm{(a^4w_s)_s}_{L^2} \big) 
\le C(\epsilon)\norm{v}_{\mc{A}}\norm{w}_{\mc{A}}\,.
\end{aligned}
\end{equation}
Furthermore, using Lemma \ref{lem:pos_op}, we have
\begin{equation}\label{eq:coercive}
\begin{aligned}
\mc{B}(v,v) &= \int_0^1\bigg( \frac{1}{\alpha(s)}\abs{(a^4v_s)_s}^2+ \abs{a^2v_s}^2\bigg)\,ds + \int_0^1\int_{-1}^1 K_\epsilon(s,t)\, [(a^4v_t)_t]^*(t)\,(a^4v_s)_s(s)\,dt\,ds \\
&\ge \frac{\eta}{2\pi\omega}\|a^{-1/2}(a^4v_s)_s\|_{L^2}^2 + \norm{a^2v_s}_{L^2}^2 - C\epsilon^{1/2}\abs{\log\epsilon}^{1/2}\|a^{-1/2}(a^4v_s)_s\|_{L^2}^2 \\
&\ge C\norm{v}_{\mc{A}}^2
\end{aligned}
\end{equation}
for $\epsilon$ sufficiently small.
We may then define a solution to \eqref{eq:v_eqn} as $v\in \mc{A}$ satisfying
\begin{equation}
\begin{aligned}
\mc{B}(v,w)
= p_0\int_{0}^1(a^4w_s)_s\,ds 
= -p_0(a^4w_s)\big|_{s=0}
\qquad \text{for all }w\in \mc{A}\,. 
\end{aligned}
\end{equation}
Existence and uniqueness follow from the Lax-Milgram lemma. In addition, we obtain the weighted $H^2$-type estimate 
\begin{equation}\label{eq:vssL2}
\norm{v}_{\mc{A}}= \|a^{-1/2}(a^4v_s)_s\|_{L^2} + \norm{a^2v_s}_{L^2} \le C\abs{p_0}
\end{equation}
for $C$ independent of $\epsilon$.

Now, rewriting equation \eqref{eq:v_eqn} as
\begin{equation}\label{eq:rewrite_v}
\begin{aligned}
(a^4v_s)_s - \alpha(s)v &= - \alpha(s)\int_{-1}^1 K_\epsilon(s,t)\, [(a^4v_t)_t]^*(t)\,dt +\alpha(s)p_0\,, 
\end{aligned}
\end{equation}
we note that the right hand side is of the form $\alpha(s)f(s)$ for $f\in L^\infty(0,1)$. In particular, using Lemma \ref{lem:at_bound}, we may estimate 
\begin{equation}
\begin{aligned}
\abs{\int_{-1}^1 K_\epsilon(s,t)\,[(a^4v_t)_t]^*(t)\,dt}
&\le 
C\bigg(\int_{-1}^1a^*(t)\abs{K_\epsilon(s,t)}^2\,dt\bigg)^{1/2}
\bigg(\int_{-1}^1 \frac{\big([(a^4v_t)_t]^*\big)^2(t)}{a^*(t)}\,dt \bigg)^{1/2}\\
&\le C\bigg(\int_{-1}^1\frac{a^*(t)}{\abs{\bR_\epsilon}^2}\,dt\bigg)^{1/2}\abs{p_0}
\le C\epsilon^{-1/2}\abs{p_0}\,.
\end{aligned}
\end{equation}
Given $f\in L^\infty(0,1)$, the equation
\begin{equation}\label{eq:max_principle}
(a^4v_s)_s - \alpha(s)v = \alpha(s)f(s) 
\end{equation}
satisfies a maximum principle: $\norm{v}_{L^\infty(0,1)}\le 2\norm{f}_{L^\infty(0,1)}$. This may be seen by the following energy argument. The solution $v\in \mc{H}^a_0$ to \eqref{eq:max_principle} is the unique minimizer of the energy
\begin{equation}
  E[v]=\frac{1}{2}\int_0^1\big(a^4v_s^2 + \alpha\,v^2\big)\,ds - \int_0^1\alpha\,f(s)v(s)\,ds\,,
\end{equation}
since, for any $u\in \mc{H}^a_0$,
\begin{equation}
\begin{aligned}
  E[v+u] &= \frac{1}{2}\int_0^1\big(a^4(v_s^2+2u_sv_s+u_s^2) + \alpha\,(v^2+2uv+u^2)\big)\,ds - \int_0^1\alpha\,f(s)(v(s)+u(s))\,ds  \\
  &= E[v] + \frac{1}{2}\int_0^1\big(a^4 u_s^2 + \alpha \,u^2\big)\,ds \ge E[v]\,,
\end{aligned}
\end{equation}
with equality only if $u=0$. Assume that $\norm{f}_{L^\infty}$ is achieved for $f>0$ (otherwise reverse the following signs accordingly). For the sake of contradiction, suppose $v(s)>2\norm{f}_{L^\infty}$ for $s\in I\subset [0,1]$. Writing the contribution of interval $I$ to the energy $E$ as 
\begin{equation}
  E_I = \frac{1}{2}\int_Ia^4v_s^2\,ds + \int_I\alpha(s) v(s)\bigg(\frac{1}{2}v(s) - f(s) \bigg)\,ds\,,
\end{equation}
we see that $E_I$ may be decreased by replacing $v(s)$ with $2f(s)$, contradicting that $v$ minimizes $E$.

Returning to \eqref{eq:rewrite_v}, we thus have 
\begin{equation}\label{eq:first_bd}
\begin{aligned}
\abs{a(s)^{-1}(a^4v_s)_s} &\le C\big(\norm{v}_{L^\infty} + C\epsilon^{-1/2}\abs{p_0} + \abs{p_0} \big) 
\le C\epsilon^{-1/2}\abs{p_0} \,.
\end{aligned}
\end{equation}

Next, dividing \eqref{eq:v_eqn} through by $\alpha(s)$ and taking a derivative, we have that $(a^4v_s)_{ss}$ satisfies
\begin{equation}\label{eq:vsss_eqn}
\begin{aligned}
\frac{1}{\alpha}(a^4v_s)_{ss} - \frac{\alpha'}{\alpha^2}(a^4v_s)_s - v_s + \int_{-1}^1 \p_sK_\epsilon(s,t)\, [(a^4v_t)_t]^*(t)\,dt &= 0\,.
\end{aligned}
\end{equation}
%
%
%
To estimate the integral term, we may combine Lemmas \ref{lem:sym_op} and \ref{lem:kernelIBP} to write 
\begin{equation}
\begin{aligned}
a^2(s)\int_{-1}^1\p_sK_\epsilon(s,t)\,[(a^4v_t)_t]^*(t)\,dt &= \int_{-1}^1 a^*(t)^2K_\epsilon(s,t)\, [(a^4v_t)_{tt}]^*(t)\,dt \\
&\qquad+ a^2(s)H_\epsilon[(a^4v_s)_s](s) - G_\epsilon[(a^4v_s)_s](s)  \,,
\end{aligned}
\end{equation}
where, by \eqref{eq:vssL2} and \eqref{eq:first_bd}, we have 
\begin{equation}\label{eq:Heps_bds}
\begin{aligned}
\|a^2H_\epsilon[(a^4v_s)_s]\|_{L^2} 
&\le C\abs{p_0}+ C\epsilon\norm{a\,(a^4v_s)_{ss}}_{L^\infty} \,, \\
\norm{a^2H_\epsilon[(a^4v_s)_s]}_{L^\infty} 
&\le C\epsilon^{-1/2}\abs{p_0}+  C\epsilon\norm{a\,(a^4v_s)_{ss}}_{L^\infty} \,,
\end{aligned}
\end{equation}
as well as
\begin{equation}\label{eq:Geps_bds}
\begin{aligned}
\norm{a^{1/2}G_\epsilon[(a^4v_s)_s]}_{L^2}&\le C\abs{\log\epsilon}\abs{p_0}\,, \\
\norm{G_\epsilon[(a^4v_s)_s]}_{L^\infty} &\le C\epsilon^{-1/2}\abs{p_0} \,.
\end{aligned}
\end{equation}

Now, to obtain a weighted $L^2$ estimate for $(a^4v_s)_{ss}$, using that $\alpha^{1/2}\alpha'(a^4v_s)_s\in L^2(0,1)$ since $\abs{\alpha\alpha'}\le C\,a_\star$, we rewrite \eqref{eq:vsss_eqn} as 
\begin{equation}\label{eq:vsss_eqn2}
\begin{aligned}
\alpha^{3/2}(a^4v_s)_{ss}  + \alpha^{1/2}\int_{-1}^1 K_\epsilon(s,t)\,\alpha^*(t)^2 & [(a^4v_t)_{tt}]^*(t)\,dt =  \alpha^{1/2}\alpha'(a^4v_s)_s + \alpha^{5/2}v_s \\
&- \alpha^{5/2} H_\epsilon[(a^4v_s)_s] + \frac{4\pi^2\omega^2}{\eta^2}\alpha^{1/2}G_\epsilon[(a^4v_s)_s]  \,,
\end{aligned}
\end{equation}
where we recall that $\alpha(s)\propto a(s)$. 
Using Lemma \ref{lem:pos_op}, we have
\begin{equation}
\int_0^1\alpha^2(s)(a^4v_s)_{ss}\int_{-1}^1 K_\epsilon(s,t)\,\alpha^*(t)^2 [(a^4v_t)_{tt}]^*(t)\,dt\,ds \ge -C\epsilon^{1/2}\abs{\log\epsilon}^{1/2}\|a^{3/2}(a^4v_s)_{ss}\|_{L^2}^2\,.
\end{equation}
Thus, combining the $L^2$ bounds \eqref{eq:vssL2}, \eqref{eq:Heps_bds}, and \eqref{eq:Geps_bds}, for $\epsilon$ sufficiently small, we obtain 
\begin{equation}\label{eq:vsss_L2}
\|a^{3/2}(a^4v_s)_{ss}\|_{L^2} \le C\abs{\log\epsilon}\abs{p_0} + C\epsilon\norm{a\,(a^4v_s)_{ss}}_{L^\infty}\,.
\end{equation}

To close the $L^2$ bound, we also need a weighted $L^\infty$ estimate for $(a^4v_s)_{ss}$. We rewrite \eqref{eq:vsss_eqn2} as 
\begin{equation}\label{eq:vsss_eq3}
\begin{aligned}
\alpha (a^4v_s)_{ss} &- \alpha^2v_s + \int_{-1}^1 K_\epsilon(s,t)\,\alpha^*(t)^2 [(a^4v_t)_{tt}]^*(t)\,dt \\
&=  \alpha'(a^4v_s)_s  - \alpha^2 H_\epsilon[(a^4v_s)_s] + \frac{4\pi^2\omega^2}{\eta^2}G_\epsilon[(a^4v_s)_s]  \,,
\end{aligned}
\end{equation}
so that by \eqref{eq:first_bd} and \eqref{eq:astar}, each term on the right hand side belongs to $L^\infty$. We have that
\begin{equation}
\begin{aligned}
\abs{\int_{-1}^1 K_\epsilon(s,t)\,\alpha^*(t)^2 [(a^4v_t)_{tt}]^*(t)\,dt} &\le C\|a^{3/2}(a^4v_s)_{ss}\|_{L^2}\bigg(\int_{-1}^1 \abs{K_\epsilon}^2a^*(t)\,dt\bigg)^{1/2}\\
&\le C\epsilon^{-1/2}\|a^{3/2}(a^4v_s)_{ss}\|_{L^2}\,, 
\end{aligned}
\end{equation}
by Lemma \ref{lem:at_bound}. 
Furthermore, since $a^2v_s\in L^2$, by the spheroidal decay \eqref{eq:spheroidal} of $a(s)$ as $s\to 1$, we have that $\abs{a^2v_s}\lesssim a^{-1}$ almost everywhere, and thus $\abs{a^4v_s}\lesssim a(s)$ almost everywhere. In particular, $a^4v_s\to 0$ along some subsequence as $s\to 1$.
Then, since $\abs{(a^4v_s)_s}\le C\,a(s)\epsilon^{-1/2}\abs{p_0}$, upon integrating in $s$, we have  
\begin{equation}
  \abs{a^4v_s} \le C\epsilon^{-1/2}\abs{p_0}\abs{\int_s^1a(t)\,dt} \le C\epsilon^{-1/2}a^3(s)\abs{p_0}\,,
\end{equation}
by \eqref{eq:spheroidal}. Therefore, $\abs{av_s}\le C\epsilon^{-1/2}\abs{p_0}$.
Using \eqref{eq:first_bd}, \eqref{eq:Heps_bds} and \eqref{eq:Geps_bds}, we then have 
\begin{equation}
\begin{aligned}
\abs{\alpha (a^4v_s)_{ss}} &\le C\big(\epsilon^{-1/2}\|a^{3/2}(a^4v_s)_{ss}\|_{L^2}+ a_\star\norm{a^{-1}(a^4v_s)_s}_{L^\infty} + \abs{a^2v_s} \\
&\qquad + \norm{a^2 H_\epsilon[(a^4v_s)_s]}_{L^\infty}+\norm{G_\epsilon[(a^4v_s)_s]}_{L^\infty}\big) \\
&\le C\epsilon^{-1/2}\abs{p_0} + C\epsilon\norm{a\,(a^4v_s)_{ss}}_{L^\infty} + C\epsilon^{-1/2}\|a^{3/2}(a^4v_s)_{ss}\|_{L^2}\,.
\end{aligned}
\end{equation}
For $\epsilon$ sufficiently small, we may absorb the right hand side $\norm{a\,(a^4v_s)_{ss}}_{L^\infty}$ into the left hand side to obtain 
\begin{equation}\label{eq:vsss_Linfty}
\norm{a\,(a^4v_s)_{ss}}_{L^\infty} \le C\epsilon^{-1/2}\abs{p_0} + C\epsilon^{-1/2}\|a^{3/2}(a^4v_s)_{ss}\|_{L^2}\,.
\end{equation}
Then, combining the weighted $L^\infty$ bound \eqref{eq:vsss_Linfty} with the weighted $L^2$ bound \eqref{eq:vsss_L2}, we have
\begin{equation}\label{eq:vsss_L2_2}
\begin{aligned}
\|a^{3/2}(a^4v_s)_{ss}\|_{L^2} &\le C\abs{\log\epsilon}\abs{p_0} + C\epsilon\big( \epsilon^{-1/2}\abs{p_0} + \epsilon^{-1/2}\|a^{3/2}(a^4v_s)_{ss}\|_{L^2} \big)\\
&\le C\abs{\log\epsilon}\abs{p_0} + C\epsilon^{1/2}\|a^{3/2}(a^4v_s)_{ss}\|_{L^2}\,.
\end{aligned}
\end{equation}
Again, for $\epsilon$ sufficiently small, we may absorb the right-hand side $\|a^{3/2}(a^4v_s)_{ss}\|_{L^2}$ into the left-hand side to obtain 
\begin{equation}
\|a^{3/2}(a^4v_s)_{ss}\|_{L^2} \le C\abs{\log\epsilon}\abs{p_0}\,.
\end{equation}
Using this $L^2$ bound in \eqref{eq:vsss_Linfty}, we thus have
\begin{equation}\label{eq:vsss_Linfty_2}
\norm{a\,(a^4v_s)_{ss}}_{L^\infty} \le C\epsilon^{-1/2}\abs{\log\epsilon}\abs{p_0} \,.
\end{equation} 

Finally, recalling that $v = p^{\rm SB} - p_0$, we obtain Theorem \ref{thm:pSB}.
\hfill \qedsymbol

\section{Residual error for 1D approximation}\label{sec:residual}
Given the estimates of Theorem \ref{thm:pSB} for $p^{\rm SB}(s)$, we now turn to the approximation $q^{\rm SB}(\bx)$. We may use the formula \eqref{eq:qSB} for $q^{\rm SB}$ along with Lemmas \ref{lem:int1}-\ref{lem:at_bound} to calculate the following residual for the boundary value $\frac{\p q^{\rm SB}}{\p \bm{n}}$ along $\Gamma_\epsilon$. 

\begin{lemma}[Residual error for $q^{\rm SB}$ normal derivative]\label{lem:qSB_residual}
Along the blood vessel surface $\Gamma_\epsilon$, the approximation $q^{\rm SB}(\bx)$ given by \eqref{eq:qSB} satisfies 
\begin{equation}\label{eq:err1}
\frac{\p q^{\rm SB}}{\p \bm{n}} = \frac{\omega}{\epsilon}\big(p^{\rm SB} - q^{\rm SB} \big) + \mc{R}_\epsilon(s,\theta)\,,
\end{equation}
where the residual $\mc{R}_\epsilon$ satisfies 
\begin{equation}\label{eq:residual1}
\abs{\mc{R}_\epsilon} \le C\min\{a^{-1}(s),\epsilon^{-1}\}\big(\norm{a\,\big(a^4p^{\rm SB}_s\big)_s}_{C^1}+ \abs{\log\epsilon}\norm{a^{-1}\big(a^4p^{\rm SB}_s\big)_s}_{L^\infty} \big)\,.
\end{equation}
Furthermore, we have
\begin{equation}\label{eq:err2}
\int_0^{2\pi}\frac{\p q^{\rm SB}}{\p \bm{n}}\,\mc{J}_\epsilon(s,\theta)\,d\theta = \eta\, \big(a^4p^{\rm SB}_s\big)_s + \overline{\mc{R}_\epsilon}(s)\,,
\end{equation}
where 
\begin{equation}\label{eq:residual2}
\abs{\overline{\mc{R}_\epsilon}(s)}\le C\epsilon \big(\norm{a\,\big(a^4p^{\rm SB}_s\big)_s}_{C^1}+ \abs{\log\epsilon}\norm{a^{-1}\big(a^4p^{\rm SB}_s\big)_s}_{L^\infty} \big)\,.
\end{equation}
\end{lemma}

\begin{proof}
Using the equation \eqref{eq:pSB}, we may use that
\begin{equation}\label{eq:pSB_form}
\omega \,p^{\rm SB}(s) = \frac{1}{2\pi a(s)}\bigg(\eta \,\big(a^4p^{\rm SB}_s\big)_s+ \omega\int_0^{2\pi}q^{\rm SB}\big|_{\Gamma_\epsilon}\,a(s)\,d\theta \bigg)
\end{equation}
to instead show that  
\begin{equation}\label{eq:to_show}
\frac{\p q^{\rm SB}}{\p \bm{n}} = \frac{\eta}{2\pi\epsilon a(s)} \big(a^4p^{\rm SB}_s\big)_s+ \frac{\omega}{\epsilon}\bigg(\frac{1}{2\pi}\int_0^{2\pi}q^{\rm SB}\big|_{\Gamma_\epsilon}\,d\theta - q^{\rm SB}\big|_{\Gamma_\epsilon} \bigg) + \mc{R}_\epsilon
\end{equation}
for $\mc{R}_\epsilon$ as in \eqref{eq:residual1}. This will imply both \eqref{eq:err1} and \eqref{eq:err2}.

Recall that, using the notation of equation \eqref{eq:SN_Gamma}, along $\Gamma_\epsilon$, we may write 
\begin{equation}
q^{\rm SB}\big|_{\Gamma_\epsilon} = \frac{\eta}{4\pi(1-\ell(\epsilon))}\int_{s-\varphi_\epsilon^{-1}(1)}^{s+\varphi_\epsilon^{-1}(1)}\frac{1}{\abs{\bR_\epsilon}}\,\big[\big(a^4p^{\rm SB}_s\big)_s\big]^*(\varphi_\epsilon(s-\bars))\,d\bars\,,
\end{equation}
where $\big[\big(a^4p^{\rm SB}_s\big)_s\big]^*$ denotes the extension-by-reflection of $\big(a^4p^{\rm SB}_s\big)_s$ about $s=0$, as in \eqref{eq:fstar}. 
We will use the above expression for $q^{\rm SB}\big|_{\Gamma_\epsilon}$ to show \eqref{eq:to_show}.

We first calculate $\frac{\p q^{\rm SB}}{\p\bm{n}}$ along $\Gamma_\epsilon$. Since the radius $\epsilon a(s)$ varies along the length of the blood vessel, we have 
\begin{equation}
\bm{n}(s,\theta) = -\frac{1}{\sqrt{1+\epsilon^2(a')^2}}\be_r(s,\theta) + \frac{\epsilon a'(s)}{\sqrt{1+\epsilon^2(a')^2}}\be_{\rm t}(s)\,.
\end{equation}
In particular, at $s=1$, the normal vector to $\Gamma_\epsilon$ points purely in the $\be_{\rm t}$ direction, while away from $s=1$, the normal vector is roughly just the radial vector $\be_r$. Then, with respect to the curved cylindrical coordinates $(r,\theta,s)$, we may write 
\begin{equation}\label{eq:qSBn}
\frac{\p q^{\rm SB}}{\p\bm{n}} = \frac{1}{\sqrt{1+\epsilon^2(a')^2}}\bigg(-\frac{\p q^{\rm SB}}{\p r}\bigg|_{r=\epsilon a(s)} + \frac{\epsilon a'(s)}{1-\epsilon a\wh\kappa}\frac{\p q^{\rm SB}}{\p s}\bigg|_{r=\epsilon a(s)} \bigg)\,.
\end{equation}

Using the derivatives \eqref{eq:Rderivs} of the expression $\bR$, we may calculate 
\begin{equation}\label{eq:qSBr}
\begin{aligned}
-\frac{\p q^{\rm SB}}{\p r}\bigg|_{r=\epsilon a(s)} &= \frac{\eta}{4\pi(1-\ell(\epsilon))}\int_{s-\varphi_\epsilon^{-1}(1)}^{s+\varphi_\epsilon^{-1}(1)} \frac{\bR_\epsilon\cdot\be_r(s,\theta)}{\abs{\bR_\epsilon}^3}\big[\big(a^4p^{\rm SB}_s\big)_s\big]^*(\varphi_\epsilon(s-\bars))\,d\bars\\
&= \frac{\eta}{4\pi(1-\ell(\epsilon))}\int_{s-\varphi_\epsilon^{-1}(1)}^{s+\varphi_\epsilon^{-1}(1)} \frac{\bars^2\bm{Q}\cdot\be_r + \epsilon a}{\abs{\bR_\epsilon}^3}\big[\big(a^4p^{\rm SB}_s\big)_s\big]^*(\varphi_\epsilon(s-\bars))\,d\bars\,,
\end{aligned}
\end{equation}
as well as 
\begin{equation}\label{eq:qSBs}
\begin{aligned}
\frac{\p q^{\rm SB}}{\p s}\bigg|_{r=\epsilon a(s)} &=
-\frac{\eta}{4\pi(1-\ell(\epsilon))}\int_{s-\varphi_\epsilon^{-1}(1)}^{s+\varphi_\epsilon^{-1}(1)}\frac{\bR_\epsilon\cdot\be_{\rm t}(1-\epsilon a\wh\kappa)}{\abs{\bR_\epsilon}^3}\big[\big(a^4p^{\rm SB}_s\big)_s\big]^*(\varphi_\epsilon(s-\bars))\,d\bars\\
&=-\frac{\eta}{4\pi(1-\ell(\epsilon))}\int_{s-\varphi_\epsilon^{-1}(1)}^{s+\varphi_\epsilon^{-1}(1)} \frac{(\bars+\bars^2\bm{Q}\cdot\be_{\rm t})(1-\epsilon a\wh\kappa)}{\abs{\bR_\epsilon}^3}\big[\big(a^4p^{\rm SB}_s\big)_s\big]^*(\varphi_\epsilon(s-\bars))\,d\bars\,.
\end{aligned}
\end{equation}
Using Lemma \ref{lem:int1}, we have
\begin{equation}
\begin{aligned}
-\frac{\p q^{\rm SB}}{\p r}\bigg|_{r=\epsilon a(s)} &= \frac{\eta}{4\pi}\int_{s-\varphi_\epsilon^{-1}(1)}^{s+\varphi_\epsilon^{-1}(1)} \frac{\epsilon a(s)}{\abs{\bR_\epsilon}^3}\big[\big(a^4p^{\rm SB}_s\big)_s\big]^*(\varphi_\epsilon(s-\bars))\,d\bars + \mc{R}^{(1)}(s,\theta)\,, \\
\frac{\p q^{\rm SB}}{\p s}\bigg|_{r=\epsilon a(s)} &= -\frac{\eta}{4\pi}\int_{s-\varphi_\epsilon^{-1}(1)}^{s+\varphi_\epsilon^{-1}(1)} \frac{\bars}{\abs{\bR_\epsilon}^3}\big[\big(a^4p^{\rm SB}_s\big)_s\big]^*(\varphi_\epsilon(s-\bars))\,d\bars + \mc{R}^{(2)}(s,\theta)\,,
\end{aligned}
\end{equation}
where, for both $j=1,2$,
\begin{equation}
|\mc{R}^{(j)}(s,\theta)| \le C\abs{\log\epsilon}\norm{\big(a^4p^{\rm SB}_s\big)_s}_{L^\infty}\,.
\end{equation}
%
Recall that for $s\in [0,1]$ and $t\in[-1,1]$, we may expand 
\begin{equation}
a(s) = a^*(t) + (s-t)\mc{R}^{(a)}(s,t) 
\end{equation}
for $\mc{R}^{(a)}(s,t)$ as in \eqref{eq:Ra_bd}. Replacing $t$ with $\varphi_\epsilon(s-\bars)$ and using that $s-\varphi_\epsilon(s-\bars)=\frac{\bars-\ell(\epsilon)s}{1-\ell(\epsilon)}$, we may then write 
\begin{equation}
\begin{aligned}
-\frac{\p q^{\rm SB}}{\p r}\bigg|_{r=\epsilon a(s)} &= \frac{\eta}{4\pi}\int_{s-\varphi_\epsilon^{-1}(1)}^{s+\varphi_\epsilon^{-1}(1)} \frac{\epsilon}{\abs{\bR_\epsilon}^3}\big[a\,\big(a^4p^{\rm SB}_s\big)_s\big]^*(\varphi_\epsilon(s-\bars))\,d\bars + \wt{\mc{R}}^{(1)}(s,\theta)\,, \\
\frac{\p q^{\rm SB}}{\p s}\bigg|_{r=\epsilon a(s)} &= -\frac{\eta a(s)^{-1}}{4\pi}\int_{s-\varphi_\epsilon^{-1}(1)}^{s+\varphi_\epsilon^{-1}(1)} \frac{\bars}{\abs{\bR_\epsilon}^3}\big[a\,\big(a^4p^{\rm SB}_s\big)_s\big]^*(\varphi_\epsilon(s-\bars))\,d\bars + \wt{\mc{R}}^{(2)}(s,\theta)\,,
\end{aligned}
\end{equation}
where, by Lemma \ref{lem:int1}, for $j=1,2$, 
\begin{equation}
|\wt{\mc{R}}^{(j)}(s,\theta)| \le C \min\{a(s)^{-1},\epsilon^{-1}\}\abs{\log\epsilon} \norm{a^{-1}\big(a^4p^{\rm SB}_s\big)_s}_{L^\infty}\,.
\end{equation}
Using Lemmas \ref{lem:int2} and \ref{lem:int3}, along $\Gamma_\epsilon$, we then have 
\begin{equation}
\begin{aligned}
&\abs{-\frac{\p q^{\rm SB}}{\p r}-\frac{\eta}{2\pi\epsilon a(s)}\big(a^4p^{\rm SB}_s\big)_s(s)} \\
&\quad \le C\min\{a(s)^{-1},\epsilon^{-1}\}\big(\norm{a\,\big(a^4p^{\rm SB}_s\big)_s}_{C^1} + \norm{a^{-1}\big(a^4p^{\rm SB}_s\big)_s}_{L^\infty}\big) + C|\wt{\mc{R}}^{(1)}| \,,
\end{aligned}
\end{equation}
while 
\begin{equation}
\abs{\frac{\p q^{\rm SB}}{\p s}} \le C\min\{a(s)^{-1},\epsilon^{-1}\}\big(\norm{a\,\big(a^4p^{\rm SB}_s\big)_s}_{C^1}+ \norm{a^{-1}\big(a^4p^{\rm SB}_s\big)_s}_{L^\infty}\big)+ C|\wt{\mc{R}}^{(2)}| \,.
\end{equation}
In total, using the expression \eqref{eq:qSBn} for $\frac{\p q^{\rm SB}}{\p\bm{n}}$, for $\epsilon$ sufficiently small, we have 
\begin{equation}\label{eq:qSBn_est}
\begin{aligned}
&\abs{\frac{\p q^{\rm SB}}{\p\bm{n}}-\frac{\eta}{2\pi\epsilon a(s)}\big(a^4p^{\rm SB}_s\big)_s(s)} \\
&\quad \le C\min\{a(s)^{-1},\epsilon^{-1}\}\big(\norm{a\,\big(a^4p^{\rm SB}_s\big)_s}_{C^1}+ \abs{\log\epsilon}\norm{a^{-1}\big(a^4p^{\rm SB}_s\big)_s}_{L^\infty} \big)\,.
\end{aligned}
\end{equation}

We next consider the expression
\begin{equation}
q^{\rm r}(s,\theta) =  \frac{\omega}{\epsilon}\bigg(\frac{1}{2\pi}\int_0^{2\pi}q^{\rm SB}\big|_{\Gamma_\epsilon}\,d\theta - q^{\rm SB}\big|_{\Gamma_\epsilon} \bigg) \,,
\end{equation}
which measures the $\theta$-dependence of $q^{\rm SB}$ on $\Gamma_\epsilon$. We may calculate
\begin{equation}
\begin{aligned}
\frac{\p q^{\rm SB}}{\p\theta}\bigg|_{r=\epsilon a(s)} &=
-\frac{\eta}{4\pi(1-\ell(\epsilon))}\int_{s-\varphi_\epsilon^{-1}(1)}^{s+\varphi_\epsilon^{-1}(1)} \frac{\epsilon a(s)\bR_\epsilon\cdot\be_\theta(s,\theta)}{\abs{\bR_\epsilon}^3}\big[\big(a^4p^{\rm SB}_s\big)_s\big]^*(\varphi_\epsilon(s-\bars))\,d\bars \\
&=-\frac{\epsilon a(s)\eta}{4\pi(1-\ell(\epsilon))}\int_{s-\varphi_\epsilon^{-1}(1)}^{s+\varphi_\epsilon^{-1}(1)}\frac{\bars^2\bm{Q}\cdot\be_\theta}{\abs{\bR_\epsilon}^3}\big[\big(a^4p^{\rm SB}_s\big)_s\big]^*(\varphi_\epsilon(s-\bars))\,d\bars\,.
\end{aligned}
\end{equation}
By Lemma \ref{lem:int1}, we have
\begin{equation}
\abs{\frac{\p q^{\rm SB}}{\p\theta}} \le C \epsilon a(s)\abs{\log\epsilon}\norm{\big(a^4p^{\rm SB}_s\big)_s}_{L^\infty}\,,
\end{equation}
and thus
\begin{equation}\label{eq:qr_est}
\abs{q^{\rm r}(s,\theta)}\le C a(s)\abs{\log\epsilon}\norm{\big(a^4p^{\rm SB}_s\big)_s}_{L^\infty}\,.
\end{equation}
Combining \eqref{eq:qSBn_est} and \eqref{eq:qr_est}, we have that 
\begin{equation}
\begin{aligned}
&\abs{\frac{\p q^{\rm SB}}{\p\bm{n}}-\frac{\eta}{2\pi\epsilon a(s)}\big(a^4p^{\rm SB}_s\big)_s(s) - q^{\rm r}(s,\theta)} \\
&\quad \le C\min\{a(s)^{-1},\epsilon^{-1}\}\big(\norm{a\,\big(a^4p^{\rm SB}_s\big)_s}_{C^1}+ \abs{\log\epsilon}\norm{a^{-1}\big(a^4p^{\rm SB}_s\big)_s}_{L^\infty} \big)\,.
\end{aligned}
\end{equation}
By the representation \eqref{eq:pSB_form} of $p^{\rm SB}$, we obtain \eqref{eq:err1}.

To show \eqref{eq:err2}, we multiply the expression \eqref{eq:to_show} by $\mc{J}_\epsilon$ and integrate in $\theta$ to obtain 
\begin{equation}
\begin{aligned}
\int_0^{2\pi}\frac{\p q^{\rm SB}}{\p\bm{n}}\,\mc{J}_\epsilon(s,\theta)\,d\theta
&= \eta \,\big(a^4p^{\rm SB}_s\big)_s + \epsilon a(s) \int_0^{2\pi}\mc{R}_\epsilon\,d\theta 
+ \int_0^{2\pi}\frac{\p q^{\rm SB}}{\p\bm{n}}\,\big(\mc{J}_\epsilon(s,\theta)- \epsilon a(s)\big)\,d\theta \,,
\end{aligned}
\end{equation}
where we have used that $q^{\rm r}$ integrates to zero in $\theta$. Using \eqref{eq:jac_est} along with \eqref{eq:qSBn}, \eqref{eq:qSBr}, and \eqref{eq:qSBs}, we have 
\begin{equation}
\begin{aligned}
&\abs{\int_0^{2\pi}\frac{\p q^{\rm SB}}{\p\bm{n}}\,\big(\mc{J}_\epsilon(s,\theta)- \epsilon a(s)\big)\,d\theta}\\
&\le C\epsilon^2\int_{s-\varphi_\epsilon^{-1}(1)}^{s+\varphi_\epsilon^{-1}(1)}\frac{a^*(\varphi_\epsilon(s-\bars))}{\abs{\bR_\epsilon}^2}\bigg|\big[a^{-1}\big(a^4p^{\rm SB}_s\big)_s\big]^*(\varphi_\epsilon(s-\bars))\bigg|\,d\bars
\le C\epsilon\norm{a^{-1}\big(a^4p^{\rm SB}_s\big)_s}_{L^\infty}\,,
\end{aligned}
\end{equation}
where we have used Lemma \ref{lem:at_bound} to obtain a bound without a prefactor of $a^{-1}(s)$. Combining the above bound with the bound \eqref{eq:err1} for $\mc{R}_\epsilon$, we obtain
\begin{equation}
\begin{aligned}
\abs{\int_0^{2\pi}\frac{\p q^{\rm SB}}{\p\bm{n}}\,\mc{J}_\epsilon(s,\theta)\,d\theta - \eta \,\big(a^4p^{\rm SB}_s\big)_s} 
\le C\epsilon \big(\norm{a\,\big(a^4p^{\rm SB}_s\big)_s}_{C^1}+ \abs{\log\epsilon}\norm{a^{-1}\big(a^4p^{\rm SB}_s\big)_s}_{L^\infty} \big)\,.
\end{aligned}
\end{equation}
\end{proof}

\section{Error estimate: 3D-1D to 1D}\label{sec:error}
Equipped with the residual bounds of Lemma \ref{lem:qSB_residual} for $q^{\rm SB}(\bx)$ as well as the estimates of Theorem \ref{thm:pSB} for $p^{\rm SB}(s)$, we may proceed to the proof of Theorem \ref{thm:3D1Dto1D}.

\begin{proof}[Proof of Theorem \ref{thm:3D1Dto1D}]
We consider the difference $(\wh q,\wh p)=(q^{\rm SB}-q,p^{\rm SB}-p)$ where $(q,p)$ satisfies \eqref{eq:darcy}-\eqref{eq:extBCs} and $(q^{\rm SB},p^{\rm SB})$ is given by \eqref{eq:qSB}, \eqref{eq:pSB}. By Lemma \ref{lem:qSB_residual}, we have that $(\wh q,\wh p)$ satisfies the PDE 
\begin{subequations}
\begin{align}
\Delta \wh q &= 0 \hspace{3.5cm} \text{in }\Omega_\epsilon  \label{eq:errDarcy} \\
\frac{\p \wh q}{\p \bm{n}} &= \frac{\omega}{\epsilon} (\wh p(s) -\wh q)+ \mc{R}_\epsilon \qquad \text{on }\Gamma_\epsilon  \label{eq:errBC1} \\
\int_0^{2\pi}\frac{\p \wh q}{\p \bm{n}}\,\mc{J}_\epsilon(s,\theta)\,d\theta &= \eta (a^4\wh p_s)_s + \overline{\mc{R}_\epsilon}(s)  \qquad \text{ on }\Gamma_\epsilon  \label{eq:errBC2} \\
\frac{\p \wh q}{\p\bm{n}} &= 0 \; \text{ on }\p\Omega_\epsilon\backslash\Gamma_\epsilon\,, \quad \wh q\to 0 \text{ as }\abs{\bx}\to\infty\,,
\end{align}
\end{subequations}
with $\wh p(0)=0$. We recall that the residuals satisfy
\begin{equation}\label{eq:recall}
\begin{aligned}
\abs{\mc{R}_\epsilon(s,\theta)} &\le C\min\{a^{-1}(s),\epsilon^{-1}\}\big(\norm{a\,\big(a^4p^{\rm SB}_s\big)_s}_{C^1}+ \abs{\log\epsilon}\norm{a^{-1}\big(a^4p^{\rm SB}_s\big)_s}_{L^\infty} \big)\,,\\
\abs{\overline{\mc{R}_\epsilon}(s)}&\le C\epsilon \big(\norm{a\,\big(a^4p^{\rm SB}_s\big)_s}_{C^1}+ \abs{\log\epsilon}\norm{a^{-1}\big(a^4p^{\rm SB}_s\big)_s}_{L^\infty} \big)\,.
\end{aligned}
\end{equation}

Now, multiplying \eqref{eq:errDarcy} by $\wh q(\bx)$ and integrating by parts, we have that $(\wh q,\wh p)$ satisfies the identity
\begin{equation}\label{eq:whq_ID}
\begin{aligned}
&\int_{\Omega_\epsilon}\abs{\nabla \wh q}^2\,d\bx = \int_{\Gamma_\epsilon}\wh q\frac{\p\wh q}{\p\bm{n}}\,dS_\epsilon 
= \int_{\Gamma_\epsilon}(\wh q-\wh p)\frac{\p\wh q}{\p\bm{n}}\,dS_\epsilon + \int_0^1\wh p(s)\int_0^{2\pi}\frac{\p\wh q}{\p\bm{n}}\,\mc{J}_\epsilon\,d\theta ds \\
&= -\frac{\omega}{\epsilon}\int_{\Gamma_\epsilon}(\wh p-\wh q)^2\,dS_\epsilon +\eta\int_0^1\wh p(a^4\wh p_s)_s\,ds + \int_{\Gamma_\epsilon}(\wh q-\wh p)\mc{R}_\epsilon\,dS_\epsilon + \int_0^1\wh p(s)\overline{\mc{R}_\epsilon}(s)\,ds\\
&= -\frac{\omega}{\epsilon}\int_{\Gamma_\epsilon}(\wh p-\wh q)^2\,dS_\epsilon -\eta\int_0^1a^4\wh p_s^2(s)\,ds 
+ \int_{\Gamma_\epsilon}(\wh q-\wh p)\mc{R}_\epsilon\,dS_\epsilon + \int_0^1\wh p(s)\overline{\mc{R}_\epsilon}(s)\,ds\,.
\end{aligned}
\end{equation}
Using that $\wh p\in \mc{H}^a_0(0,1)$ satisfies the weighted Poincar\'e inequality \eqref{eq:Ha_poincare}, we may estimate 
\begin{equation}
\begin{aligned}
&\int_{\Omega_\epsilon}\abs{\nabla \wh q}^2\,d\bx + \frac{\omega}{\epsilon}\int_{\Gamma_\epsilon}\abs{\wh p-\wh q}^2\,dS_\epsilon +\eta\int_0^1a^4\wh p_s^2(s)\,ds \\
&\quad \le 
C\bigg(\int_{\Gamma_\epsilon}\abs{\wh q-\wh p}^2dS_\epsilon\bigg)^{1/2}\bigg(\int_{\Gamma_\epsilon}\abs{\mc{R}_\epsilon}^2dS_\epsilon\bigg)^{1/2} 
 + C\bigg(\int_0^1\abs{\wh p}^2ds\bigg)^{1/2}\bigg(\int_0^1\abs{\overline{\mc{R}_\epsilon}(s)}^2ds\bigg)^{1/2}\\
&\quad \le \frac{\omega}{2\epsilon}\int_{\Gamma_\epsilon}\abs{\wh q-\wh p}^2dS_\epsilon + C(\omega)\epsilon\int_{\Gamma_\epsilon}\abs{\mc{R}_\epsilon}^2dS_\epsilon
+ \frac{\eta}{2}\int_0^1a^4\wh p_s^2\,ds + C(\eta)\int_0^1\abs{\overline{\mc{R}_\epsilon}(s)}^2ds \,,
\end{aligned}
\end{equation}
where we have used Young's inequality to split terms. Absorbing the first and third terms into the left hand side, we have 
\begin{equation}
\begin{aligned}
&\int_{\Omega_\epsilon}\abs{\nabla \wh q}^2\,d\bx + \frac{\omega}{2\epsilon}\int_{\Gamma_\epsilon}(\wh p-\wh q)^2\,dS_\epsilon +\frac{\eta}{2}\int_0^1a^4\wh p_s^2(s)\,ds \\
&\quad \le C\int_0^1\abs{\overline{\mc{R}_\epsilon}(s)}^2ds + C\epsilon\int_0^1\int_0^{2\pi}\abs{\mc{R}_\epsilon}^2 \epsilon a(s) d\theta\,ds + C\epsilon\int_0^1\int_0^{2\pi}\abs{\mc{R}_\epsilon}^2 \epsilon^2 d\theta\,ds \\
&\quad \le C\bigg(\epsilon^2+\epsilon^2\int_0^1a^{-1}(s)\,ds\bigg) \big(\norm{a\,\big(a^4p^{\rm SB}_s\big)_s}_{C^1}^2+ \abs{\log\epsilon}^2\norm{a^{-1}\big(a^4p^{\rm SB}_s\big)_s}_{L^\infty}^2 \big) \,,
\end{aligned}
\end{equation}
where we have used \eqref{eq:jac_est} to split up $dS_\epsilon = \mc{J}_\epsilon(s,\theta)\,d\theta ds$ in the second line and have used the bounds \eqref{eq:recall} for $\overline{\mc{R}_\epsilon}$ and $\mc{R}_\epsilon$ in the third line. Note that by the spheroidal endpoint condition \eqref{eq:spheroidal}, $a^{-1}(s)$ is integrable over $[0,1]$, and the value is independent of $\epsilon$.

Using Theorem \ref{thm:pSB}, we thus have 
\begin{equation}
\begin{aligned}
&\norm{\nabla \wh q}_{L^2(\Omega_\epsilon)}+ \epsilon^{-1/2}\norm{\wh p -\wh q}_{L^2(\Gamma_\epsilon)}+ \norm{a^2\wh p_s}_{L^2(0,1)} \\
&\qquad \le C\epsilon\big(\norm{a\,\big(a^4p^{\rm SB}_s\big)_s}_{C^1}+ \abs{\log\epsilon}\norm{a^{-1}\big(a^4p^{\rm SB}_s\big)_s}_{L^\infty} \big)\\
&\qquad \le C \epsilon^{1/2}\abs{\log\epsilon}\,\abs{p_0}\,.
\end{aligned}
\end{equation}
By the Sobolev inequality \eqref{eq:sob} on $\Omega_\epsilon$ and the weighted Poincar\'e inequality \eqref{eq:Ha_poincare}, we obtain Theorem \ref{thm:3D1Dto1D}.
\end{proof}

\appendix
\section{Well-posedness of 3D-1D system}\label{sec:app}
Here we briefly discuss the well-posedness theory for the 3D-1D system \eqref{eq:darcy}-\eqref{eq:extBCs} stated in Theorem \ref{thm:3d1d}, which is also explored in Part II \cite{partII} of our program. 

Since the domain $\Omega_\epsilon=\R^3_+\backslash\overline{\mc{V}_\epsilon}$ is unbounded, and we expect $q(\bx)$ satisfying \eqref{eq:darcy} to decay like $\mc{G}_N(\bx,0)\sim\frac{1}{\abs{\bx}}$ as $\abs{\bx}\to\infty$, we do not seek $q\in L^2(\Omega_\epsilon)$, but we do expect $q\in L^6(\Omega_\epsilon)$. Following \cite[Chap. II.6-II.10]{galdi2011introduction}, we will work with functions belonging to the homogeneous Sobolev space 
\begin{equation}
D^{1,2}(\Omega_\epsilon) = L^6(\Omega_\epsilon)\cap\dot H^1(\Omega_\epsilon) = \{ u\in L^6(\Omega_\epsilon)\,:\, \nabla u\in L^2(\Omega_\epsilon)\}\,.
\end{equation}
Note that by the Sobolev inequality
\begin{equation}\label{eq:sob}
\norm{u}_{L^6(\Omega_\epsilon)}\le C\norm{\nabla u}_{L^2(\Omega_\epsilon)}\,,
\end{equation}
where $C$ is independent of $\epsilon$ (see \cite{free_ends,closed_loop}), we have that $D^{1,2}(\Omega_\epsilon)$ is a Hilbert space with norm
\begin{equation}
\norm{u}_{D^{1,2}(\Omega_\epsilon)}:= \norm{\nabla u}_{L^2(\Omega_\epsilon)}\,.
\end{equation}

Turning to the system \eqref{eq:darcy}-\eqref{eq:extBCs}, we will consider $\overline p = p-p_0$ to lift the boundary condition for $p$ to the interior of the vessel. The condition \eqref{eq:robin} along $\Gamma_\epsilon$ then becomes
\begin{align}
\frac{\p q}{\p\bm{n}} &= \frac{\omega}{\epsilon}\big(\overline p(s)-q) + \frac{\omega}{\epsilon}p_0\,. \label{eq:BC1do}
%
\end{align}

Recalling the definition \eqref{eq:Ha_def} of the weighted space $\mc{H}^a(0,1)$ as well as the space $\mc{H}^a_0(0,1)= \{ u\in \mc{H}^a(0,1)\,:\, u\big|_{s=0}=0 \}$, we define the bilinear form $\mc{B}:\big( D^{1,2}(\Omega_\epsilon)\times \mc{H}^a_0(0,1)\big)\times \big( D^{1,2}(\Omega_\epsilon)\times \mc{H}^a_0(0,1)\big)\to \R$ as
\begin{equation}
\mc{B}\big((q,\overline p),(g,\gamma)\big) := \int_{\Omega_\epsilon}\nabla q\cdot\nabla g\,d\bx +\frac{\omega}{\epsilon}\int_{\Gamma_\epsilon}(\gamma(s)-g)(\overline p(s)- q)\,dS_\epsilon +\eta\int_0^1\overline p_s\gamma_s\,a^4(s)\,ds \,. 
\end{equation}
Note that by \cite[section 4]{partII}, we have that $\mc{B}(\cdot,\cdot)$ is bounded and coercive on $\big( D^{1,2}(\Omega_\epsilon)\times \mc{H}^a_0(0,1)\big)\times \big( D^{1,2}(\Omega_\epsilon)\times \mc{H}^a_0(0,1)\big)$. 
We may then define a weak solution to the 3D-1D system as $(q,\overline p)\in D^{1,2}(\Omega_\epsilon)\times \mc{H}^a_0(0,1)$ satisfying
\begin{equation}\label{eq:weak}
\mc{B}\big((q,\overline p),(g,\gamma)\big)= \frac{\omega}{\epsilon}\int_{\Gamma_\epsilon}(g-\gamma(s))\,p_0\,dS_\epsilon
\end{equation}
for all $(g,\gamma)\in D^{1,2}(\Omega_\epsilon)\times \mc{H}^a_0(0,1)$.
The weak formulation \eqref{eq:weak} may be formally derived from \eqref{eq:darcy}, \eqref{eq:extBCs}, and \eqref{eq:BC1do} as follows. Given $(g,\gamma)\in D^{1,2}(\Omega_\epsilon)\times \mc{H}^a_0(0,1)$, we may multiply \eqref{eq:darcy} by $g$ and integrate by parts to obtain
\begin{equation}
\begin{aligned}
\int_{\Omega_\epsilon}&\nabla q\cdot\nabla g\,d\bx =\int_{\Gamma_\epsilon}g\,\frac{\p q}{\p\bm{n}}\,dS_\epsilon 
= \int_{\Gamma_\epsilon}(g-\gamma(s))\,\frac{\p q}{\p\bm{n}}\,dS_\epsilon + \int_0^1\gamma(s)\int_0^{2\pi}\frac{\p q}{\p\bm{n}}\,\mc{J}_\epsilon\,d\theta ds\\
&= \frac{\omega}{\epsilon}\int_{\Gamma_\epsilon}(g-\gamma(s))(\overline p(s)- q)\,dS_\epsilon + \frac{\omega}{\epsilon}\int_{\Gamma_\epsilon}(g-\gamma(s))\,p_0\,dS_\epsilon  + \eta\int_0^1\gamma(s)\p_s\big(a^4(s)\overline p_s\big)\, ds\\
&= -\frac{\omega}{\epsilon}\int_{\Gamma_\epsilon}(\gamma(s)-g)(\overline p(s)- q)\,dS_\epsilon - \eta\int_0^1\gamma_s(s)\overline p_s(s) a^4(s)\, ds  + \frac{\omega}{\epsilon}\int_{\Gamma_\epsilon}(g-\gamma(s))\,p_0\,dS_\epsilon  \,.
\end{aligned}
\end{equation}
 Here, in the first line, we have added and subtracted $\gamma(s)$; in the second line, we have used the boundary conditions \eqref{eq:BC1do} and \eqref{eq:RtoN}; and, in the final equality, we have integrated by parts in $s$. 

Existence and uniqueness of a weak solution \eqref{eq:weak} follows from the Lax-Milgram lemma. In addition, taking $(g,\gamma)=(q,\overline p)$, we obtain the energy estimate 
\begin{equation}
\begin{aligned}
&\int_{\Omega_\epsilon}\abs{\nabla q}^2\,d\bx +\frac{\omega}{\epsilon}\int_{\Gamma_\epsilon}\abs{\overline p(s)-q}^2\,dS_\epsilon +\eta\int_0^1a^4(s)\overline p_s^2\,ds 
\le \frac{\omega}{\epsilon}|p_0|\abs{\Gamma_\epsilon}^{1/2}\norm{q-\overline p}_{L^2(\Gamma_\epsilon)}\,. 
\end{aligned}
\end{equation}
Using Young's inequality and that $\abs{\Gamma_\epsilon}^{1/2}\le C\epsilon^{1/2}$, we may split up and hide terms from the right-hand side in the left-hand side to obtain the estimate 
\begin{equation}
\int_{\Omega_\epsilon}\abs{\nabla q}^2\,d\bx +\frac{\omega}{\epsilon}\int_{\Gamma_\epsilon}\abs{\overline p(s)- q}^2\,dS_\epsilon +\eta\int_0^1a^4(s)\overline p_s^2\,ds \le C\abs{p_0}^2 \,.
\end{equation}
Using that $\overline p\in \mc{H}^a_0$ satisfies the weighted Poincar\'e inequality \eqref{eq:Ha_poincare} and recalling that $\overline p = p-p_0$, we obtain the estimate \eqref{eq:energy} of Theorem \ref{thm:3d1d}. 
\hfill \qedsymbol



\subsubsection*{Acknowledgments} We thank Yoichiro Mori for guidance in formulating the 3D-3D and 3D-1D models and for helpful initial discussions. LO acknowledges support from NSF grant DMS-2406003 and from the VCRGE with funding from the Wisconsin Alumni Research Foundation. SS acknowledges support from NSF grant DMS-2037851.

\bibliographystyle{abbrv} 
\bibliography{BloodBib}

\end{document}